\documentclass[11pt]{amsart}
\usepackage{amssymb,latexsym,amscd,array}
\usepackage[leqno]{amsmath}
\usepackage{tikz}
\usepackage{tikz-cd}
\usepackage{enumerate}
\usepackage{exscale}
\usepackage[centertags]{amsmath}
\usepackage{amsthm}
\usepackage{dsfont}
\usepackage[all]{xy}
\usepackage{lscape}
\usepackage{pdflscape}
\usepackage{color}
\usepackage{graphicx}
\usepackage[colorlinks,plainpages,backref]{hyperref}
\usepackage{url}
\usepackage{lineno}
\usepackage[neveradjust]{paralist}
\usepackage{todonotes}
\usepackage{csquotes}
\usepackage{mathabx}
\usepackage{multirow}
\usepackage[mathscr]{euscript}
\usepackage{thmtools}
\usepackage{comment}

\newtheorem{thm}{Theorem}[section]

\declaretheorem[
  style=definition,
  title=Example,
  qed={$\diamond$},
  refname={example,examples},
  Refname={Example,Examples},
  sharenumber=thm,
]{exa}

\declaretheorem[
  style=definition,
  title=Definition,
  qed={$\diamond$},
  sharenumber=thm,
]{dfn}

\declaretheorem[
  style=definition,
  title=Convention,
  qed={$\diamond$},
  sharenumber=thm,
]{cnv}

\declaretheorem[
  style=definition,
  title=Notation,
  qed={$\diamond$},
  sharenumber=thm,
]{ntn}

\declaretheorem[
  style=definition,
  title=Remark,
  qed={$\diamond$},
  sharenumber=thm,
]{rmk}

\declaretheorem[
  title=Question,
  qed={$\diamond$},
  sharenumber=thm,
]{qun}

\declaretheorem[
  title=Conjecture,
  qed={$\diamond$},
  sharenumber=thm,
]{conj}

\declaretheorem[
  title=Problem,
  qed={$\diamond$},
  sharenumber=thm,
]{prb}

   \theoremstyle{definition}
   \theoremstyle{remark}

\newcommand{\CC}{{\mathbb C}}
\newcommand{\DD}{{\mathbb D}}
\newcommand{\FF}{{\mathbb F}}
\newcommand{\kk}{{\Bbbk}}
\newcommand{\KK}{{\mathbb K}}
\newcommand{\NN}{{\mathbb N}}
\newcommand{\PP}{{\mathbb P}}
\newcommand{\QQ}{{\mathbb Q}}
\newcommand{\RR}{{\mathbb R}}
\newcommand{\TT}{{\mathbb T}}
\newcommand{\VV}{{\mathbb V}}
\newcommand{\ZZ}{{\mathbb Z}}

\newcommand{\HH}{{\mathrm{HH}}}

\DeclareMathOperator{\ann}{ann}
\DeclareMathOperator{\ara}{ara}
\DeclareMathOperator{\CM}{CM}
\DeclareMathOperator{\codim}{codim}
\DeclareMathOperator{\depth}{depth}
\DeclareMathOperator{\dRdepth}{dR-depth}

\DeclareMathOperator{\dlim}{\varinjlim}

\DeclareMathOperator{\charac}{char}

\DeclareMathOperator{\ch}{ch}
\DeclareMathOperator{\Ch}{Ch}
\DeclareMathOperator{\cd}{cd}
\DeclareMathOperator{\ecd}{ecd}
\DeclareMathOperator{\embdim}{embdim}
\DeclareMathOperator{\Ext}{Ext}
\DeclareMathOperator{\FL}{FL}
\DeclareMathOperator{\FS}{FS}
\DeclareMathOperator{\funcF}{{\mathrm F}}
\DeclareMathOperator{\gr}{gr}

\DeclareMathOperator{\hight}{ht}
\DeclareMathOperator{\Hom}{Hom}
\DeclareMathOperator{\id}{id}
\DeclareMathOperator{\image}{im}
\DeclareMathOperator{\Jac}{Jac}

\DeclareMathOperator{\lcd}{lcd}
\DeclareMathOperator{\Lie}{Lie}

\DeclareMathOperator{\mdim}{mdim}
\DeclareMathOperator{\Pic}{Pic}
\DeclareMathOperator{\Proj}{Proj}

\DeclareMathOperator{\reg}{reg}
\DeclareMathOperator{\RGamma}{{\mathrm R}\Gamma}
\DeclareMathOperator{\rk}{rk}

\DeclareMathOperator{\Spec}{Spec}

\DeclareMathOperator{\Tor}{Tor}

\DeclareMathOperator{\Var}{Var}
\DeclareMathOperator{\vol}{vol}
\DeclareMathOperator{\Supp}{Supp}
\DeclareMathOperator{\injdim}{injdim}
\DeclareMathOperator{\End}{End}
\DeclareMathOperator{\Ass}{Ass}
\DeclareMathOperator{\rank}{rank}

\newcommand{\bsx}{{\boldsymbol x}}
\newcommand{\bsg}{{\boldsymbol g}}
\newcommand{\bsdel}{{\boldsymbol \del}}

\newcommand{\boldzero}{{\mathbf{0}}}

\newcommand{\bolda}{{\mathbf{a}}}

\newcommand{\boldu}{{\mathbf{u}}}

\newcommand{\del}{\partial}

\newcommand{\an}{{\mathrm{an}}}
\newcommand{\de}{{\mathrm d}}
\newcommand{\dR}{{\mathrm{dR}}}
\newcommand{\et}{{\mathrm{et}}}
\newcommand{\FH}{{\rm{FH}}}

\newcommand{\red}{{\mathrm{red}}}

\newcommand{\Sing}{{\mathrm{Sing}}}
\newcommand{\st}{{\mathrm{st}}}

\newcommand{\eps}{\varepsilon}

\newcommand{\Gr}{{\mathrm{Gr}}}
\newcommand{\ideal}[1]{{( #1 )}}
\newcommand{\into}{\hookrightarrow}

\newcommand{\onto}{\twoheadrightarrow}
\renewcommand{\to}{\longrightarrow}

\newcommand{\calC}{\mathcal{C}}
\newcommand{\calD}{\mathcal{D}}
\newcommand{\calE}{\mathcal{E}}
\newcommand{\calF}{\mathcal{F}}

\newcommand{\calH}{\mathcal{H}}

\newcommand{\calL}{\mathcal{L}}

\newcommand{\calO}{\mathcal{O}}

\newcommand{\fraka}{{\mathfrak{a}}}
\newcommand{\frakb}{{\mathfrak{b}}}
\newcommand{\frakg}{{\mathfrak{g}}}
\newcommand{\frakm}{{\mathfrak{m}}}
\newcommand{\frakp}{{\mathfrak{p}}}

\newcommand{\fraky}{{\mathfrak{y}}}
\newcommand{\frakP}{{\mathfrak{P}}}

\newcommand{\frakX}{{\mathfrak{X}}}

\newcommand{\pointt}{{\mathfrak{t}}}
\newcommand{\pointx}{{\mathfrak{x}}}
\newcommand{\pointy}{{\mathfrak{y}}}

\newcommand{\wz}[1]{\todo{WZ: \enquote{#1}}}

\numberwithin{equation}{subsubsection}

\begin{document}
\title[Local cohomology]{Local cohomology---an invitation}
\author[U.~Walther]{Uli Walther}
\address{U.~Walther : Purdue University, Dept.\ of Mathematics, 150 N.\ University St., West Lafayette, IN 47907, USA}
\email{walther@math.purdue.edu}

\author[W.~Zhang]{Wenliang Zhang}
\address{W.~Zhang: Department of Mathematics, Statistics, and Computer Science, University of Illinois, Chicago, IL 60607, USA}
\email{wlzhang@uic.edu}

\thanks{UW acknowledges support through Simons Foundation
  Collaboration Grant for Mathematicians~\#580839, and through NSF Grant
  DMS-2100288. WZ acknowledges the support by the NSF through
  DMS-1752081.}
\begin{abstract}
This article is part introduction and part survey to the 
mathematical area centered around local cohomology.
\end{abstract}


\maketitle

\setcounter{tocdepth}{3}

\tableofcontents

%

This article is a mixture of an introduction to local cohomology, and
a survey of the recent advances in the area, {with a view
towards} 
 relations to other parts of mathematics. It thus proceeds at times
rather carefully, with definitions and examples, and sometimes is more
cursory, aiming to give the reader an impression about certain parts
of the mathematical landscape. As such, it is more than a reference
list but less than a monograph. One possible use we envision is as a
guide for a novice, such as a beginning graduate student, to get an
idea what the general thrust of local cohomology is, and where one can
read more about certain topics.

While the article is rather much longer than originally anticipated,
several active areas that interact with local cohomology have been
left out. {For instance, we refer the reader to \cite{Hartshorne-lc-notes, Schenzel-Dualisierende,SchenzelUseofLocalCohomology} for connections with dualizing complexes which are not discussed in this article.} What we have put into the article is driven by personal
preferences and lack of expertise; we apologize to those
offended by our choices.

\smallskip

Over time, several excellent survey articles on local cohomology and
related themes have been written, and we strongly recommend the reader
study the following ones. One should name \cite{L-lc-survey}
on the state of the art 20 years since,
the article \cite{BetancourtWittZhang-survey}
specifically geared at Lyubeznik numbers, and the survey
\cite{Hochster-SurveyLC-CiA2020}.

In the more expository direction, we and many others have been
fortunate to be able to study Hochster's unpublished notes (available
on his website) and Huneke's point of view in
\cite{Huneke-lc-notes}. These notes come with our highest
recommendations and have strongly influenced us and this article. For
a treatment de-emphasizing Noetherianness we point at
 \cite{SchenzelSimon-NotNoetherian}. 

We close this thread of thoughts with mentioning the books concerned
with local cohomology as main subject: the original account of
Grothendieck as recorded by Hartshorne \cite{Hartshorne-lc-notes}, the
classic \cite{BrodmannSharp} by Brodmann and Sharp, and the outgrowth
\cite{24h} of a summer school on local cohomology.

\smallskip

Some words on the prerequisites for reading this article are in
order. Inasmuch as pure commutative algebra is concerned, we imagine
the reader be familiar with the contents of the book by Atiyah and
Macdonald \cite{A+M} or an appropriate subset of the book
by Eisenbud \cite{Eisenbud}. For homological
algebra one should know about injective and projective resolutions,
Ext and Tor and the principles of derived functors, and perhaps a bit
about spectral sequences at the level of Rotman
\cite{Rotman}. Hartshorne's opus \cite{Hartshorne-book} covers all
that is needed on varieties, schemes and sheaves in chapters 1-3. 

\section*{Acknowledgments}

We are grateful to Josep Alvarez, Robin Hartshorne, Jack Jeffries,
Peter Schenzel, Craig Huneke, Karl Schwede, Kazuma Shimomoto, and Matteo Varbaro for
telling us about corrections and suggestions in the manuscript.

Our main intellectual debt and gratitude is due to Gennady
Lyubeznik, our both advisor. We also happily acknowledge the
impact our many teachers and collaborators have had on our
understanding of the subjects discussed in this article. Special thanks
go to Mel Hochster and Anurag~K.~Singh for patience, insights, and
friendship.

\section{Introduction}
\begin{ntn}\label{ntn-setup}
  Throughout,  $A$ will denote a commutative
Noetherian ring. On occasion, $A$ will be
assumed to be local; then its maximal ideal is denoted by $\frakm$ and the
residue field by $\kk$.

We reserve the symbol $R$ for the case that $A$ is regular, while $M$
will generally denote a module over $A$.
\end{ntn}

\begin{dfn}\label{dfn-lc}
For an ideal $I\subseteq A$ the (left-exact) \emph{section functor
  with support in $I$} (also called the \emph{$I$-torsion functor})
$\Gamma_I(-)$ and the \emph{local cohomology functors $H^\bullet_I(-)$
  with support in $I$} are
  \[
    \Gamma_I\colon M\leadsto \{m\in M|\exists \ell\in\NN, I^\ell
    m=0\}
  \]
  and its right derived functors $H^\bullet_I(-)$.
Since $\Gamma_I(-)$ is left exact, $\Gamma_I(-)$ agrees with
$H^0_I(-)$. 
\end{dfn}

Local cohomology was invented by Grothendieck, at least in part, for
the purpose of proving Lefschetz and Barth type theorems (comparisons between a
smooth ambient variety and a possibly singular subvariety). The idea rests
on the fact, already exploited by Serre in \cite{Serre-FAC}, that
the geometry of projective varieties is encoded in the algebra of its
coordinate ring. Grothendieck makes it clear in his Harvard seminar
that, for this purpose, studying general properties of the concept of local
cohomological dimension is of great importance
\cite[p.~79]{Hartshorne-lc-notes}.

\begin{dfn}\label{dfn-lcd}
  The \emph{local cohomological dimension} $\lcd_A(I)$  of the
  $A$-ideal $I$ is
  \[
  \lcd_A(I)=\max\{k\in\NN\mid H^k_I(A)\neq 0\}.
  \]
One can show, using long exact sequences and direct limits, that
$H^{>\lcd_A(I)}_I(M)$ vanishes for every $A$-module $M$.
\end{dfn}

It is an essential feature of  the theory of local cohomology and
its applications that there are several different ways of {calculating}
$H^k_I(M)$ for any $A$-module $M$, all compatible with natural
functors. We review briefly three other approaches; for a more
complete account we refer to \cite{24h}.

\subsection{Koszul cohomology}
Let $x\in A$ be a single element and consider the multiplication map
$A\stackrel{x}{\to}A$
by $x$, also referred to as the \emph{cohomological Koszul complex}
$K^\bullet(A;x)$, so the displayed map is a morphism from position 0
to position 1 in the complex. We write $H^i(A;x)$ for the cohomology
modules of this complex.

Replacing $x$ with its own powers, one arrives at a tower of
commutative diagrams
\begin{eqnarray}\label{eqn-tower-x}
\begin{tikzcd}[ampersand replacement=\&,]
  A\ar[r,"x"]\ar[d,"1"]\&A\ar[d,"x"]\\
  A\ar[r,"x^2"]\ar[d,"1"]\&A\ar[d,"x"]\\
  A\ar[r,"x^3"]\ar[d,"1"]\&A\ar[d,"x"]\\
  \vdots\&\vdots
\end{tikzcd}
\end{eqnarray}
which induces maps on the cohomology level, $x\colon H^i(A;x^\ell)\to
H^i(A;x^{\ell+1})$ and hence a direct system of cohomology modules 
over the index set $\NN$. It is an instructive exercise (using the
fact that $\NN$ is an index set that satisfies: for all $ n,n'\in\NN$
there is $N\in\NN$ exceeding both $n,n'$) to check that the direct
limit $\dlim_\ell H^k(A;x^\ell)$ agrees with the local cohomology
module $H^k_{\ideal{x}}(A)$.
  
If $M$ is an $A$-module and the ideal $I$ is generated by
$x_1,\ldots,x_m$ then to each such generating set there is a
\emph{cohomological Koszul complex}
\[
K^\bullet(M;x_1,\ldots,x_m):=M\otimes_A\bigotimes_{i=i}^mK^\bullet(A;x_i)
\]
whose cohomology modules are denoted
$H^\bullet(M;x_1,\ldots,x_m)$. Again, one can verify that replacing
each $x_i$ by powers of themselves leads to a tower of complexes whose
direct limit has a cohomology that functorially equals the local
cohomology $H^\bullet_I(M)$. In particular, it is independent of the
chosen generating set for $I$.

\subsection{The \v Cech complex}

Inspection shows that the direct limit of the tower
$A\stackrel{x}{\to}A\stackrel{x}{\to}A\stackrel{x}{\to}\cdots$ is
functorially equal to the localization $A[x^{-1}]$ which we also write
as $A_x$. Thus, the limit complex to the tower \eqref{eqn-tower-x} is
the localization complex $A\to A_x$.  In greater generality, the
module that appears in the limit complex $\check
C^\bullet(M;x_1,\ldots,x_m)$ of the tower
$K^\bullet(M;x_1,\ldots,x_m)\to K^\bullet(M;x_1^2,\ldots,x_m^2)\to
K^\bullet(M;x_1^3,\ldots,x_m^3)\to \cdots$ in cohomological degree $k$
is the direct sum of all localizations of $M$ at $k$ of the $m$
elements $x_1,\ldots,x_m$. Hence,
\[
\check C^\bullet(M;x_1,\ldots,x_m)=\dlim_\ell
K^\bullet(M;x_1^\ell,\ldots,x_m^\ell)
\]
and a corresponding statement links the cohomology modules on both sides.

The point of view of the \v Cech complex provides a useful link to
projective geometry. Indeed,
suppose $I\subseteq R=\KK[x_1,\ldots,x_n]$ is the homogeneous ideal 
defining the projective variety $X$ in
$\PP^n_\KK$. Then the \emph{cohomological dimension} $\cd(U)$
of  $U=\PP^n_\KK\setminus X$, the largest integer $k$ for which
$H^k(U,-)$ is not the zero functor on the category of quasi-coherent
sheaves on $U$,  equals $\lcd_R(I)-1$. This follows from
the exact sequence
\begin{gather}\label{eqn-cd-Serre}
0\to\Gamma_I(M)\to M\to \bigoplus_{k\in\ZZ}\Gamma(\PP^n_\KK\setminus
X,\tilde M(k))\to H^1_I(M)\to 0
\end{gather}
and the isomorphisms $\bigoplus_{k\in\ZZ}H^i(\PP^n_\KK\setminus
X,\tilde M(k))=H^{i+1}_I(M)$ for any $R$-module $M$ with associated
quasi-coherent sheaf $\tilde M$.

\subsection{Limits of Ext-modules}

Again, let $I=\ideal{x_1,\ldots,x_m}$ be an ideal of $A$.  The natural
projections $A/I^{\ell+1}\to A/I^\ell$ lead to a natural tower of
morphisms $\Ext^k_A(A/I,M) \to \Ext^k_A(A/I^2,M) \to\Ext^k_A(A/I^3,M)
\to\cdots$. An exercise involving $\delta$-functors (also known as
connected sequences of functors) shows that the
direct limit of this system functorially agrees with $H^k_I(M)$.

We have thus the functorial isomorphisms
\[
H^k_I(M)\simeq \dlim_\ell H^k(M;x_1^\ell,\ldots,x_m^\ell)\simeq H^k\check
C^\bullet(M;x_1,\ldots,x_m)\simeq \dlim_\ell\Ext^k_A(A/I^\ell,M)
\]
for all choices of generating sets $x_1,\ldots,x_m$ for $I$.

\begin{rmk}\label{rmk-other-powers}
  \begin{asparaenum}
  \item The derived functor version of local cohomology shows that
    $H^\bullet_I(-)$ and $H^\bullet_J(-)$ are the same functor
    whenever $I$ and $J$ have the same radical.
\item It follows easily from the \v Cech complex interpretation that
  local cohomology satisfies a local-to-global principle: for any
  multiplicatively closed subset $S$ of $A$ one has $S^{-1}\cdot
  H^i_I(M)=H^i_{I(S^{-1}A)}(S^{-1}M)$, and so in particular
  $H^i_I(M)=0$ if and only if $H^i_{IA_\frakp}(M_\frakp)=0$ for all
  $\frakp\in\Spec A$.
\item If $I$ is (up to radical) a complete intersection in the
  localized ring $A_\frakp$, then $H^k_I(A)\otimes_A A_\frakp$ is zero
  unless $k=\hight(IA_\frakp)$. If $R$ is a regular
  local ring and $I$ reduced then $I$ is a complete
  intersection in every smooth point. It follows that for
  equidimensional $I$ the support of
  $H^k_I(R)$ with $k>\hight(I)$ only contains primes $\frakp$ contained
  in the singular locus of $I$.
  \item It is in general a difficult question to predict how the
    natural maps $\Ext^k_A(A/I^\ell,M)\to H^k_I(M)$ and
    $H^k(M;x_1^\ell,\ldots,x_m^\ell)\to H^k_I(M)$ behave; some
    information can be found in
    \cite{EisenbudMustataStillman,SinghWalther-pure,BBLSZ2,MaSinghWalther-Dutta}.
  \item If $\phi\colon A'\to A$ is a ring morphism, and if $M$ is an
    $A$-module and $I'$ an ideal of $A'$, then there is a functorial
    isomorphism between $H^k_{I'}(\phi_*M)$ and
    $\phi_*(H^k_{I'A}(M))$, where $\phi_*$ denotes restriction of
    scalars from $A$ to $A'$. The easiest way to see this is by comparison
    of the two \v Cech complexes involved.
  \end{asparaenum}
\end{rmk}

\begin{rmk}
\label{rmk: cofinal sequences}
Let $I$ be an ideal of a Noetherian commutative ring $A$. A sequence
of ideals $\{I_k\}$ is called {\it cofinal} with the sequence of
powers $\{I^k\}$ if, for all
    $k\in\NN$, there are $\ell, \ell'\in\NN$ such that both $I_{\ell}\subseteq I^k$ and $I^{\ell'} \subseteq I_{k}$.
    
Sequences $\{I_k\}$ cofinal with $\{I^k\}$ are of interest in the study of local cohomology since 
\[\varinjlim_{k} \Ext^i_R(R/I_k,M)\cong \varinjlim_{k} \Ext^i_R(R/I^k,M)=H^i_I(M).\] 
This provides one with the flexibility of using sequences of ideals
other than $\{I^n\}$.
In characteristic $p>0$, the sequence of ideals defined next plays an
extraordinary part in the story.

Let $A$ be a Noetherian commutative ring of prime characteristic $p$
and $I$ be an ideal of $A$. The \emph{$e$-th Frobenius power} of $I$, denoted by $I^{[p^e]}$, is defined to be the ideal generated by the $p^e$-th powers of all elements of $I$. Since the Frobenius endomorphism $A\xrightarrow{a\mapsto a^p}A$ is a ring homomorphism, $I^{[p^e]}=(f^{p^e}_1,\dots,f^{p^e}_t)$ for every set of generators $\{f_1,\dots,f_t\}$ of $I$.

It is straightforward to check that $\{I^{[p^e]}\}$ is cofinal with
$\{I^k\}$ since $A$ is Noetherian and thus $I$ is finitely generated. 
\end{rmk}

\subsection{Local duality}

Matlis duality over a complete local ring $(A,\frakm,\kk)$ provides a
one-to-one correspondence between the Artinian and the Noetherian
modules over $A$; in both directions it is given by the functor
\[
D(M):= \Hom_A(M,E_A(\kk))
\]
of homomorphisms into the injective hull of the residue
field.\footnote{Strictly speaking, one should write $D_A(-)$, but in
  all cases the underlying ring will be understood from the context.}
Of course, one can in principle  apply $D(-)$ to any module, but the property
$D(D(M))=M$ is likely to fail when $M$ does not enjoy any finiteness
condition.

A natural question is what the result of applying $D(-)$ to
$H^i_\frakm(A)$ should be or, more generally, how to describe
$D(H^i_\frakm(M))$ for Noetherian $A$-modules $M$. It turns out that
when $A$ ``lends itself to duality'', then this question has a
pleasing answer:
\begin{thm}
  Suppose $(A,\frakm,\kk)$ is a local Gorenstein ring. Then
  \[
  D(H^i_\frakm(M))\cong \Ext^{\dim(A)-i}_A(M,A)
  \]
  for every finitely generated $A$-module $M$.
\end{thm}
The original version is due to Grothendieck
\cite{Hartshorne-lc-notes}, and then expanded in Hartshorne's opus
\cite{Hartshorne-RD}.  As it turns out, there are extensions of local
duality to Cohen--Macaulay rings with a dualizing module, and yet more
generally to rings with a dualizing complex. Duality on formal or
non-Noetherian schemes and other generalizations are discussed in
\cite{TarrioJeremiasLipman}.

In particular, Chapter 4 of \cite{Hartshorne-RD}
contains a discussion on Cousin complexes and their connection to
local cohomology, that we do not have the space to give justice to.
Further accounts in this direction can be found in
\cite{Sharp-Cousin,Sharp-Cousin+lc, Schenzel-Dualisierende,Lipman-Guanajuato}.

\section{Finiteness and vanishing}
\label{sec-finiteness}

\subsection{Finiteness properties.}

In general, local cohomology modules are not finitely generated. For
instance, the Grothendieck nonvanishing theorem says:

\begin{thm}
Let $(A,\frakm)$ be a Noetherian local ring and $M$ be a nonzero finitely generated $A$-module. Then $H^{\dim(M)}_{\frakm}(M)\neq 0$. 

Moreover, if $\dim(M)>0$, then $H^{\dim(M)}_{\frakm}(M)$ is not finitely generated.
\end{thm}

Finiteness is more unusual yet than this theorem indicates. For
example, over a ring $R$ of polynomials over $\CC$, a local
cohomology module $H^k_I(R)$ is  a finite $R$-module precisely if
$I=0$ and $k=0$, or if $H^k_I(R)=0$. 
This lack of finite generation prompted people to look at other types
of finiteness properties, and in this section we survey various
fruitful avenues of research that pertain to finiteness.

\medskip

In \cite[expos\'e 13, 1.2]{SGA2} Grothendieck conjectured that, if $I$ is an ideal in a Noetherian local ring $A$, then $\Hom_A(A/I, H^j_I(A))$ is finitely generated. Hartshorne refined this finiteness of $\Hom_A(A/I,H^j_I(A))$ and introduced the notion of {\it cofinite modules} in \cite{Hartshorne-AffineDuality}.

\begin{dfn}
Let $A$ be a Noetherian commutative ring and $I\subseteq A$ an ideal. An $A$-module $M$ is called {\it $I$-cofinite} if $\Supp_A(M)\subseteq V(I)$ and $\Ext^i_A(A/I,M)$ is finitely generated for all $i$.
\end{dfn}

In \cite{Hartshorne-AffineDuality} Hartshorne constructed the following example which answered Grothendieck's conjecture on finiteness of $\Hom_A(A/I,H^j_I(M))$ in the negative.

\begin{exa}
Let $\kk$ be a field and put $A=\frac{\kk[x,y,u,v]}{(xu-yv)}$. Set $\fraka=\ideal{x,y}$ and $\frakm=\ideal{x,y,u,v}$. Then $\Hom_A(A/\frakm, H^2_\fraka(A))$ is not finitely generated and hence neither is $\Hom_A(A/\fraka, H^2_\fraka(A))$.

We note in passing, that while the socle dimension of $H^2_\fraka(A)$
is infinite, it is nonetheless a finitely generated module over the
ring of $\kk$-linear differential operators on $A$,
\cite{Jenchieh-TAMS12}.
\end{exa}

The ring $A$ in Hartshorne's example is not regular; one may ask
whether local cohomology modules $H^i_I(R)$ of a Noetherian regular
ring $R$ are $I$-cofinite. Huneke and Koh showed in \cite{HunekeKoh-cofiniteness} that this is not the case even for a polynomial ring over a field.

\begin{exa}
\label{exa-2x3-socle}
Let $\kk$ be a field of characteristic 0 and let $R=\kk[x_{1,1},\dots,
  x_{2,3}]$ be the polynomial ring over $\kk$ in 6 variables. Set $I$
to be the ideal generated by the $2\times 2$ minors of the matrix
$(x_{ij})$.

The geometric origins and connections of this example, including a
discussion of the interaction of the relevant local cohomology groups
with de Rham cohomology and $D$-modules, can be found in Examples
\ref{exa-2x3}, \ref{exa-2x3-betti} and Remark \ref{rmk-2x3} below. In
particular, Example \ref{exa-2x3-betti} discusses that $H^3_I(R)$ is
isomorphic to the injective hull of $\kk$ over $R$, which means that
$\Hom_R(R/I,H^3_I(R))$ is the injective hull of $\kk$ over $R/I$ and
thus surely not finitely generated.
\end{exa}

Huneke and Koh further proved in \cite{HunekeKoh-cofiniteness} that:
\begin{thm}
Let $R$ be a regular local ring and $I$ be an ideal in $R$. Set $b$ to
be the biggest height of any minimal prime of $I$ and set
$c=\lcd_R(I)$, compare Definition \ref{dfn-lcd}.
\begin{enumerate}
\item If $R$ contains a field of characteristic $p>0$ and if $j>b$ is an integer such that $\Hom_R(R/I,H^j_I(R))$ is finitely generated, then $H^j_I(R)=0$.
\item If $R$ contains $\QQ$ 
then $\Hom_R(R/I,H^c_I(R))$ is not finitely generated.
\end{enumerate}
\end{thm}

In Example \ref{exa-2x3-socle}, it turns out that the socle
$\Hom_R(R/\frakm,H^3_I(R))$ of $H^3_I(R)$ is finitely generated. It
is natural to ask whether the socle of local cohomology of a
Noetherian regular ring is always finitely generated; as a matter of
fact this was precisely \cite[Conjecture 4.3]{HunekeProblemsLC}.

In \cite{HunekeProblemsLC}, Huneke proposed a number of problems on
 local cohomology which guided the study of local cohomology modules
 for decades.

\begin{prb}[Huneke's List]~
\begin{enumerate}
\item[1.] When is $H^j_I(M)=0$?
\item[2.] When is $H^j_I(M)$ finitely generated?
\item[3.] When is $H^j_I(M)$ Artinian?
\item[4.] If $M$ is finitely generated, is the number of associated primes of $H^j_I(M)$ always finite?
\end{enumerate}
\end{prb}
Huneke remarked that ``all of these problems are connected with
another question:
{\it
\begin{enumerate}
\item[5.] What annihilates the local cohomology module $H^j_I(M)$? 
\end{enumerate}}

More concretely, Huneke conjectured:
\begin{conj}[Conjectures 4.4 and 5.2 in \cite{HunekeProblemsLC}]
\label{Huneke Conjecture on Bass and associated primes}
Let $R$ be a regular local ring and $I$ be an ideal. Then
\begin{enumerate}
\item the Bass numbers $\Ext^i_{R_{\frakp}}(\kappa(\frakp), H^j_I(R_{\frakp}))$ are finite for all $i$, $j$, and prime ideals $\frakp$, and 
\item the number of associated primes of $H^j_I(R)$ is finite for all $j$.
\end{enumerate}
\end{conj}

Later in \cite{L-Dmods} Lyubeznik conjectured further that the finiteness of associated primes holds for local cohomology of all Noetherian regular rings. Substantial progress has been made on these conjectures. If the regular ring has prime characteristic $p>0$, then these conjectures were completely settled by Huneke and Sharp in \cite{HunekeSharp}; in equi-characteristic 0, Lyubeznik proved these conjectures for two large classes of regular rings in \cite{L-Dmods}; for complete unramified regular local rings of mixed characteristic, these conjectures were first settled by Lyubeznik in \cite{LyubeznikUnramified} (different proofs can also be found in \cite{Betancourt-IJM13} and \cite{BBLSZ-2014}). The finiteness of associated primes of local cohomology was also proved in \cite{BBLSZ-2014} for smooth $\ZZ$-algebras. We summarize these results as follows.

\begin{thm}
\label{thm: finiteness known regular cases}
Assume that $R$ is 
\begin{enumerate}
\item a Noetherian regular ring of characteristic $p>0$, or
\item a complete regular local ring containing a field of characteristic $0$, or
\item regular of finite type over a field of characteristic $0$, or
\item an unramified regular local ring of mixed characteristic, or
\item a smooth $\ZZ$-algebra.
\end{enumerate}
Then the Bass numbers and the number of associated primes of $H^j_I(R)$ are finite for every ideal $I$ of $R$ and every integer $j$.
\end{thm} 

\begin{rmk}
When $R$ is a smooth $\ZZ$-algebra, then finiteness of Bass numbers
was not addressed in \cite{BBLSZ-2014}. However, one can conclude
readily from the unramified case in \cite{LyubeznikUnramified} as
follows. The Zariski-local structure theorem for smooth morphism says
that $\ZZ\to R$ factors as a composition of a polynomial extension and
a finite etale morphism, which implies that locally $R$ is an
unramified regular local ring of mixed characteristic. Since the
finiteness of Bass numbers is a local problem, the desired conclusion  follows from the results in \cite{LyubeznikUnramified}.
\end{rmk}

Conjecture \ref{Huneke Conjecture on Bass and associated primes} is
still open when $R$ is a ramified regular local ring of mixed
characteristic. Theorem \ref{thm: finiteness known regular cases}(1)
was proved in \cite{HunekeSharp} using properties of the Frobenius endomorphism; this approach was later conceptualized by Lyubeznik to his theory of $F$-modules in \cite{LyubeznikFModules}. The proof of Theorem \ref{thm: finiteness known regular cases}(2)-(5) uses $D$-modules ({\it i.e.} modules over the ring of differential operators). Both $F$-modules and $D$-modules will be discussed in the sequel.

For a non-regular Noetherian ring $A$, if $\dim(A)\leq 3$
(\cite{Marley-finiteness-ass-primes}), or if $A$ is a 4-dimensional
excellent normal local domain (\cite{HunekeKatzMarley}), then the
number of associated primes of $H^j_I(M)$ is finite for every finitely
generated $A$-module $M$, for every ideal $I$ and for all integers $j$. Once the
restriction on $\dim(A)$ is removed, then the number of associated
primes of local cohomology modules can be infinite; such examples have been discovered in \cite{Singh-p-torsion-elements, Katzman-infinite-ass-primes, SinghSwanson-ass-primes}. Note that all these examples are hypersurfaces; the hypersurface in \cite{SinghSwanson-ass-primes} has rational singularities. 

\smallskip

As local cohomology modules may have infinitely many associated primes in general, one may ask a weaker question (\cite[p.~3195]{HunekeKatzMarley}):
\begin{qun}
\label{question: close support}
Let $A$ be a Noetherian ring, $I$ be an ideal of $A$ and $M$ be a finitely generated $A$-module. Does $H^j_I(M)$ have only finitely many {\it minimal} associated primes? Or equivalently, is the support of $H^j_I(M)$ Zariski-closed?
\end{qun}

It is stated in \cite{HunekeKatzMarley} that ``this question is of central
importance in the study of cohomological dimension and understanding the local-global properties of local cohomology". 

When $\dim(A)\leq 4$, then Question \ref{question: close support} has a positive answer due to \cite{HunekeKatzMarley}. If $\mu(I)$ denotes the number of generators of $I$ and $A$ has prime characteristic $p$, it was proved and attributed to Lyubeznik in \cite{KatzmanSupportTopLC} that $H^{\mu(I)}_I(A)$ has a Zariski-closed support. When $A=R/(f)$ where $R$ is a Noetherian ring of prime characteristic $p$ with isolated singular closed points, it was proved independently in \cite{KatzmanZhangSupportLC} and in \cite{HochsterNunezSupportLC} that $H^j_I(A)$ has a Zariski-closed support for every ideal $I$ and integer $j$. 

\medskip

A classical result in commutative algebra ({\it c.f.} \cite[Theorem 3.1.17]{BrunsHerzogCMBook}) says that if $A$ is a Noetherian local ring and $M$ is a finitely generated $A$-module $M$ that has finite injective dimension, then 
\[\dim(M) \leq \injdim_A(M)=\depth(A)\] 
where $\injdim_A(M)$ denotes the injective dimension of
$M$ over $A$. Interestingly, for local cohomology modules over regular rings,
the inequality seems to be reversed. More precisely, the following was proved in \cite{HunekeSharp} and \cite{L-Dmods}
\begin{thm}\label{thm-injdim}
Assume that $R$ is 
\begin{enumerate}
\item a Noetherian regular ring of characteristic $p>0$, or
\item a complete regular local ring of characteristic $0$, or
\item regular of finite type over a field of characteristic $0$.
\end{enumerate}
Then 
\[\injdim_R(H^j_I(R))\leq \dim(\Supp_R(H^j_I(R)))\]
for every ideal $I$ and integer $j$.
\end{thm}

In \cite{PuthenpurakalInjectiveResolutionLC}, Puthenpurakal showed
that if $R=\kk[x_1,\dots,x_n]$ where $\kk$ is a field of
characteristic 0 then $\injdim_R(H^j_I(R))= \dim(\Supp_R(H^j_I(R)))$
for every ideal $I$. Later this was strengthened in \cite[Theorem
  1.2]{ZhangInjectiveDim} as follows: assume that either $R$ is a
regular ring of finite type over an infinite field of prime
characteristic $p$ and $M$ is an $F$-finite $F$-module, or
$R=\kk[x_1,\dots,x_n]$ where $\kk$ is a field of characteristic 0 and
$M$ is a holonomic\footnote{The notions of holonomic $D$-modules and
  $F$-finite $F$-modules will be explained in the sequel; local
  cohomology modules with argument $R$ and $R$ as discussed here are
  primary examples of those.}  $D$-module. Then
\[\injdim_R(M)= \dim(\Supp_R(M)).\]

Subsequently \cite{SwitalaZhangInjectiveDim} proved that, if M is
either a holonomic $D$-module over a formal power series ring $R$ with
coefficients in a field of characteristic 0, or an $F$-finite
$F$-module over a Noetherian regular ring $R$ of prime characteristic $p$, then
\[\dim(\Supp_R(M))-1\leq \injdim_R(M)\leq \dim(\Supp_R(M)).\]

When the regular ring $R$ does not contain a field, the bounds on injective dimension of local cohomology modules of $R$ are different. In \cite{ZhouJA1998}, Zhou proved that, {if $(R, \frakm)$ is an unramified regular local ring of mixed characteristic and $I$ is an ideal of $R$, then $\injdim_R(H^j_I(R))\leq \dim(\Supp_R(H^j_I(R)))+1$ and $\injdim_R(H^i_{\frakm}H^j_i(R))\leq 1$.} Moreover, it may be the case that $\injdim_R(H^i_{\frakm}H^j_i(R))=1$, as shown in \cite{BetancourtHernandezPerezWitt-TAMS19,DattaSwitalaZhang}.

\subsection{Vanishing}\label{subsec-vanishing}
Problem 1 in Huneke's list of problems in \cite{HunekeProblemsLC}
asks: when is $H^j_I(M)=0$? Vanishing results on local cohomology
modules have a long and rich history. Note that $H^j_I(M)=0$ for all
$j>t$ and all $A$-modules $M$ if and only if $H^j_I(A)=0$ for
$j>t$. Recall the notion of local cohomological dimension from
Definition \ref{dfn-lcd}. For a Noetherian local ring $A$, we set
\[
\mdim(A)=\min\{\dim(A/Q)\mid Q \text{ is a minimal prime of }A\}
\]
and we write $\embdim(A)$ for the embedding dimension (the number of
generators of the maximal ideal) of a local ring $A$. For an
ideal $I$ of $A$, we set
\[
c_A(I)=\embdim(A)-\mdim(A/I).
\]
Note that if $A$ is regular then $c_A(I)$ is called the big height, {\it i.e.} the biggest height of any minimal prime ideal of $I$.

We now summarize the most versatile vanishing theorems on local cohomology.
\begin{itemize}
\item (Grothendieck Vanishing) Let $A$ be a Noetherian ring and $M$ be a finitely generated $R$-module. Then $H^j_I(M)=0$ for all integers $j>\dim(M)$ and ideals $I$. In particular, this implies that $\lcd_A(I)\leq \dim(A)$ for all ideals $I$.

\item (Hartshorne--Lichtenbaum Vanishing) Let $(A,\frakm)$ be a Noetherian local ring and $I$ be an ideal of $R$. Then $\lcd_A(I)\leq \dim(A)-1$ if and only if $\dim(\hat{A}/(I\hat{A}+P))>0$ for every minimal prime $P$ of $\hat{A}$ such that $\dim(\hat{A}/P)=\dim(A)$, where $\hat{A}$ denotes the completion of $A$. In particular, this implies that if $A$ is a complete local domain and $\sqrt{I}\neq \frakm$ then $\lcd_A(I)\leq \dim(A)-1$. \emph{cf.} \cite{Hartshorne-CDAV}.

\item (Faltings Vanishing) Let $R$ be a complete equi-characteristic
regular local ring with a separably closed residue field. Then
\[
\lcd_R(I)\leq \dim(R)-\left\lfloor\frac{\dim(R)-1}{c_R(I)}\right\rfloor,
\]
{\it cf.} \cite{FaltingsCrelle1980}.
\footnote{This is the
floor function $\lfloor x\rfloor=\max\{k\in\ZZ, k\le x\}$.}
This bound is sharp according to \cite{Lyubeznik-PAMS1985}.

\item (Second Vanishing Theorem) Let $R$ be a complete regular local ring that contains a separably closed coefficient field and $I$ be an ideal. Then $\lcd_R(I)\leq \dim(R)-2$ if and only if $\dim(R/I)\geq 2$ {\it and} the punctured spectrum of $R/I$ is connected.\\ 
A version of this vanishing theorem for projective varieties was first obtained by Hartshorne in \cite[Theorem 7.5]{Hartshorne-CDAV} who coined the name `Second Vanishing Theorem'. The local version stated here was left as a problem by Hartshorne in \cite[p.~445]{Hartshorne-CDAV}. Subsequently, this theorem was proved in prime characteristic in \cite{PeskineSzpiro-IHES73}, in equi-characteristic 0 in \cite{Ogus-LCDAV} (a unified proof for equi-characteristic regular local rings can be found in \cite{HunekeLyubeznik}), and for unramified regular local rings in mixed characteristic in \cite{ZhangVanishingLCMixedChar}.

\item (Peskine--Szpiro Vanishing) Let $(R,\frakm)$ be a Noetherian
regular local ring of prime characteristic $p$ and $I$ be an
ideal. Then $\lcd_R(I)\leq \dim(R)-\depth(R/I)$, \emph{cf.} \cite{PeskineSzpiro-IHES73}.

\item (Vanishing via action of Frobenius) Let $(R,\frakm)$ be a regular local ring of prime characteristic $p$ and $I$ be an ideal. Set $d=\dim(R)$. Then $H^j_I(R)=0$ if and only if the Frobenius endomorphism on $H^{d-j}_{\frakm}(R/I)$ is nilpotent. {\it cf.} \cite{Lyubeznik-Compositio06}.
\end{itemize}

There have been various extensions of the vanishing theorems mentioned
above. Most notably, \cite{HunekeLyubeznik} initiated an investigation
on finding bounds of local cohomological dimension under topological
and/or geometric assumptions. For instance, \cite[Theorem 3.8]{HunekeLyubeznik} asserts that if $A$ is a complete local ring containing a field and $I$ is a formally geometrically irreducible ideal such that $0<c_A(I)<\dim(A)$ then
\begin{equation}
\label{HL type bound cd}
\lcd_A(I)\leq \dim(A)-1-\left\lfloor\frac{\dim(A)-2}{c_A(I)}\right\rfloor.
\end{equation}
Furthermore, if $A/I$ is normal then
\[
\lcd_A(I)\leq \dim(A)-\left\lfloor\frac{\dim(A)+1}{c_A(I)+1}\right\rfloor-
\left\lfloor\frac{\dim(A)}{c_A(I)+1}\right\rfloor.
\]
The bound on cohomological dimension in (\ref{HL type bound cd}) was later extended to reducible ideals in \cite{Lyubeznik-AdvMath2007} as follows.
\begin{thm}
\label{Lyubeznik vanishing simplicial complex}
Let $(A,\frakm,\kk)$ be a $d$-dimensional local ring containing $\kk$.  Assume
$d>1$. Let $c$ be a positive integer, let
$t=\left\lfloor(d-2)/c\right\rfloor$ and
$v=d-1-\left\lfloor(d-2)/c\right\rfloor$. Let $I$ be an ideal of $A$
with $c(I\hat{A})\leq c$. Let $B$ be the completion of the strict
Henselization of the completion of $A$. Let $I_1,\dots,I_n$ be the minimal primes of $IB$ and let $P_1,\dots,P_m$ be the primes of $B$ such that $\dim(B/P_i)=d$. Let $\Delta_i$ be the simplicial complex on $n$ vertices $\{1,2,\dots,n\}$ such that a simplex $\{j_0,\dots,j_s\}$ belongs to $\Delta_i$ if and only if $I_{j_0}+\cdots +I_{j_s}+P_i$ is not $\frakm B$-primary. Let $\tilde{H}_{t-1}(\Delta_i;k)$ be the $(t-1)$st singular homology group of $\Delta_i$ with coefficients in $k$. Then $\lcd_A(I)\leq v$ if and only if $\tilde{H}_{t-1}(\Delta_i;k)=0$ for every $i$.
\end{thm}

We ought to point out that the simplicial complex introduced in
Theorem \ref{Lyubeznik vanishing simplicial complex} has spurred a
line of research on connectedness dimension, {\it
cf.}  \cite{KLZ-PAMS2016}, \cite{DaoTakagi-Compositio16}, \cite{BetancourtSpiroffWitt}, \cite{Varbaro-CiA2019}.

Also, \cite{DaoTakagi-Compositio16} shows that the same bound as in (\ref{HL type bound cd}) holds when $A$ is a complete regular local ring containing a field such that $A/I$ has positive dimension and satisfies Serre's condition $(S_2)$. This result is in the spirit of a question raised by Huneke in \cite{HunekeProblemsLC}.

\begin{qun}[Huneke]
\label{Huneke-question-Serre-conditions}
Let $R$ be a complete regular local ring with separably closed residue field and $I$ be an ideal of $R$. Assume that $R/I$ satisfies Serre's conditions $(S_i)$ and $(R_j)$. What is the maximal possible cohomological dimension for such an ideal?
\end{qun}

In the same spirit as Huneke's Question \ref{Huneke-question-Serre-conditions}, one may ask about the possibility of an
implication 
\begin{gather}\label{eqn-Varbaro-cnj}
[\depth_R(R/I)\geq t]\quad\stackrel{?}{\Longrightarrow}
\quad[\lcd_R(I)\le\dim(R)-t].
\end{gather}

In prime characteristic $p$, such implication holds due to
Peskine--Szpiro Vanishing. On the other hand, Peskine--Szpiro
Vanishing can fail in characteristic 0.

\begin{exa}
\label{exa-2x3}
  Let $R$ be the polynomial ring in the variables
  $x_{1,1},\ldots,x_{2,3}$ over the field $\KK$, localized at
  $\bsx=\ideal{x_{1,1},\ldots,x_{2,3}}$. Let $I$ be the ideal of
  maximal minors of the matrix $(x_{i,j})$. Then $I$ is the radical
  ideal associated to 
  the 4-dimensional locus $V$ of the $2\times 3$ matrices of rank one, which agrees
  with the image of the map $\KK^2\times \KK^3\to \KK^{2\times 3}$ that
  sends $((s,t),(x,y,z))$ to $(xs,ys,zs,xt,yt,zt)$. In particular, $I$
  is the prime ideal associated to the image of the Segre embedding of
  $\PP^1_\KK\times \PP^2_\KK\into \PP^5_\KK$. 

  Thus $I$ is $3$-generated of height $2=6-4$, and in fact $R/I$ is
  Cohen--Macaulay of depth $4$. Since the origin is the only singular
  point of $V$, the local cohomology groups $H^k_I(R)$ are supported
  at the origin for $k\neq 2$ and zero for
  $k\not\in\{2,3\}$. Cohen--Macaulayness of $R/I$ forces the vanishing of
  $H^k_I(R)$ for $k\neq 2$ in prime characteristic, but if the
  characteristic of $\KK$ is zero then $H^3_I(R)$ is actually nonzero.  For a
  computational discussion involving $D$-modules see
  \cite{Oaku-Duke,W-lcD,OT2,W-M2book}.  We will return to this
  situation in Example \ref{exa-2x3-betti}.
\end{exa}

In Example \ref{exa-2x3}, $\depth(R/I)=4$ but $H^3_I(R)\neq 0$. This
shows that the implication \eqref{eqn-Varbaro-cnj} can fail in
characteristic 0 when $t\geq 4$. When $t\le2$, the implication
\eqref{eqn-Varbaro-cnj} holds due the Second Vanishing Theorem and the
Hartshorne--Lichtenbaum Theorem. The case $t=3$ is not completely
settled, but there have been positive results. In
\cite{Varbaro-Compositio13}, continuing his work on the number of
defining equations in \cite{Varbaro-TAMS12}, Varbaro proved that if a
homogeneous ideal $I$ in a polynomial ring $R=\kk[x_1,\dots,x_n]$ over
a field $\kk$ satisfies $\depth(R/I)\geq 3$ then $\lcd_R(I)\leq
n-3$. He also conjectured:
\begin{conj}[Varbaro]
\label{Varbaro conj}
Let $R$ be a regular local ring containing a field and $I$ be an ideal of $R$. If $\depth(R/I)\geq 3$, then $\lcd_R(I)\leq \dim(R)-3$.
\end{conj}

The fact that Implication \eqref{eqn-Varbaro-cnj} can fail for complex
projective threefolds (Example \ref{exa-2x3}) raises the question what
exact features are responsible for failure when $t=3$. Clearly, more
knowledge about the singularity is required than just $\depth_R(R/I)$.

Dao and Takagi prove Conjecture \ref{Varbaro conj} in
\cite{DaoTakagi-Compositio16} when $R$ is essentially of finite type
over a field. More specifically, they show the following facts about
the inequality $\lcd_R(I)\le \dim(R)-3$. Suppose $R$ is a regular
local ring essentially of finite type over its algebraically closed
residue field of characteristic zero. Take an ideal $I$ such that
$R/I$ has depth $2$ or more and $H^2_\frakm(R/I)$ is a $\KK$-vector
space (\emph{i.e.}, it is killed by $\frakm$).  Then $\lcd_R(I)\le
\dim R-3$ if and only if the torsion group of $\Pic(\Spec(R/I))$ is
finitely generated on the punctured completed spectrum.  In case that
the depth of $R/I$ is at least $4$, one even has $\lcd_R(I)\le
\dim(R)-4$ if and only if the Picard group is torsion on the punctured
completed spectrum of $R/I$.  In Example \ref{exa-2x3}, the depth of
$R/I$ is four, but the Picard group on the punctured spectrum is not
torsion but $\ZZ$.  Conjecture \ref{Varbaro conj} remains
open in general.

\medskip

Both the Hartshorne--Lichtenbaum Vanishing Theorem and the Second
Vanishing Theorem may viewed as topological criteria for vanishing and
have applications to topology of algebraic varieties ({\it cf.}
\cite{L-ecd} and \cite{Huneke-lc-notes}). It would be desirable to
have an analogue of the Second Vanishing Theorem for non-regular rings. In \cite[p.~144]{L-lc-survey} Lyubeznik asked the following questions.
\begin{qun}
\label{Lyubeznik question second vanishing}
Let $(A, \frakm)$ be a complete local domain with a separably closed residue field. 
\begin{enumerate}
\item Find necessary and sufficient conditions on $I$ such that $\lcd_A(I)\leq \dim(A)-2$.
\item Let $I$ be a prime ideal. Is it true that $\lcd_A(I)\leq \dim(A)-2$ if and only if $(P + I)$ is not primary to the maximal ideal for any prime ideal $P$ of height 1?
\end{enumerate}
\end{qun}
Question \ref{Lyubeznik question second vanishing}(1) remains open. It turns out that Question \ref{Lyubeznik question second vanishing}(2) has a negative answer due to \cite[Proposition 7.7]{HochsterZhang-TAMS2018}: 

\begin{exa}
Let $A=\frac{\CC[[x,y,z,u,v]]}{(x^3+y^3+z^3,z^2-ux-vy)}$ and $I=(x,y,z)$. Then
\begin{enumerate}
\item $\dim(A)=3$ and $\hight(I)=1$;
\item $I+P$ is {\it not} primary to the maximal ideal for every height-1 prime ideal $P$;
\item $H^2_I(A)\neq 0$. 
\end{enumerate}
\end{exa}

Given the connections between local cohomology and sheaf cohomology ({\it cf.} (\ref{eqn-cd-Serre})), vanishing of sheaf cohomology can be interpreted in terms of local cohomology. The classical Kodaira Vanishing Theorem asserts that: {\it If $X$ is smooth projective variety over a field $\KK$ of characteristic 0, then $H^i(X,\calO(j))=0$ for $i<\dim(X)$ and all $j<0$.} This result has an equivalent formulation in terms of local cohomology: {\it If $R$ is a standard\footnote{A standard graded algebra over a field $\KK$ is a graded quotient of a polynomial ring over $\KK$ with the standard grading.} graded domain over a field $\KK$ of characteristic 0 such that $\Proj(R)$ is smooth, then $H^j_{\frakm}(R)_{<0}=0$ for all $j<\dim(R)$, where $\frakm$ is the homogeneous maximal ideal of $R$.} For ideal-theoretic interpretations and connections with tight closure and Frobenius, we refer the interested reader to \cite{HunekeSmith-KodairaVanishing-1997, Smith-FujitaConjecture-JAG1997}.

The Kodaira Vanishing Theorem fails for singular varieties in
characteristic 0 and also fails for smooth varieties in characteristic
$p$. It is proved in \cite{BBLSZ2} that, if one focuses on the range
$i<\codim(\Sing(X))$, then the Kodaira Vanishing Theorem can be extended to thickenings of local complete intersections. More precisely:

\begin{thm}
\label{Kodaira for lci thickenings}
Let $X$ be a closed local complete intersection subvariety of $\PP^n_\KK$
over a field $\KK$ of characteristic 0 and let $I$ be its defining ideal. Let $X_t$ denote the scheme defined by $I^t$. Then
\[H^i(X_t,\calO_{X_t}(j))=0\]
for all $i<\codim(\Sing(X))$, all $t\geq 1$, and all $j<0$.

Or, equivalently, let $S=\KK[x_0,\dots,x_n]$ and $I$ be as above. Then
\[H^{\ell}_{\frakm}(S/I^t)_{<0}=0\]
for $\ell<\codim(\Sing(X))+1$ and all $t\geq 1$.
\end{thm} 

A natural question is whether the restriction on $\codim(\Sing(X))$ can be relaxed or even removed. The following example from \cite{BBLSZ3} shows that this is not the case.

\begin{exa}
\label{lci exa with linear vanishing}
Let $R=\KK[x,y,u,v,w]$ where $\KK$ is a field of characteristic 0. Fix an integer $c\geq 2$ and set $I:=(uy-vx,vy-wx)+(u,v,w)^c$. Then one can check that 
\begin{enumerate}
\item $X=\Proj(R/I)$ is local complete intersection in $\PP^4_{\KK}$;
\item $H^2_{\frakm}(R/I^t)_{-ct+1}\neq 0$;
\item $H^2_{\frakm}(R/I^t)_{\leq -ct}=0$.
\end{enumerate}
\end{exa}

Example \ref{lci exa with linear vanishing} indicates that, if one removes the restriction on the homological degree by $\codim(\Sing(X))$, the best vanishing result one can hope for is an asymptotic vanishing bounded by a linear function of $t$. Such an asymptotic vanishing turns out to be true, as shown in \cite{BBLSZ3}.

\begin{thm}
\label{linear vanishing for lci}
Let $X$ be a closed local complete intersection subscheme of $\PP^n$ over a field of arbitrary characteristic . Then there exists an integer $c\geq 0$ such that for each $t\geq 1$ and $i<\dim(X)$, one has
\[H^i(X_t,\calO_{X_t}(j))=0,\ \forall\ j<-ct.\]
\end{thm}

When $\Proj(R/I)$ is a local complete intersection (here $R=\KK[x_0,\dots,x_n]$ and $I$ is a homogeneous ideal of $R$), the local cohomology modules $H^j_\frakm(R/I^t)$ have finite length for $j<\dim(R/I)$ and consequently $H^j_\frakm(R/I^t)_{\ell \ll 0}=0$. This is one of the underlying reasons for the vanishing in Theorems \ref{Kodaira for lci thickenings} and \ref{linear vanishing for lci}. Once the local complete intersection assumption is dropped, $H^j_\frakm(R/I^t)$ may not have finite length and hence the vanishing may fail. However, since $H^j_\frakm(R/I^t)$ are Artinian (even when $j=\dim(R/I)$), the socles $\Hom_R(R/\frakm, H^j_\frakm(R/I^t))$ are finite dimensional and vanish in all sufficiently negative degrees. Therefore, one can ask:

\begin{qun}
Let $R=\KK[x_0,\dots,x_n]$ and $I$ be a homogeneous ideal of $R$. For each $j\geq 0$, does there exist an integer $c$ such that 
\[\Hom_R(R/\frakm, H^j_\frakm(R/I^t))_{\ell}=0\]
for all $t\geq 1$ and all $\ell<-ct$?
\end{qun}

For related questions and applications, we refer the interested reader to \cite{Zhang2020}.
\subsection{Annihilation of local cohomology}

We now turn to the question: what annihilates the local cohomology
module $H^j_I(M)$? 

If $R$ is a Noetherian regular ring of prime characteristic $p$, then Huneke and Koh proved in \cite{HunekeKoh-cofiniteness} that $\ann_R(H^j_I(R))\neq 0$ if and only if $H^j_I(R)=0$. The same conclusion for Noetherian regular rings of characteristic 0 was established implicitly in \cite{L-Dmods}. The aforementioned result due to Huneke--Koh was later generalized to strongly $F$-regular domains in \cite{BoixEghbali-Annihilators}. Inspired by the results due to Huneke--Koh and Lyubeznik, Lynch \cite{Lynch-Annihilators} conjectured that $\dim(A/\ann_A(H^{\delta}_I(A)))=\dim(A/H^0_I(A))$ for every Noetherian local ring $A$, where $\delta=\lcd_A(I)$. This conjecture turns out to be false in general, {\it cf.} \cite{Bahmanpour-CiA2017} and \cite{SinghWalther-CiA2020}. Note that the rings in the counterexamples in \cite{Bahmanpour-CiA2017} and \cite{SinghWalther-CiA2020} are not equidimensional. In \cite[Question 6]{Hochster-SurveyLC-CiA2020} Hochster asks the following.

\begin{qun}
\label{Hochster question faithful}
If $A$ is a Noetherian local domain and $I$ is an ideal of cohomological dimension $c$, is $H^c_I(A)$ a faithful $A$-module?
\end{qun}

In \cite{HochsterJeffries-MRL2019}, Hochster and Jeffries answer this question in the affirmative in the following cases:
\begin{itemize}
\item $\ch(A)=p>0$ and $c$ equals the \emph{arithmetic rank} of $I$,
  see Subsection \ref{subsec-ara} below;
\item $A$ is a pure subring of a regular ring containing a field.
\end{itemize}
In \cite{DattaSwitalaZhang}, Datta, Switala and Zhang answer Question
\ref{Hochster question faithful} in the negative by the following
(equidimensional) example.
\begin{exa}
\label{Reisner example}
Let $R=\ZZ_2[x_0,\dots,x_5]$ and let $I$ be the ideal of $R$ generated by the $10$ monomials
\[
\{x_0x_1x_2, x_0x_1x_3, x_0x_2x_4, x_0x_3x_5, x_0x_4x_5, x_1x_2x_5, x_1x_3x_4, x_1x_4x_5, x_2x_3x_4, x_2x_3x_5\}.
\]
Then $\cd(I)=4$, but $\ann_R(H^4_I(R))$ is the ideal generated
by $2\in R$.
\end{exa}

When $(A,\frakm)$ is a local ring, the annihilation of
$H^j_{\frakm}(A)$ is particularly interesting for $j<\dim(A)$, and has
a wide range of applications. We recall that an element $x\in
A^{\circ}$ is called a {\it uniform local cohomology annihilator} of
$A$ if $xH^j_{\frakm}(A)=0$ for $j<\dim(A)$, where
$A^{\circ}=A\setminus \bigcup_{\frakp\in\min(A)} \frakp$. Since
$H^j_{\frakm}(A)$ may not be finitely generated, it is not clear
whether such a uniform annihilator should exist. Surprisingly,
in \cite{ZhouUniformAnn-2007} Zhou proved that if $A$ is an excellent
local ring then $A$ admits a uniform local cohomology annihilator if
and only if $A$ is equidimensional. If $x$ is a uniform local
cohomology annihilator then $A_x$ is Cohen--Macaulay ({\it
cf.} \cite{HochsterHuneke-JAMS1990}, \cite{ZhouUniformAnn-2006}); in
fact, there is a deep connection between the existence of uniform
local cohomology annihilators and the Cohen--Macaulay locus. To explain
this connection, we need to recall some definitions from \cite{Huneke-Invent1992}. For a Noetherian ring $A$, a finite complex of finitely generated free $A$-modules 
\[G_{\bullet}: 0\to G_n\xrightarrow{f_n}G_{n-1}\to \cdots \to G_1\xrightarrow{f_1}G_0\]
is said 
\begin{itemize}
\item to satisfy the \emph{standard condition on rank} if
  $\rank(f_i)+\rank(f_{i-1})=\rank(G_{i-1})$ for $1\leq i\leq n$ and
  $\rank(f_n)=\rank(G_n)$, where the rank of a map is the determinantal rank;
\item to satisfy the \emph{standard condition on height} if
  $\hight(I(f_i))\geq i$ for all $i$, where $I(f_j)$ is the ideal generated by the rank-size minors of $f_j$ which is viewed as a matrix.
\end{itemize}

For a Noetherian ring $A$, we denote by $\CM(A)$ the set of elements $x\in A$ such that for all finite complexes $G_{\bullet}$ of finitely generated free $A$-modules
satisfying the standard conditions on rank and height, $xH_i(G_{\bullet})=0$ for $i\geq 1$. Huneke conjectured in \cite{Huneke-Invent1992} that if $A$ is an equidimensional excellent Noetherian ring then $\CM(A)$ is not contained in any minimal prime of $A$. Zhou proved this conjecture in \cite{ZhouUniformAnn-2007} by showing the following theorem.
\begin{thm}
Let $A$ be an excellent local ring $A$. Then $A$ admits a uniform local cohomology annihilator if and only if $\CM(A)$ is not contained in any minimal prime of $A$.
\end{thm}

One may consider the uniform annihilation of local cohomology in a different direction. 
\begin{qun}
\label{question uniform Frobenius ann}
Let $(A,\frakm)$ be a Noetherian local ring of characteristic $p$ and $I$ be an ideal of $A$. Does there exist a constant $B$ such that
\[\frakm^{Bp^e}H^0_{\frakm}(A/I^{[p^e]})=0\]
for all $e\geq 1$?
\end{qun}

The  special case of Question \ref{question uniform Frobenius ann} where $I$ is primary to a prime ideal of height $\dim(A)-1$ was explicitly asked by Hochster and Huneke in \cite{HochsterHuneke-JAMS1990}; a positive answer to this special case would have significant consequences in tight closure theory, especially to the notion of $F$-regularity. Question \ref{question uniform Frobenius ann} is wide open to the best of our knowledge. The graded analog, when $A$ is a standard graded ring over a field of characteristic $p$ and $I$ is homogeneous, has also attracted attention. When $\dim(A/I)=1$, the graded version was settled independently in \cite{HunekeSaturation-CiA2000} and \cite{Vraciu-JA2000}. Let $A$ be a standard graded ring over a field and let $M$ be a finitely generated graded $A$-module. We set
\[
a_j(M):=\max\{\ell\mid H^j_{\frakm}(M)_{\ell}\neq 0\}
\]
for each integer $j$. Since $H^j_{\frakm}(M)$ is Artinian,
$a_j(M)<\infty$ for each $j$. Hence for a homogeneous ideal $J$, if
$a_0(A/J)\leq t$, then $\frakm^{t+1}H^0_{\frakm}(A/J)=0$. Consequently, if
there is an integer $B$ such that $a_0(A/I^{[p^e]})\leq Bp^e$ for all
$e$, then $\frakm^{Bp^e+1}H^0_{\frakm}(A/I^{[p^e]})=0$ for all
$e$. In general, it is an open question whether there exists an integer
$B$ (independent of $e$) such that $a_0(R/I^{[p^e]})\leq Bp^e$ for all $e$. On the other hand,
$a_0(M)$ may be considered as a partial Castelnuovo--Mumford regularity since the regularity  of $M$ is defined as
\[
\reg(M):=\max\{a_j(M)+j\mid 0\leq j\leq \dim(M)\},
\]
so that $a_0(M)\leq \reg(M)$. Therefore, one may ask for a stronger conclusion on the linear growth of $\reg(A/I^{[p^e]})$ with respect to $p^e$. Indeed, the following was asked in \cite[p.~212]{KatzmanComplexity-JA1998}.

\begin{qun}
Let $A$ be a standard graded ring over a field of characteristic $p$ and $I$ be a homogeneous ideal. Does there exist a constant $C$ such that
\[\reg(A/I^{[p^e]})\leq Cp^e\]
for all $e$?
\end{qun}

Some progress has been made: for cases of small singular locus see \cite{ChardinRegularity-2004}, \cite{BrennerLinearBound-MMJ2005} and \cite{ZhangRegularity-2015}; for rings of finite Frobenius representation type, see \cite{KatzmanSchwedeSinghZhang}.

\medskip

At the crux of the homological conjectures stands the the
existence of big Cohen--Macaulay algebras: the  assertion that each Noetherian complete local domain $(A,\frakm)$ admits an algebra (not necessarily Noetherian) in which every system of parameters of $A$ becomes a regular sequence. A beautiful result of Hochster--Huneke in \cite{HochsterHuneke-AnnMath1992} says that if $A$ is an excellent Noetherian local domain of characteristic $p$ then its absolute integer closure\footnote{The absolute integral closure of an integral domain $A$ is defined to be the integral closure of $A$ in the algebraic closure of the field of fractions of $A$.} $A^+$ is a big Cohen--Macaulay $A$-algebra. In \cite{HunekeLyubeznik-AdvMath2007}, Huneke and Lyubeznik gave a much simpler proof using annihilation of local cohomology. 

The following result, proved in \cite[Lemma~2.2]{HunekeLyubeznik-AdvMath2007}, has been referred to as `the equational lemma'.

\begin{thm}
\label{thm: equational lemma}
Let $A$ be a commutative Noetherian domain that contains a field of characteristic $p$, let $\KK$ be its field of fractions and $\overline{\KK}$ be the algebraic closure of $\KK$. Let $I$ be an ideal of $A$ and let $\alpha$ be an element in $H^i_I(R)$ such that the elements\footnote{Here $\alpha^p$ denotes $f(\alpha)$ where $f$ is the natural action of Frobenius on $H^i_I(A)$ induced by the Frobenius endomorphism on $A$.} $\alpha, \alpha^p,\dots, \alpha^{p^t},\dots$ belong to a finitely generated submodule of $H^i_I(A)$. Then there is a module-finite extension $A'$ of $A$ inside $\overline{\KK}$ such that the natural map $H^i_I(A)\to H^i_I(A')$ induced by $A\to A'$ sends $\alpha$ to 0.
\end{thm}

Since the module-finite extension $A'$ is constructed using the equations satisfied by $\alpha$, Theorem \ref{thm: equational lemma} is referred in the literature as an ``equational lemma". Using Theorem \ref{thm: equational lemma}, Huneke and Lyubeznik proved the following in \cite[Theorem~2.1, Corollary~2.3]{HunekeLyubeznik-AdvMath2007}:

\begin{thm}
Let $(A,\frakm)$ be a commutative Noetherian domain that contains a field of characteristic $p$, let $\KK$ be its field of fractions and $\overline{\KK}$ be the algebraic closure of $\KK$. Assume, furthermore, that $A$ is a homomorphic image of a  Gorenstein local ring. For every module-finite extension $A'$ of $A$
 inside $\overline{\KK}$, there exists module-finite extension $A'\subseteq A''$ inside $\overline{\KK}$ such that the natural maps
 \[H^i_{\frakm}(A')\to H^i_{\frakm}(A'')\]
 are the zero map for each $i<\dim(A)$.
 
 In particular, 
 \begin{enumerate}
 \item $H^i_{\frakm}(A^+)=0$ for $i<\dim(A)$;
 \item every system of parameter of $A$ is a regular sequence on $A^+$.
 \end{enumerate}
 \end{thm}
 
 The Huneke--Lyubeznik equational lemma (or equivalently, the technique
 of annihilating local cohomology with finite extensions) in characteristic $p$ has found many applications, for instance \cite{Bhatt-ANT2012} and \cite{BlickleSchwedeTucker-AJM2015}. In equi-characteristic 0, such annihilation of local cohomology is not possible once the dimension is at least 3: every module-finite extension of a normal domain must split in equi-characteristic 0. The situation in mixed characteristic has long been a mystery. However, in a very surprising turn of events, Bhatt proved in \cite[Theorem 5.1]{BhattRPlusMixedChar} the following:
 
 \begin{thm}
 Let $(A,\frakm)$ be an excellent Noetherian local domain with mixed
 characteristic $(0,p)$ and let $A^+$ be an absolute integral closure of $A$. Then
 \begin{enumerate}
 \item $H^i_{\frakm}(A^+/pA^+)=0$ for $i<\dim(A/pA)$ and $H^i_{\frakm}(A^+)=0$ for $i<\dim(A)$.
 \item Every system of parameters of $A$ is a Koszul regular sequence\footnote{A sequence of elements $z_1\dots,z_t$ in a commutative ring $C$ is a Koszul regular sequence if $H_i(K_{\bullet}(C;z_1,\dots,z_t))=0$ for $i>0$ where $K_{\bullet}(C;z_1,\dots,z_t)$ is the Koszul complex of $C$ on $z_1,\dots,z_t$.} on $A^+$.
\item If $A$ admits a dualizing complex, then there exists a module-finite extension $A\to B$ with $H^i_{\frakm}(A/pA)\to H^i_{\frakm}(B/pB)$ being the 0 map for all $i<\dim(A/pA)$.
 \end{enumerate} 
 \end{thm}  

{For other connections between annihilators of local cohomology modules and homological conjectures, we refer the reader to \cite{RobertsApplicationsDualizingComplexes1976, SchenzelCohomologicalAnnihiltors1982}.}

\section{$D$- and $F$-structure}
\label{sec-DF}

In this section we discuss some special structures that
local cohomology have. In positive characteristic the Frobenius
endomorphism is the main tool, while in any case they have a
structure over the ring of differential operators. 

\subsection{$\calD$-modules}
Following Grothendieck's approach in \cite{EGA4-4}, we reproduce the definition of differential operators as follows.
Let $A$ be a commutative ring. The  \emph{differential operators}
\[
\calD(A)=\bigcup_{j\in\NN}\calD_j(A)
\]
on $A$ (which is to say, the differential operators from $A$ to $A$)
are classified by their \emph{order} $j$ (a natural number), and defined
inductively as follows. The differential operators $\calD_0(A)$ of
\emph{order zero} are precisely the multiplication maps $\tilde{a}\colon A\to
A$ where $a\in A$; for each positive integer $j$, the differential
operators $\calD_j(A)$ of \emph{order less than or equal to $j$} are those
additive maps $P\colon A\to A$ for which the commutator
\[
[\tilde{a},P]\ =\ \tilde{a}\circ P-P\circ\tilde{a}
\]
is a differential operator on $A$ of order less than or equal to
$j-1$. If $P'$ and $P''$ are differential operators of orders at most
$j'$ and $j''$ respectively, then $P'\circ P''$ is again a differential
operator and its  order is at most $j'+j''$. Thus, the differential
operators on $R$ form an $\NN$-filtered subring $\calD(R)$ of $\End_\ZZ(R)$,
and the order filtration is (by definition) increasing and exhaustive.

When $A$ is an algebra over the central subring $\kk$, we define
$\calD(A,\kk)$ to be the subring of $\calD(A)$ consisting of those
elements of $\calD(A)$ that are $\kk$-linear. Thus,
$\calD(A,\ZZ)=\calD(A)$ and $\calD(A,\kk)=\calD(A)\cap\End_\kk(A)$.
It turns out that if $A$ is an algebra over a perfect field $\FF$ of prime characteristic, then $\calD(A,\FF)=\calD(A)$, see, for example, \cite[Example~5.1~(c)]{LyubeznikFModules}.

By a $\calD(A,\kk)$-module, we mean a \emph{left} $\calD(A,\kk)$-module,
unless we expressly indicate a right module. The standard example of a
$\calD(A,\kk)$-module is $A$ itself. Using the quotient rule,
localizations $A'$ of $A$ also carry a natural
$\calD(A,\kk)$-structure and the formal quotient rule induces a natural map $\calD(A,\kk)\to \calD(A',\kk)$. Suppose $\fraka$ is an ideal of $A$. The \v Cech complex on a generating set for $\fraka$ is a complex of $\calD(A,\kk)$-modules; it then follows that each local cohomology module $H^k_{\fraka}(A)$ is a $\calD(A,\kk)$-module.

More generally, if $M$ is a $\calD(A,\kk)$-module, then each local
cohomology module $H^k_{\fraka}(M)$ is also a
$\calD(A,\kk)$-module. This was used by Kashiwara as early as 1970 as inductive
tool in algebraic analysis
via reduction of dimension \cite{Kashi-thesis} and was
introduced to commutative algebra in~\cite[Examples 2.1~(iv)]{L-Dmods}.)

\medskip

If $R$ is a polynomial or formal power series ring in the variables
$x_1,\dots,x_n$ over a commutative ring $\kk$, then
$\frac{1}{{t_i}!}\frac{\partial^{t_i}}{\partial x_i^{t_i}}$ can be
viewed as a differential operator on $R$ even if the integer $t_i!$ is
not invertible. In these cases, $\calD(R,\kk)$ is the free $R$-module
with basis elements
\[
\frac{1}{{t_1}!}\frac{\partial^{t_1}}{\partial x_1^{t_1}}\ \cdots\ \frac{1}{{t_n}!}\frac{\partial^{t_n}}{\partial x_n^{t_n}}
\qquad\text{ for } \ (t_1,\dots,t_n)\in\NN^n\,,
\]
see~\cite[Th\'eor\`eme~16.11.2]{EGA4-4}. When $R$ is a polynomial ring
or formal power series ring over a field $\kk$ of characteristic 0,
then the ring of differential operators 
\[
\calD(R,\kk)=R\langle \frac{\partial }{\partial
  x_1},\dots,\frac{\partial}{\partial x_n}\rangle,
\]
is known as the \emph{Weyl algebra}, a simple ring in the sense that
it has no non-trivial two-sided ideals.

If $\kk$ is a field and if $A$ is a singular $\kk$-algebra then the
structure of $\calD(R,A)$ can be very complicated, even in
characteristic zero. For example, the ring of differential operators
on the cone over an elliptic curve is not Noetherian and also not
generated by homogeneous operators of bounded finite degree,
\cite{BGG-Ell,Hsiao-big}. In most cases, differential operators on
singular spaces are completely mysterious, except for toric varieties, Stanley-Reisner rings and hyperplane arrangements, see \cite{Musson-tori,Jones-toric,MSTW,Traves-orbi, TrippDiffOperStanleyReisnerRings, HolmDiffOperHyperplaneArr}.

\subsubsection{Characteristic 0}

As references for background reading in this section we recommend
\cite{BjorkBook,Kashi-thesis,Kashi-book,KashiwaraShapira,HTT}. 

Let $\kk$ denote a field of characteristic 0 and fix $n\in\NN$. Let
$R$ denote either $\kk[x_1,\dots,x_n]$ or $\kk[[x_1,\dots,x_n]]$, and
let $\calD$ denote $\calD(R,\kk)$, unless specified otherwise. The
partial differential operator $\frac{\partial }{\partial x_i}$ is
denoted by $\partial_i$ for each variable $x_i$.

Note that here the order of $r\partial^{e_1}_1\cdots \partial^{e_n}_n$
($r \in R$) equals simply $\sum_ie_i$. The order (\emph{i.e.}, the
filtration level) of an element $\sum_{c_{\alpha,\beta}\neq 0}
c_{\alpha,\beta}\bsx^\alpha\bsdel^\beta\in \calD$ is the maximum
of the orders $|\beta|$ of its terms $\bsx^\alpha\bsdel^\beta$.
Then we have $\calD_j=\{r\partial^{e_1}_1\cdots \partial^{e_n}_n\mid
r\in R,~\sum_ie_i\leq j\}$, an increasing and exhaustive filtration of
$\calD$, called the {\it order filtration} of $\calD$.

Using the order filtration $\{\calD_j\}$, one can form the associated
graded ring,
\[
\gr(\calD):=\calD_0\oplus
\frac{\calD_{1}}{\calD_0}\oplus\cdots.
\]
Since the only nonzero commutators of pairs of generators in $\calD$
are the $[\del_i,x_i]=1\in\calD_0$, it follows that $\gr(\calD)$ is
isomorphic to a (commutative) ring of polynomials
$R[\xi_1\dots,\xi_n]$ where $\xi_i$ is the image of $\partial_i$ in
$\calD_1/\calD_0$. Note that $\gr(\calD)$ is naturally the coordinate
ring on the cotangent space of $\kk^n$, if $R$ is a
ring of polynomials. We use this to construct varieties from
$\calD$-modules as follows.

\begin{dfn}
Let $M$ be a $\calD$-module. A \emph{filtration} of $M$ \emph{with respect to the order filtration} $\{\calD_j\}$ is a sequence of $R$-submodules $\{F_iM\}$ such that
\begin{enumerate}
\item $F_0M\subseteq F_1M\subseteq \cdots \subseteq F_iM\subseteq F_{i+1}M\subseteq \cdots$;
\item $\bigcup_i F_iM=M$;
\item $\calD_j\cdot F_iM\subseteq F_{i+j}M$.
\end{enumerate}

Such filtration is called a \emph{good filtration} if the associated graded module $\gr^{F}(M):=F_0M\oplus\frac{F_{1}M}{F_0M}\oplus \cdots$ is finitely generated over $\gr(\calD)$.
\end{dfn}

Every finitely generated $\calD$-module admits a good filtration $\{F_iM\}$; for instance, if $M$ can be generated by $m_1,\dots,m_d$, then setting $F_iM:=\sum_j\calD_im_j$ produces a good filtration of $M$. Set $J$ to be the radical of $\ann_{\gr(\calD)}(\gr^{F}(M))$. This ideal $J$ is independent of the good filtration $\{F_iM\}$ (\emph{cf.} \cite[1.3.4]{BjorkBook}, \cite[11.1]{CoutinhoBook}), and is called the {\it characteristic ideal} of $M$. The characteristic ideal of $M$ induces the notion of dimension of $M$ (as a $\calD$-module) and characteristic variety of $M$.

\begin{dfn}
Let $M$ be a $\calD$-module with good filtration and let $J$ be its
characteristic ideal. The dimension of $M$ is defined as
\[d(M):=\dim(\gr(\calD)/J).\]
The \emph{characteristic variety} $\Ch(M)$ of $M$ is defined as the
subvariety of $\Spec(\gr(\calD))$ defined by $J$. The set of the
irreducible components of $\Ch(M)$, paired with their multiplicities
in $\gr(M)$  is called the \emph{characteristic cycle} of $M$.
\end{dfn}

It turns out that dimensions cannot be small:

\begin{thm}[Bernstein Inequality]
Let $M$ be a nonzero finitely generated $\calD$-module. Then
\[n\leq d(M)\leq 2n.\]
\end{thm}
The nonzero modules of minimal dimension form a category with many
good features.
\begin{dfn}
A finitely generated $\calD$-module $M$ is called \emph{holonomic} if $d(M)=n$ or $M=0$.
\end{dfn}

\begin{exa}
  \begin{asparaenum}
    \item
Set $F_iR=R$ for all $i\in\NN$. Then one can check that $\{F_iR\}$ is a good filtration on $R$ and $\gr^F(R)\cong R$. Hence 
\[
J=\sqrt{\ann_{\gr(\calD)}(\gr^F(R))}=\ideal{\xi_1,\dots,\xi_n}.
\] 
This shows that $d(R)=n$. Therefore, $R$ is a holonomic $\calD$-module.
\item
Denote $H^n_{\frakm}(R)$ by $E$ and set $\eta=\left[\frac{1}{x_1\cdots
    x_n}\right]$, the class of the given fraction inside $E$. Set
$F_iE=\calD_i\cdot \eta$. Then one can check that $\{F_iE\}$ is a good
filtration of $E$ and $\gr^F(E)\cong \kk[\xi_1,\dots,\xi_n]$ where
$\xi_i$ denotes the image of $\partial_i$ in $\calD_1/\calD_0$. Hence
\[
J=\sqrt{\ann_{\gr(\calD)}(\gr^F(E))}=\ideal{x_1,\dots,x_n}.
\] 
This shows that $d(E)=n$. Therefore, $E=H^n_{\frakm}(R)$ is a
holonomic $\calD$-module.
\end{asparaenum}
\end{exa}

We collect next some of the basic properties of holonomic $\calD$-modules.
\begin{thm}
\label{thm: properties of holonomic}
\begin{enumerate}
\item Holonomic $\calD$-modules form an Abelian subcategory of the category of $\calD$-modules that is closed under the formation of submodules, quotient modules and extensions (\cite[1.5.2]{BjorkBook}). 
\item If $M$ is holonomic, then so is the localization $M_f$ for every
  $f\in R$ (\cite[3.4.1]{BjorkBook}). Consequently, each local
  cohomology module $H^j_I(M)$ of $M$ is holonomic.
\item Each holonomic $\calD$-module admits a finite filtration in the
  category of $\calD$-modules in which each composition factor is a simple
  $\calD$-module (\cite[2.7.13]{BjorkBook}).
\item A simple holonomic $\calD$-module has only one associated prime
  (\cite[3.3.16]{BjorkBook}).
\end{enumerate}
\end{thm}

Certain finiteness properties of $H^j_I(R)$ are enjoyed by arbitrary
holonomic $\calD$-modules. In the following list, the first is a
special case of Kashiwara equivalence; the latter were established in
\cite[Theorem 2.4]{L-Dmods}.

\begin{thm}
Let $R=\kk[[x_1,\dots,x_n]]$ and let $\frakm$ denote the maximal
ideal. Let $M$ be a finitely generated $\calD$-module.
\begin{enumerate}
\item If $\dim(\Supp_R(M))=0$, then $M$ is a direct sum of copies of $\calD/\calD \frakm$.
\item $\injdim_R(M)\leq \dim(\Supp_R(M))$.
\item If $M$ is finitely generated (as a $\calD$-module), then $M$ has finitely many associated primes (as an $R$-module).
\item If $M$ is holonomic, then the Bass numbers of $M$ are finite.
\end{enumerate}
\end{thm}

Similar statements hold when $R=\kk[x_1,\ldots,x_n]$.

\begin{rmk}
Let $S=\kk[y_1,\dots,y_{2n}]$ be the polynomial ring over $\kk$ in $2n$ variables. When $R=\kk[x_1,\dots,x_n]$, we have seen that $\gr(\calD)\cong S$. The Poisson bracket on $S$ is defined as follows:
\[\{f,g\}=\sum_{i=1}^n(\frac{\partial f}{\partial y_{n+i}}\frac{\partial g}{\partial y_i}-\frac{\partial g}{\partial y_{n+i}}\frac{\partial f}{\partial y_i}).\]
An ideal $\fraka$ of $S$ is said to be {\it closed under the Poisson bracket} if $\{f,g\}\in \fraka$ whenever $f,g\in \fraka$.

The Poisson bracket is closely related to symplectic structures on $\CC^{2n}$ and involutive subvarieties of $\CC^{2n}$. A \emph{symplectic structure} $\omega$ on $\CC^{2n}$ is a non-degenerate skew-symmetric form; the standard one is given by
\[\begin{bmatrix}0& -I_n\\I_n&0 \end{bmatrix}\]
where $I_n$ is the $n\times n$ identity matrix. Fix a symplectic structure $\omega$ on $\CC^{2n}$. Given any subspace $W$ of $\CC^{2n}$, its skew-orthogonal complement is defined as 
\[W^{\perp}:=\{\vec{v}\in \CC^{2n}\mid
\omega(\vec{w},\vec{v})=0\ \forall\ \vec{w}\in W\}.\]

A subspace $W$ is called {\it involutive} if $W^{\perp}\subseteq W$. A subvariety $X$ of $\CC^{2n}$ is called {\it involutive} if the tangent space $T_{x}X\subseteq \CC^{2n}$ is a involutive subspace for every smooth point $x\in X$. One can show that an affine variety $X\subseteq \CC^{2n}$ is involutive with respect to the standard symplectic structure on $\CC^{2n}$ if and only if its (radical) defining ideal $I(X)$ is closed under the Poisson bracket. 
\end{rmk}

The following was conjectured in \cite{GuilleminQuillenSternberg1970} by Guillemin--Quillen--Sternberg and proved in \cite{KashiwaraKawaiSato1973} for sheaves of differential operators with holomorphic coefficients on a complex analytic manifold by Kashiwara--Kawai--Sato. The first algebraic proof was discovered by Gabber in \cite{GabberIntegrability1981}.
\begin{thm}
Let $R=\kk[x_1,\dots,x_n]$ and $M$ be a holonomic $\calD$-module. Then the characteristic ideal $J$ of $M$ is closed under the Poisson bracket on $gr(\calD)$.
\end{thm}

Again, let $R$ be either $\kk[x_1,\dots,x_n]$ or
$\kk[[x_1,\dots,x_n]]$. Then each
$\calD$-module $M$ admits a (global) \emph{de Rham complex}.  This is a complex
of length $n$, denoted $ \Omega_R^{\bullet}\otimes M$ (or simply
$\Omega_R^{\bullet}$ in the case $M = R$), whose objects are
$R$-modules but whose differentials are merely $\kk$-linear.  It is
defined as follows \cite[\S 1.6]{BjorkBook}: for $0 \leq i \leq n$,
$\Omega^i_R\otimes M$ is a direct sum of $n \choose i$ copies of $M$,
indexed by $i$-tuples $1 \leq j_1 < \cdots < j_i \leq n$.  The summand
corresponding to such an $i$-tuple will be written $M \, \de x_{j_1}
\wedge \cdots \wedge \de x_{j_i}$. The $\kk$-linear differentials
$\de^i: \Omega_R^i \otimes M\rightarrow \Omega_R^{i+1}\otimes M$ are defined by 
\[
\de^i(m \,\de x_{j_1} \wedge \cdots \wedge \de x_{j_i}) = \sum_{s=1}^n
\partial_s(m)\, \de x_s \wedge \de x_{j_1} \wedge \cdots \wedge \de x_{j_i},
\]
with the usual exterior algebra conventions for rearranging the wedge
terms, and extended by linearity to the direct sum. We remark that in
the polynomial case we are simply using the usual K\"{a}hler
differentials to build this complex, whereas in the formal power
series case, we are using the $\mathfrak{m}$-adically continuous
differentials (since in this case the usual module $\Omega^1_{R/\kk}$ of
K\"{a}hler differentials is not finitely generated over $R$).
An alternative way is to view $\Omega^\bullet_R\otimes M$ as a
representative of $\omega_R\otimes^L_\calD M$, where $\omega_R$ is the
\emph{right} $\calD$-module $\calD/\ideal{\del_1,\ldots,\del_n}\calD$
which is as $R$-module simply $R$.

The cohomology objects $H^i(M \otimes \Omega_R^{\bullet})$ are
$\kk$-spaces and called the \emph{de Rham cohomology spaces} of the left $\calD$-module $M$, and are denoted $H^i_{\dR}(M)$. The simplest de Rham cohomology spaces (the $0$th and $n$th) of $M$ take the form
\begin{align*}
H^0_{\dR}(M) &= \{m \in M \mid \partial_1(m) = \cdots = \partial_n(m) = 0\} \subseteq M\\
 H^n_{\dR}(M) &= M/(\partial_1\cdot(M) + \cdots + \partial_n\cdot(M)).
\end{align*}

The de Rham cohomology spaces are not finite dimensional in
general, even for finitely generated $M$. The following theorem is (for the Weyl algebra) a special
case of fact that the $\calD$-module theoretic direct image functor
preserves holonomicity, \cite[Section 3.2]{HTT}. It can be found in
\cite[1.6.1]{BjorkBook}) for the polynomial case and in
\cite[Prop. 2.2]{vdEcoker1985} for the formal power series case.

\begin{thm}
\label{dR finite dim for holonomic}
Let $M$ be a holonomic $\calD$-module. The de Rham cohomology spaces $H^i_{\dR}(M)$ are finite-dimensional over $\kk$ for all $i$.
\end{thm} 

Let $E$ denote $H^n_{\frakm}(R)$. If $R=\kk[[x_1,\dots,x_n]]$, then we
use $D(-)$ to denote $\Hom_R(-,E)$ (this is the Matlis dual; it should not be
confused with the holonomic duality functor $\DD$ which is quite
different). If $R=\kk[x_1,\dots,x_n]$, we consider the following
``natural'' grading on $R$ and on $\calD$: 
\[
\deg(x_i)=1,\ \deg(\partial_i)=-1,\ i=1,\dots,n.
\]
Note that this is really a grading on $\calD$ since the relations
$[\del_i,x_i]=1$ are homogeneous of degree zero.
Then $E$ inherits a grading from setting $\deg([\frac{1}{x_1\dots
    x_n}])=-n$. In this graded setting, we use $\sideset{^*}{_R}\Hom$
to denote the graded $\Hom$ and use $D^*(-)$ to denote
$\sideset{^*}{_R}\Hom(-,E)$ (the graded Matlis dual).

It turns out that $D(-)$ is a functor on the category of
$\calD$-modules, that is compatible with de Rham cohomology. The following theorem is a combination of \cite[Theorem 5.1]{SwitalaCompos2017} and \cite[Theorem A]{SwitalaZhangAdvMath2018}.
\begin{thm}
\label{thm: SwitalaZhang}
\begin{enumerate}
\item Let $R=\kk[[x_1,\dots,x_n]]$ and $M$ be a holonomic $\calD$-module. Then
\[H^{i}_{\dR}(M)^\vee \cong H^{n-i}_{\dR}(D(M)),\ i=1,\dots,n,\]
where $(-)^\vee$ denotes the $\kk$-dual of a $\kk$-vector space.
\item Let $R=\kk[x_1,\dots,x_n]$ and $M$ be a graded $\calD$-module. Assume that $\dim_{\kk}(H^i_{\dR}(M))<\infty$. Then
\[(H^{i}_{\dR}(M))^\vee \cong H^{n-i}_{\dR}(D^*(M)).\]
\end{enumerate}
\end{thm}

As shown in \cite[Example 3.14]{SwitalaZhangAdvMath2018}, $D(M)$ may not be holonomic even if $M$ is. The duality statements in Theorem \ref{thm: SwitalaZhang} show that the (graded) Matlis duals of holonomic $\calD$-modules still have finite dimensional de Rham cohomology.

\begin{rmk}\label{rmk-GreenleesMay}
  The idea of applying Matlis duality to local cohomology modules
  already appears in the work of Ogus and Hartshorne. For example,
  Proposition 2.2 in \cite{Ogus-LCDAV} states that in a local
  Gorenstein ring $A$ with dualizing functor $D(-)$, the dual
  $D(H^i_I(A))$ of the local cohomology module $H^i_I(A)$ is equal
  to the local cohomology module
  $H^{\dim(A)-i}_\frakP(\frakX,\calO_{\frakX})$ where $\frakX$ is the
  completion of $\Spec(A)$ along $I$, and $\frakP$ its closed point.

  In much greater generality, \emph{Greenlees--May duality}
  \cite{GreenleesMay} states that (the derived functor of sections with
  support in $I$) $R\Gamma_I(-)$ and (the derived functor of completion along $I$) 
  $L\Lambda^I(-)$ are adjoint functors. See also
  \cite{TarrioJeremiasLipman, Lipman-Guanajuato}.
\end{rmk}

\smallskip

We briefly discuss algorithmic aspects.
The Weyl algebra is both left and right Noetherian and has a
Poincar\'e--Birkhoff--Witt basis of a polynomial ring in $2n$
variables; this makes it possible to
extend the usual Gr\"obner basis techniques to
$\calD$-modules, see for example \cite{Galligo-Linz85}.

When $R$ is a polynomial
  ring over the rational numbers, algorithms have been formulated that
  compute: 
\begin{enumerate}
\item\label{item-alg1} the local cohomology modules $H^i_I(R)$ in
  \cite{W-lcD}, but see also \cite{Oaku-Duke,OT2,BerkeschLeykin};
\item\label{item-alg2} the characteristic cycles and Bass numbers of $H^j_I(R)$ when $I$ is a monomial ideal in \cite{MontanerJPPA2000,MontanerJPPA2004};
\item\label{item-alg3} an algorithm to compute the support of local
  cohomology modules in \cite{LeykinMontaner2006}.
\end{enumerate} 
In a nutshell, the algorithms are based on the fact that the modules
that appear in a \v Cech complex $\check
C^{\bullet}(R;f_1,\ldots,f_m)$ are holonomic and sums of modules
generated by fractions of the form $(f_{i_1}\cdots f_{i_t})^{e}$ for
sufficiently small $e\in\ZZ$. In general, $e=-n$ is sufficient by
\cite{Saito-ASPM09}, but in the spirit of computability, it is
desirable to know the largest $e$ that may be used. This number turns
out to be the smallest integer root of the \emph{Bernstein--Sato
  polynomial} $b_f(s)$ of the polynomial $f$ in question. Indeed, as
was shown by Bernstein in \cite{Bernstein-bfu}, for every polynomial
$f\in R$ there is a linear differential operator $P$ depending
polynomially on the additional variable $s$ such that
\[
P(x_1,\ldots,x_n,\del_1,\ldots,\del_n,s)\bullet f^{s+1}=b_{P,f}(s)\cdot f^s,
\]
where $0\neq b_{P,f}(s)\in\kk[s]$ with $\kk$ a field of definition for
$f$.  Since $\kk[s]$ is a PID, Bernstein's theorem implies there is a
monic generator for the ideal of all $b_{P,f}(s)$ that arise this way;
this then is called the \emph{Bernstein--Sato polynomial} $b_f(s)$. It was shown
to factor over the rational numbers in
\cite{Malgrange-isolee,Kashiwara-bfu} and is a fascinating invariant
of $f$ as it relates to monodromy of the Milnor fiber, multiplier
ideals, (Igusa, topological, motivic) zeta functions, the
log-canonical threshold and various other geometric notions with
differential background. See \cite{Kollar-Proc97, Walther-survey} for
more details and \cite{AlvarezMontanerHunekeBetancourt} for a generalization of Bernstein-Sato polynomials to direct summands of polynomial rings. 

The polynomial $b_f(s)$ can be computed as the intersection of a left
ideal (derived from $f_1,\ldots,f_k$) inside a Weyl algebra with one
more variable $t$, with a ``diagonal subring'' $\QQ[t\del_t]$. The
idea of how to compute this intersection, and then to give a presentation
for the corresponding localization $R_f$, is due to Oaku. In
\cite{W-lcD} it was realized how to read off the $D$-structure of the
resulting local cohomology $H^1_f(R)$ and the process was scaled up to
non-principal ideals. The algorithm in \cite{OT2} is different in
nature and exploits the fact that local cohomology can be seen as
certain Tor-modules along the geometric diagonal in $2n$-space. It is,
however, still based on the computation of certain $b$-functions that
generalize the notion of a Bernstein--Sato polynomial.  To understand
conceptually how exactly the singularity structure of $I$ influences
the structure of the $\calD$-module $H^k_I(R)$ remains a question of
great interest.

\subsubsection{$\calD$-modules and group actions}
\label{subsubsec-inv}

We start with discussing the ideal determining the space of
matrices of bounded rank, and then outline more recent developments
that consider more general actions by Lie groups.

Let for now $\KK$ be a field, choose natural numbers $m\le n$ and set $R = \KK
[x_{ij}\mid 1\le i\le m,1\le j\le n]$. Let $I_{m,n,t}$ be the ideal
generated by the $t$-minors of the matrix $(x_{ij})$. Then $R/I$ is
Cohen--Macaulay and $I$ has height $(m - t + 1)(n - t + 1)$, compare
\cite{Bruns-Berkeley87,BrunsSchwanzl-BLMS90}.

Thus, in characteristic $p>0$ one has vanishing $H^k_{I_{m,n,t}}(R)$ for any
$k\neq (m - t + 1)(n - t + 1)$, because of the Frobenius (via the
Peskine-Szpiro vanishing result in Subsection \ref{subsec-vanishing}). In
characteristic zero, by \cite{BrunsSchwanzl-BLMS90}, $\lcd_R( I ) = mn
- t^2 + 1$. Therefore, $\lcd_R(I) - \depth(I,R) = (m + n - 2t)(t -
1)>0$, unless $m = n = t$ or $t = 1$.

Bruns and Schw\"anzl also proved in all characteristics that a
determinantal variety is cut out set-theoretically by $mn-t^2+1$
equations, and no fewer. In fact, these equations can be chosen to be
homogeneous; their methods rest on results involving \'etale
cohomology.  In particular, $I_{m,n,t}$ is a set-theoretic complete
intersection if and only if $n = m = t$. The same questions for the
case of symmetric and skew-symmetric matrices were answered completely
in
\cite{Barile-JA95} by Barile.  In many but not all cases the number of
defining equations agree  with the
local cohomological dimension.

Consider now the integral version of $I_{m,n,p}$ inside
$R_\ZZ=\ZZ[x_{ij}\mid 1\le i\le m,1\le j\le n]$. By \cite{LSW},
$H^k_{I_{m.n,t}}(R_\ZZ)$ is a vector space over $\QQ$ when $k$ exceeds
the height of $I_{m,n,t}$. Similar results are shown for the case of
generic matrices that are symmetric or anti-symmetric. As a corollary,
$H^{mn-t^2+1}_\fraka(A)$ vanishes for \emph{every} commutative ring
$A$ of dimension less than $mn$ where $\fraka$ is the ideal of
$t$-minors of \emph{any} $m\times n$ matrix over $A$. The initial version of
this result ($m=2=n-1=t$) appeared in \cite{HunekeKatzMarley}.

If $\KK$ is algebraically closed, Barile and Macchia study in
\cite{BarileMacchia-MR3957102} the number of elements needed  to
generate the ideal of $t$-minors
of a matrix $X$ up to radical, if the entries of $X$ outside some fixed
$t\times t$-submatrix are algebraically dependent over $\KK$. They
prove that this number drops at
least by one with respect to the generic case; under suitable
assumptions, it drops at least by $k$ if $X$ has $k$ zero
entries. 

\begin{ntn}
We now specialize the base field to $\CC$ and let $G$ be a connected linear 
algebraic group acting on a smooth connected complex algebraic
variety $X$. 
\end{ntn}

Suppose $R=\CC[x_1,\ldots,x_n]$ and $G$ is an algebraic Lie group
acting algebraically on $X=\CC^n$. There is a natural map
\[
\psi\colon \frakg \to
Der(\CC^n)
\]
from the Lie
algebra to the global vector fields on $\CC^n$,
\emph{i.e.}, the derivations inside the Weyl algebra
$\calD=\calD(R,\CC)$. 

The induced action $\star$ of $G$ on $R$ can be extended to an action on
$\calD$ that we also denote by $\star$. 
If $M$ is a $\calD$-module with a $G$-action, it is
\emph{equivariant} if the actions of $G$ on $\calD$ and $M$ are
compatible:
\[
(g\star P)\bullet (g\star m)=g\star(P\bullet m)
\]
for all $g\in G$, $P\in\calD$, $m\in M$.

Differentiating the $G$-action on $M$ one obtains an action of
$\frakg$ on $M$. One can now ask whether the Lie algebra element
$\gamma$ acts on $M$ via differentiation of the $G$-action the same
way that $\psi(\gamma)$ acts on $M$ as element of $\calD$. This is
not necessarily the case.
\begin{exa}
  Let $G=\CC^*$ act on $\CC$ by standard multiplication. The Lie
  algebra $\Lie(G)$ has an equivariant generator $\gamma$ that via $\psi$
  becomes $x\del_x\in\calD$.
  
  Let $M=\calD/\ideal{x\del_x-\lambda}$, with $G$-action inherited
  from the standard $G$-action on $\calD$: $g\star x=g^{-1}x$,
  $g\star \del_x=g\del_x$. Since $x\del_x-\lambda$ is $g$-invariant,
  this is indeed a $G$-action on $M$. Since $1\in \calD$ is $G$-invariant, the
  effect of $\gamma$ on $\bar 1\in M$ should be zero. On the other
  hand, $\psi(\gamma)\cdot\bar 1=\bar\lambda$. Thus, the two actions
  agree if and only if $\lambda=0$.
  
  Now note that there are other ways to act with $G$ on $M$. Indeed,
  a $\CC^*$-action is the same as the choice of a $\ZZ$-grading on
  $M$. Our choice above was $\deg(\bar 1)=0$; we now consider the
  choice $\deg(\bar 1)=k\in\ZZ$. This corresponds to $g\star \bar
  1=g^k\bar 1$, so that $\gamma$ must act on $\bar 1$ as
  multiplication by $k$. We conclude that the two actions of
  $\gamma$ agree if and only if $\lambda$ is an integer and the
  degree of $\bar 1$ is $\lambda$.
\end{exa}

\begin{dfn}
  The $\calD$-module $M$ is \emph{strongly equivariant} if the
  differential action of $G$ on $M$ agrees with the effect of $\psi$
  on $M$. In other words, $\gamma\star m=\psi(\gamma)m$ for
  all $\gamma\in\frakg,m\in M$.
\end{dfn}
\begin{rmk}
  Strong $G$-equivariance of a group acting on a variety $X$ can be
  also phrased as follows, see \cite[Dfn.\ 11.5.2]{HTT}: let $\pi$ and
  $\mu$ be the projection and multiplication maps
\begin{eqnarray*}
\pi \colon G \times X&\to & X,\\
\mu \colon G \times X&\to & X,
\end{eqnarray*}
respectively. Then $M$ is strongly equivariant if there is a
$\calD_{G\times X}$-isomorphism 
\[
\tau\colon \pi^*M\to \mu^*M
\]
that satisfies the usual compatibility conditions on $G\times G\times
X$, see \cite[Prop.\ 2.6]{vanDenBergh-notes}. If such $\tau$ exists,
it is unique.
\end{rmk}

Strongly $G$-equivariant $\calD_X$-modules are rather special
$D$-modules.  A $G$-equivariant morphism of smooth varieties with $G$-action
 automatically preserves $G$-equivariance under direct
and inverse images (since $G$ is connected, see
\cite[before Prop.\ 3.1.2]{vanDenBergh-Adv99}). If $G$ has finitely many
orbits on $X$, strong equivariance implies that the underlying
$\calD$-module is \emph{regular holonomic}; this is a growth condition
of the solution sheaf of the module and a critical component of the
Riemann--Hilbert correspondence. In this case, the simple and strongly
equivariant $\calD_X$-modules are labeled by pairs consisting of a
$G$-orbit $G/H$ and a finite-dimensional irreducible representation of
the component group of $H$ (in other words, a simple $G$-equivariant
local system on the orbit), \cite[Prop.\ 11.6.1]{HTT}. For example, if
$(\CC^*)^n$ acts on $\CC^n$, these simple modules are the modules
$H^{|S|}_{I_S}(R)$ where $R=\CC[x_1,\ldots,x_n]$, $S\in 2^{[n]}$ and
$I_S=\ideal{\{x_s\mid s\in S\}}$.

If $I$ is an ideal of $R$ and $Y$ the corresponding variety, then $I$
is $G$-stable if and only if $Y$ is. In this case, the localization
of a strongly equivariant module $M$ at an equivariant $g\in R$ is
also strongly equivariant. It follows that all local cohomology
modules $H^i_I(M)$ are as well. In particular, this holds when $M$ is $R$
or a local cohomology module of $R$ obtained in this way.

In \cite{RaicuWeymanWitt-Adv14}, the authors initiated the study of the
$GL$-equivariant decomposition of the local cohomology modules of determinantal ideals in
characteristic zero. The main result of the paper is a complete and
explicit description of the character of this representation. An
important consequence is a complete and explicit description of exactly 
which local cohomology modules $H^j_{I_{m,n,t}}(R)$ vanish and which
do not, in the case $t=n$. This was then refined and extended to
Pfaffians in \cite{RaicuWeymanWitt-Adv14}. The restriction $t=n$ was
removed in \cite{RaicuWeyman-ANT14}. Generalizations to symmetric and
skew-symmetric matrices were published in \cite{RaicuWeyman-JLMS16,Nang-JA12}.

In \cite{Raicu-Compositio16}, Raicu obtains results on the structure
of the $G$-invariant simple $D$-modules and their characters for
rank-preserving actions on matrices, extending work of Nang
\cite{Nang-Adv08,Nang-JA12}. Remarkably, for the case of symmetric
matrices, this provides a correction to a conjecture of Levasseur. Raicu's
methods produce composition factors for certain local cohomology
modules. In \cite{LoerinczRaicu-lambda} then this was taken the
furthest, to give character formul\ae\ for iterated local cohomology modules.

A more general approach was
used in
\cite{LoerinczRaicu-lambda,LorinczWalther-equi} in order
to study decompositions and categories of equivariant modules in the category of
$D$-modules, specifically with regards to quivers. These arise when
 when $G$ acts on $X$ with finitely many orbits and more particularly
 when $X$ is a spherical vector space and $G$ is reductive and
 connected. This leads to the study of the ``representation type'' of
 the underlying quiver (shown to be finite or tame) and the quivers are
 described explicitly for all irreducible $G$-spherical vector spaces
 of connected reductive groups using the classification of Kac. An
 early paper on this regarding the determinantal case was \cite{Nang-Japan04}.
More recently, cases of  exceptional representations and their
quivers have been studied:
\cite{Perlmann-JA20, LoerinczRaicuWeyman-CiA19}.

\begin{rmk}
  Invariant theory has also recently been aimed at singularity
invariants such as multiplier and test ideals 
\cite{HenriquesVarbaro16}, and 
$F$-pure thresholds \cite{MillerSinghVarbaro14}.
\end{rmk}

\subsubsection{Coefficient fields of arbitrary characteristic}

Let here $\kk$ be a field and set $R=\kk[x_1,\dots,x_n]$ or
$R=\kk[[x_1,\dots,x_n]]$. We have seen that $\calD=\calD(R,\kk)$ is
the free $R$-module with basis
\[
\frac{1}{{t_1}!}\frac{\partial^{t_1}}{\partial x_1^{t_1}}\ \cdots\ \frac{1}{{t_n}!}\frac{\partial^{t_n}}{\partial x_n^{t_n}}
\qquad\text{ for } \ (t_1,\dots,t_n)\in\NN^n\,
\]
When $\ch(\kk)=p>0$, the ring $\calD$ is no longer left or right
Noetherian. However, some desirable properties of $\calD$-modules in
characteristic 0 extend to the finite characteristic case. 

\begin{thm} Let $R=\kk[x_1,\dots,x_n]$ or $R=\kk[[x_1,\dots,x_n]]$,
  where $\kk$ is a field. Then
\begin{enumerate}
\item $\injdim_R(M)\leq \dim(\Supp_R(M))$ for every $\calD$-module $M$ (\cite{LyubeznikInjDim2000}).
\item $R_f$ has finite length in the category of $\calD$-modules for each $f\in R$ (\cite{LyubeznikFinitenessCharFree2000}, \cite{LyubeznikJPPA2011}).

Consequently, local cohomology modules $H^j_I(R)$ have finite length in the category of $\calD$-modules.
\end{enumerate}
\end{thm}

In general, it is a difficult problem to calculate the length
$\ell_{\calD}(H^j_I(R))$, or even just $\ell_{\calD}(R_f)$. Some
results are in \cite{Torrelli-intHom} and \cite{Bitoun-IMRN20}. The
following upper bounds were obtained in
\cite{KatzmanMaSmirnovZhang2018}.
\begin{thm}
Let $\kk$ be a field and $R=\kk[x_1,\dots,x_n]$. 
\begin{enumerate}
\item For each $f\in R$, 
\[\ell_{\calD}(R_f)\leq (\deg(f)+1)^n.\]
\item Assume an ideal $I$ can be generated by $f_1,\dots,f_t$. Then
\[\ell_{\calD}(H^j_I(R))\leq \sum_{1\leq i_1\cdots\leq i_j\leq t}(\deg(f_{i_1})+\cdots +\deg(f_{i_j})+1)^n-1.\]
\end{enumerate}
\end{thm}

\begin{exa}
\label{example: length Fermat cubic}
Let $R=\kk[x_1,x_2,x_3]$ and $f=x^3_1+x^3_2+x^3_3$. Then
\[
\ell_{\calD}(H^1_{(f)}(R))=\begin{cases}
1 \text{ if }  \ch(\kk)\equiv 2\ (\textrm{mod}\ 3);\\
1 \text{ if }  \ch(\kk)=3;\\
2 \text{ if } \ch(\kk)\equiv 1\ (\textrm{mod}\ 3);\\
2 \text{ if } \ch(\kk)=0.
\end{cases}
\]
\end{exa}

If $\ch(\kk)=0$ and $R=\kk[x_1,\dots,x_n]$ or
$R=\kk[[x_1,\dots,x_n]]$, then $R_f$ can always be generated by
$1/f^n$ as a $\calD$-module, but may not be generated by $1/f$.  For
instance, let $f,R$ be as in Example \ref{example: length Fermat
  cubic}, then $1/f$ generates a proper $\calD$-submodule of $R_f$ in
characteristic 0. On the other hand, in characteristic $p$, the
situation is quite different as shown in
\cite{AlvarezBlickleLyubeznik}, and generalized to rings of $F$-finite
representation type in \cite{TakagiTakahashi-MRL08}.

\begin{thm}
Let $\kk$ be a field of characteristic $p>0$ and let $R=\kk[x_1,\dots,x_n]$ or $R=\kk[[x_1,\dots,x_n]]$. Then $R_f$ can be generated by $1/f$ as a $\calD$-module for every $f\in R$.
\end{thm}

We have seen that, when $\kk$ is a field, $R=\kk[x_1,\dots,x_n]$, and
$M$ is a $\calD$-module, then $\injdim_R(M)\leq \dim(\Supp_R(M))$. Thus, if $\dim(\Supp_R(M))=0$, then $M$ must be an injective $R$-module. Let $I$ be a homogeneous ideal of $R$ and assume that $\Supp_R(H^j_I(R))=\{\frakm\}$ where $\frakm=(x_1.\dots,x_n)$. Then $H^j_I(R)\cong \oplus H^n_{\frakm}(R)^{\mu_j}$, a direct sum of finitely many copies of $H^n_{\frakm}(R)$. Since both $H^j_I(R)$ and $H^n_{\frakm}(R)$ are graded, a natural question is whether this isomorphism is degree-preserving. To answer this question, the notion of Eulerian graded $\calD$-modules was introduced in \cite{MaZhangEulerian}.

Recall that $R=\kk[x_1,\dots,x_n]$ and $\calD=\calD(R,\kk)$ are naturally graded via: 
\[\deg(x_i)=1,\qquad \deg(\partial_i)=-1.\]

\begin{dfn}
Denote the operator $\frac{1}{t_i!}\frac{\partial^{t_i}}{\partial
  x^{t_i}_i}$ by $\partial^{[t_i]}_i$.

  The \emph{$t$-th Euler operator} $E_t$ is defined as
\[E_t:=\sum_{{t_1+t_2+\cdots+t_n=t\atop t_1\geq 0,\dots,t_n\geq 0}}x_1^{t_1}\cdots x_n^{t_n}\partial_1^{[t_1]}\cdots\partial_n^{[t_n]}.\]
In particular $E_1$ is the usual Euler operator $\sum_{i=1}^nx_i\partial_i$.

A graded $\calD$-module $M$ is called {\it Eulerian}, if each homogeneous element $z\in M$ satisfies
\[
E_t\cdot z=\binom{\deg(z)}{t}\cdot z
\]
for every $t\geq 1$.
\end{dfn} 

We collect some basic properties of Eulerian graded $\calD$-modules as follows.
\begin{thm}
\label{Eulerian properties}
Let $M$ be an Eulerian graded $\calD$-module. Then
\begin{enumerate}
\item Graded $\calD$-submodules of $M$ and graded $\calD$-quotients of $M$ are Eulerian.
\item If $S$ is a homogeneous multiplicative system in $R$, then $S^{-1}M$ is Eulerian. In particular, $M_g$ is Eulerian for every homogeneous $g\in R$.
\item The local cohomology modules $H^j_I(M)$ are Eulerian for every homogeneous ideal $I$.
\item The degree-shift $M(\ell)$ is Eulerian if and only if $\ell=0$.
\end{enumerate}
\end{thm}

It follows from Theorem \ref{Eulerian properties} that, if $\Supp_R(H^j_I(R))=\{\frakm\}$ for a homogeneous ideal $I$, then $H^j_I(R)\cong \oplus H^n_{\frakm}(R)^{\mu_j}$ is a degree-preserving isomorphism. Consequently, 
\begin{equation}
\label{vanishing above socle degree}
H^j_I(R)_{\geq -n+1}=0.
\end{equation}
This turns out to be a source of vanishing results for sheaf cohomology. For example, (\ref{vanishing above socle degree}) is one of the ingredients in \cite{BBLSZ2} to prove Theorem \ref{Kodaira for lci thickenings} which is an extension of Kodaira vanishing to a non-reduced setting.

Extensions of Eulerian $\calD$-modules may not be Eulerian as shown in \cite[Remark 3.6]{MaZhangEulerian}. In \cite{PuthenpurakalNagoya2015} the notion of generalized Eulerian $\calD$-module in characteristic 0 was introduced as follows. Fix integers $w_1,\dots,w_n$ and set
\[\deg(x_i)=w_i\ \deg(\partial_i)=-w_i\] 
A graded $\calD$-module $M$ is called {\it generalized Eulerian} if, for every homogeneous element $m\in M$, there is an integer $a$ (which may depend on $m$) such that
\[(E_1-\deg(m))^a\cdot m=0.\]
It was shown that the category of generalized Eulerian $\calD$-modules is closed under extension. This notion of generalized Eulerian $\calD$-modules turns out to be useful in calculating de Rham cohomology of local cohomology modules in characteristic 0 ({\it cf.} \cite{PuthenpurakalNagoya2015}, \cite{PuthenpurakalSingh2019}, \cite{RWZ-lambda}).

In characteristic $p$, the fact that $H^j_I(R)\cong \oplus
H^n_{\frakm}(R)^{\mu_j}$ is a degree-preserving isomorphism when
$\Supp_R(H^j_I(R))=\{\frakm\}$ was also established in
\cite{YiZhangGradedFModules} using $F$-modules, a technique that we
discuss next.

\subsection{$F$-modules}
Let $A$ be a Noetherian commutative ring of characteristic $p$. Then $A$ is equipped with the Frobenius endomorphism
\[F:A\xrightarrow{a\mapsto a^p}A.\] 
The Frobenius endomorphism plays a very important role in the study of rings of characteristic $p$. For instance, in \cite{KunzRegularRing1969}, regularity of $A$ is characterized by the flatness of the Frobenius endomorphism. 

\begin{dfn}[Peskine--Szpiro functor]
\label{dfn: PS functor}
Let $A$ be a Noetherian commutative ring of characteristic $p$. For each $A$-module $M$, denote
by $F_*M$ the $A$-bimodule whose underlying Abelian group is the same as
$M$, whose left $A$-module structure is the usual one: $a\cdot z=az$ for each $z\in F_*M$, and whose right $A$-module structure is given via
the Frobenius $F$:  $z\cdot a:=a^pz$ for each $z\in F_*M$.

The Peskine--Szpiro functor $F_A(-)$ from the category of left $A$-modules to itself is defined via
\[F_A(M):=F_*A\otimes_A M\]
for each $A$-module $M$, where the tensor product uses the right
$A$-structure on $F_*A$.

Geometrically, consider the morphism of spectra induced by the
Frobenius $F\colon A\to A$. Then the right $A$-module structure of
$F_*(M)$ is obtained via restriction of scalars along $F$, and hence
agrees with the pushforward of $M$. On the other hand, $F_A(M)$ is the
pullback of a module under the Frobenius.
\end{dfn}

If $A$ is regular, then it follows from \cite{KunzRegularRing1969} that $F_*A$ is a flat $A$-module and hence $F_A(-)$ is an exact functor.

\begin{rmk}
  \label{remark: Huneke-Sharp approach}
  \begin{asparaenum}
    Let $R$ be a Noetherian regular ring of characteristic $p$ and $I$ be an ideal of $R$.
\item We have 
\begin{align*}F_R(R^m)&\cong R^m,\\
F_R(R/I)&\cong R/I^{[p]}.\end{align*}
Here $I^{[p]}$ is the Frobenius power from Remark \ref{rmk: cofinal sequences}

\item Moreover,
\[F_R(\Ext^j_R(R/I,R))\cong \Ext^j_R(F_R(R/I),F_R(R))\cong \Ext^j_R(R/I^{[p]},R).\]
The natural surjection $R/I^{[p]}\to R/I$ induces
\[
\beta:\Ext^j_R(R/I_,R)\to \Ext^j_R(R/I^{[p]},R)
\]
and by iteration produces a directed system
\[\Ext^j_R(R/I_,R)\xrightarrow{\beta} \Ext^j_R(R/I^{[p]},R)\xrightarrow{F_R(\beta)} \Ext^j_R(R/I^{[p^2]},R)\cdots\]
which agrees with
\[\Ext^j_R(R/I_,R)\xrightarrow{\beta} F_R(\Ext^j_R(R/I,R))\xrightarrow{F_R(\beta)} F^2_R(\Ext^j_R(R/I,R))\cdots\]
Since $\{I^{[p^e]}\}_{e\geq 0}$ and $\{I^t\}_{t\geq 0}$ are cofinal
(that is, the two families of ideals define the same topology on the
ring), the direct limit of this direct system is $H^j_I(R)$.

\item The previous items suggest that $H^j_I(R)$ may be built from the finitely generated $R$-module $\Ext^j_R(R/I_,R)$ using Frobenius, and hence it is natural to expect some properties of $H^j_I(R)$ to be reflected in $\Ext^j_R(R/I_,R)$. Indeed, it was proved in \cite{HunekeSharp} that 
\[\Ass_R(H^j_I(R))\subseteq \Ass_R(\Ext^j_R(R/I_,R)),\ \mu^i_{\frakp}(H^j_I(R))\leq \mu^i_{\frakp}(\Ext^j_R(R/I_,R))\]
for every prime ideal $\frakp$, where $\mu^i_{\frakp}(M)$ denotes the
$i$-th Bass number of an $R$-module $M$ with respect to
$\frakp$. (This was generalized to rings of $F$-finite representation
type in \cite{TakagiTakahashi-MRL08}).
\end{asparaenum}  
Based on the idea of building $H^j_I(R)$ using $\Ext^j_R(R/I_,R)$,
\cite{KatzmanZhangSupportLC} describes a practical algorithm to
calculate the support of $H^j_I(R)$; this algorithm has been
implemented in \emph{Macaulay2} \cite{Macaulay2}.
\end{rmk}

\subsubsection{$F$-modules}
\label{subsubsec-Fmods}
In order to conceptualize the approach in \cite{HunekeSharp},
Lyubeznik introduced the theory of $F$-modules in
\cite{LyubeznikFModules}. Throughout \ref{subsubsec-Fmods}, $R$ is a regular
(not necessarily local) Noetherian ring of characteristic $p>0$, and
$I$ is an ideal of $R$.

\begin{dfn}
 An
\emph{$F$-module} over $R$ (or \emph{$F_R$-module}) is a pair $(M, \theta_M)$
where $M$ is an $R$-module and $\theta_M: M \xrightarrow{\sim} F_R(M)$
is an $R$-module isomorphism, called the \emph{structure
  morphism}. (When the underlying ring is understood, we sometimes refer simply to $M$ as an ``$F$-module''.) The category of $F_R$-modules will be denoted by $\calF_R$ (or $\calF$ when $R$ is clear from the context).

If $R$ is graded, a \emph{graded $F$-module} is an $F$-module $M$ such that $M$ is graded and the structure isomorphism $M\to F_R(M)$ is degree-preserving.
\end{dfn}

\begin{exa}
One can check that $F_*R\otimes_R R\xrightarrow{r'\otimes r\mapsto
  r'r^p}R$ is an $R$-linear isomorphism. Hence $R$ is an $F$-module;
consequently so are all free $R$-modules.

Given any $g\in R$, one can check that $F_*R\otimes_R R_g\xrightarrow{r'\otimes \frac{r}{g^t}\mapsto \frac{r'r^p}{g^{tp}}}R_g$ is an $R$-linear isomorphism. Hence $R_g$ is an $F$-module.

When $R=\kk[x_1,\dots,x_n]$ with standard grading, then for each graded $R$-module $M$ we define a grading on $F_R(M)=F_*R\otimes_RM$ via
\[\deg(r'\otimes m)=\deg(r')+p\deg(m)\]
for all homogeneous $r'\in R$ and $m\in M$.

In this setting, $F_*R\otimes_R R\xrightarrow{r'\otimes r\mapsto
  r'r^p}R$ is a degree-preserving $R$-linear isomorphism and so $R$ is a graded $F$-module. Likewise, if $g\in R$ is homogeneous, then $F_*R\otimes_R R_g\xrightarrow{r'\otimes \frac{r}{g^t}\mapsto \frac{r'r^p}{g^{tp}}}R_g$ is a degree-preserving $R$-linear isomorphism and hence $R_g$ is a graded $F$-module.
\end{exa}

\begin{dfn}
Let $(M, \theta_M)$ be an $F$-module. We say that $M$ is \emph{$F$-finite} if there exists a finitely generated $R$-module $M'$ and an $R$-linear map $\beta: M' \rightarrow F_R(M')$ such that 
\begin{gather}\label{eqn-F-Finite}
\varinjlim(M' \xrightarrow{\beta} F_R(M') \xrightarrow{F^*\beta} F^2_R (M') \to \cdots) \cong M,
\end{gather}
and the structure morphism $\theta_M$ is induced by taking the direct limit over $\ell$ of $F^\ell_R(\beta): F^\ell_R (M') \rightarrow F^{\ell+1}_R (M')$. In this case we call $M'$ a \emph{generator} of $M$ and $\beta$ a \emph{generating morphism}. A generator $M'$ of an $F$-finite $F$-module $M$ is called a {\it root} if the generating morphism $\beta: M' \rightarrow F_R(M')$ is injective.

A graded $F$-finite $F$-module is defined to be an $F$-finite
$F$-module for which the modules and morphisms in \eqref{eqn-F-Finite}
can be chosen to be homogeneous.
\end{dfn}

\begin{exa}
From Remark \ref{remark: Huneke-Sharp approach}, one can see that
every local cohomology module $H^j_I(R)$ is an $F$-finite $F$-module since it is the direct limit of
\[\Ext^j_R(R/I_,R)\xrightarrow{\beta} F_R(\Ext^j_R(R/I,R))\xrightarrow{F_R(\beta)} F^2_R(\Ext^j_R(R/I,R))\cdots\]
and $\Ext^j_R(R/I_,R)$ is finitely generated.

When $R=\kk[x_1,\dots,x_n]$ and $I$ is a homogeneous ideal of $R$, the local cohomology modules $H^j_I(R)$ are graded $F$-finite $F$-modules.
\end{exa}

There is a fruitful analogy between ($F$-finite) $F$-modules and (holonomic) $\calD$-modules. We collect some basic properties of $F$-modules, which are parallel to those of $\calD$-modules, as follows.
\begin{thm}
\label{thm: properties F-modules}
Let $R$ be a Noetherian regular ring of characteristic $p>0$.
\begin{enumerate}
\item If $M$ is an $F$-module, then $\injdim_R(M)\leq
  \dim(\Supp_R(M))$, \cite[1.4]{LyubeznikFModules}.
\item $F$-finite $F$-modules form a full Abelian subcategory of the
  category of $R$-modules that is
  closed under the formation of submodules, quotient modules, and
  extensions, \cite[2.8]{LyubeznikFModules}.
\item If $M$ is an $F$-finite $F$-module, then so is the localization
  $M_g$ for each $g\in R$, \cite[2.9]{LyubeznikFModules}.
\item A simple $F$-module has a unique associated prime, \cite[2.12]{LyubeznikFModules}.
\item $F$-finite $F$-modules have finite length in the category of $F$-modules, \cite[3.2]{LyubeznikFModules}.
\end{enumerate}
\end{thm}

\begin{rmk}
  The theory of $F$-modules plays a crucial role in the extension of the
  Riemann--Hilbert correspondence to characteristic $p$ by 
  Emerton and Kisin \cite{Emerton-Kisin}, which is beyond the scope of this survey.
\end{rmk}

\subsubsection{$A\{f\}$-modules: action of Frobenius}
\label{subsec-f-action}

Let $A$ be a Noetherian commutative ring of characteristic $p$. We will use $A\{f\}$ to denote the associative $A$-algebra with one generator $f$ and relations $fa=a^pf$ for all $a\in A$. 

\begin{rmk}\label{rmk-F-action}
\label{characterization of Frobenius modules}
Let $M$ be an $A$-module $M$. The following are equivalent. 
\begin{enumerate}
\item $M$ is an $A\{f\}$-module.
\item $M$ admits an additive map $f:M\to M$ such that $f(am)=a^pf(m)$ for every $a\in A$ and $m\in M$; this $f$ is called a \emph{Frobenius action} on $M$.
\item $M$ admits an $A$-linear map $M\to F_*M$
   where $F\colon A\to A$ is the Frobenius endomorphism on $A$.
\item $M$ admits an $A$-linear map $F_*A\otimes_AM\to M$ where $F:A\to A$ is the Frobenius endomorphism on $A$.
\end{enumerate}
In (2), we still use $f$ to denote the Frobenius action since multiplication on the left by $f$ on $M$ is indeed a Frobenius action for each $A\{f\}$-module $M$. 

Of course, the standard example of a Frobenius action is $A$ with the
$p$-th power map. Note that the image $f(M)$ is in general just a
group, but acquires the structure of a $\kk$-space when $\kk$ is perfect.

The Frobenius on $A$ induces a natural Frobenius action on each $H^i_{\fraka}(A)$ for every ideal $\fraka$; hence $H^i_{\fraka}(A)$ is an $A\{f\}$-module. In this paper, we always consider $H^i_{\fraka}(A)$ as an $A\{f\}$-module with the Frobenius action $f$ induced by the Frobenius endomorphism on $A$. For this reason, some authors denote by $F$ (instead of $f$) the Frobenius action on $H^i_{\fraka}(A)$ induced by the Frobenius endomorphism on $A$.
\end{rmk}

\begin{dfn}
Given an $A\{f\}$-module $M$ with Frobenius action $f:M\to
M$, the intersection
\[
M_\st:=\bigcap_{t\geq 1}f^t(M)
\]
is called the \emph{$f$-stable part} of $M$.

An element $z\in M$ is called \emph{$f$-nilpotent} if $f^t(z)=0$ for some integer $t$. 

An $A\{f\}$-module $M$ is called \emph{$f$-torsion} if every element
in $M$ is in the kernel of some iterate of $f$, and it is called \emph{$f$-nilpotent} if there is an integer $t$ such that $f^t(M)=0$.
\end{dfn}

\begin{rmk}
When $M=H^i_{\fraka}(A)$ is a local cohomology module of $A$, the notions of $f$-torsion and $f$-nilpotent are also denoted by $F$-torsion and $F$-nilpotent, respectively, since the Frobenius action $f$ is induced by the Frobenius endomorphism on $A$.

Assume $(A,\frakm,\kk)$ is a local ring and $x_1\dots,x_d$ is a full system of parameters. Then the Frobenius action $f$ on $H^d_{\frakm}(A)$ can be described as follows. Let $\eta=[\frac{a}{x^{n_1}_1\cdots x^{n_d}_d}]$ be an element in $H^d_{\frakm}(A)$, then 
\[f(\eta)=[\frac{a^p}{x^{n_1p}_1\cdots x^{n_dp}_d}].\]
\end{rmk}

An $A\{f\}$-module that is also an Artinian $A$-module is called a {\it cofinite} $A\{f\}$-module. Cofinite $A\{f\}$-modules enjoy an amazing property. 

\begin{thm}
\label{nilpotency on cofinite modules}
Let $A$ be a local ring of characteristic $p>0$. Assume that $M$ is an $f$-torsion cofinite $A\{f\}$-module. Then $M$ must be $f$-nilpotent. 
\end{thm}

Theorem \ref{nilpotency on cofinite modules} was first proved by
Hartshorne and Speiser in \cite{HartshorneSpeiser-Annals77}.
There, Hartshorne and Speiser created a version of some of Ogus'
results from \cite{Ogus-LCDAV} in characteristic $p>0$. Their
motivating question was to determine when the cohomology of every
coherent sheaf on the complement of a projective variety be a finite
dimensional vector space.  Hartshorne and Speiser use the Frobenius
endomorphism on $\calO_{\hat X}$ to supply the information given by
the connection used by Ogus in characteristic zero, and
$\ZZ/p$-\'etale cohomology turns up in place of de Rham cohomology.
Theorem \ref{nilpotency on cofinite modules} was
later
generalized by Lyubeznik in \cite{LyubeznikFModules}  (using the
$\calH_{R,A}$-functor discussed in the sequel). It has found applications in \cite{KatzmanLyubeznikZhang-JA2009, BlickleSchwedeTakagiZhang-MathAnn2010, BlickleBoeckle-Crelle11} in the study of singularities and invariants defined by Frobenius.

Theorem 4.6 in \cite{Lyubeznik-Compositio06} 
reads as follows: if $\kk$ is an algebraically closed field of
positive characteristic, and if $(A,\frakm,\kk)$ is a complete local
ring with connected punctured spectrum and $\kk\subseteq A$, then
$H^1_\frakm(A)$ is $f$-torsion.  Lyubeznik derives this via a
comparison with local cohomology in a complete regular local ring that
surjects onto $A$.  In \cite{SinghWalther-TAMS08}, this result is
sharpened to a numerical statement over an algebraically closed
coefficient field: the number of connected components
of the punctured spectrum of $A$ is one more than the dimension of
the $f$-stable part of $H^1_\frakm(A)$.

A general study of Frobenius operators started with
\cite{LyubeznikSmith-TAMS01} and later was carried out by various
authors: aside from Sharp's article \cite{Sharp-TAMS07} we should point at
\cite{Sharp-PAMS07} by the same author, \cite{BlickleBoeckle-Crelle11}
who develop the notion of \emph{Cartier modules} (which are
approximately modules with a Frobenius action), and \cite{Gabber-04}. The
article \cite{KatzmanSchwedeSinghZhang} contains positive results on
finiteness dual to \cite{Schwede-TAMS11} as well as examples of
failure.

\begin{dfn} Let $(A,\frakm,\kk)$ be a local ring of characteristic
  $p>0$. 
Given a cofinite $A\{f\}$-module $W$, a prime ideal $\frakp$ is called a {\it special prime} of $W$ if it is the annihilator of an $A\{f\}$-submodule of $W$.
\end{dfn} 

It is proved in \cite[Corollary 3.7]{Sharp-TAMS07} and \cite[Theorem
  3.6]{EnescuHochster-ANT08} that if the Frobenius action $f\colon M\to M$
on the $A\{f\}$-module $M$ is injective then $M$ admits only finitely many special primes. This will be useful when we discuss the $F$-module length of local cohomology modules in the sequel.

\begin{dfn}[\cite{EnescuHochster-ANT08}]
Let $(A,\frakm)$ be a Noetherian local ring of characteristic $p$. Let $f:H^j_{\frakm}(A)\to H^j_{\frakm}(A)$ denote the Frobenius action induced by the Frobenius on $A$. 

A submodule $N$ of $H^j_{\frakm}(A)$ is called {\it $F$-stable} if $f(N)\subseteq N$. 

The ring $A$ is called {\it $\FH$-finite} if $H^j_{\frakm}(A)$ admits
only finitely many $F$-stable submodules for each $0\leq j\leq \dim(A)$.

Also, $A$ is called \emph{$F$-injective} if the natural Frobenius action $f:H^j_{\frakm}(A)\to H^j_{\frakm}(A)$ is injective for each integer $j\leq\dim(A)$.
\end{dfn}

The Frobenius action on local cohomology modules connects with a very
important type of singularities, that of $F$-rationality, which we  recall next.

\begin{dfn}
\label{dfn-tightcl}
  Let $A$ be a Noetherian
  ring of characteristic $p$, let $A^{\circ}$ denote the complement of the union of minimal primes in $A$ and let $\fraka$ be an ideal of $A$. An element $a\in A$ is in the {\it tight closure} of $\fraka$ if there is a $c\in A^\circ$ such that $ca^{p^e}\in \fraka^{[p^e]}$ for all $e\gg 0$. Let $\fraka^*$ denote the set of elements $a\in A$ that are in the tight closure of $\fraka$; it is an ideal of $A$. An ideal $\fraka$ is called {\it tightly closed} if $\fraka=\fraka^*$.

A local ring $A$ is called {\it $F$-rational} if $\fraka=\fraka^*$ for every parameter ideal $\fraka$.
\end{dfn}

In her work to relate $F$-rationality (an algebraic notion) to
rational singularity (a geometric notion), Smith \cite{Smith-AJM97}
proves the following characterization of $F$-rationality using a Frobenius action.

\begin{thm}
Let $(A,\frakm)$ be a $d$-dimensional excellent local domain of characteristic $p$. Then $A$ is $F$-rational if and only if $A$ is normal, Cohen--Macaulay, and $H^d_{\frakm}(A)$ contains no non-trivial $F$-stable submodules. 
\end{thm}

In independent work of
Smith, Mehta--Srinivas, and Hara, $F$-rationality was shown to be the
algebraic counterpart to the notion of rational singularities
\cite{MehtaSrinivas-Asian97,Smith-AJM97,Hara-AJM98}. The purpose of
these studies was to establish a parallelism between the concept of a
rational singularity in characteristic zero, and invariants based on
the Frobenius for its models in finite (large) characteristic. The development
of such connections has a fascinating and distinguished history, and
we recommend the recent and excellent survey article
\cite{TakagiWatanabe-Sugaku18} by two experts in the field.

A related construction goes back to \cite{Enescu-PAMS03}. For an element $x\in A$ and a parameter ideal $I$ of $A$ let $I(x)$ be
the ideal of elements $c\in A$ that multiply $x^{p^e}$ into
$I^{[p^e]}$ for all large $e$ ({\it cf.} Definition \ref{dfn-tightcl}). Enescu shows in 
\cite{Enescu-PAMS03} that if $A$ is $F$-injective and
Cohen--Macaulay, then the set of maximal elements in $\{I(x):x\not\in
I\}$ does not depend on $I$, is finite and consists only of prime
ideals. These are called \emph{$F$-stable primes}, and the collection of them is denoted by
$\FS(R)$. Enescu shows further that for an $F$-injective
Cohen--Macaulay complete local ring $A$, the $F$-stable primes can be
expressed in terms of $F$-unstability, introduced by Fedder and
Watanabe. Enescu and Sharp continued the  study of
properties of $F$-stable primes in
\cite{Sharp-TAMS07,Enescu-JA09}.

Along with $\FH$-finiteness goes another property of rings that will
come back to us later:
\begin{dfn}\label{dfn-F-pure}
  $A$ is called {\it $F$-pure} if $(A\xrightarrow{a\mapsto
  a^p}A)\otimes_A M$ is injective for all $A$-modules $M$.
\end{dfn}

\begin{rmk}
  For background to this remark we
  refer to the excellent article \cite{TakagiWatanabe-Sugaku18}.

  A standard question on ``deformation'' in commutative algebra is to
  ask ``If a quotient $A/\ideal{x}$ of $A$ by a regular element has a nice property, is $A$
  forced to share it?''.
  
  It turns out that $F$-purity does not deform in this sense,
  \cite{Fedder-TAMS83,Singh-JPAA99}. 
  The
  reader familiar with the concepts of $F$-regularity and
  $F$-rationality may know that $F$-rationality deforms
  \cite{HH-TAMS94} while
  $F$-regularity does not
  \cite{Anurag-AJM99} although it does so for $\QQ$-Gorenstein rings \cite{HH-TAMS94,AberbachKatzmanMcCrimmon-JA98}. Very recently, Polstra and Simpson proved in \cite{PolstraSimpsonFPurityDeformsQGorenstein} that $F$-purity deforms in $\QQ$-Gorenstein rings.
  
It is still an open question whether $F$-injectivity deforms, but some
progress has been made. Fedder showed in \cite{Fedder-TAMS83} that
$F$-injectivity deforms when the ring is Cohen--Macaulay. In
\cite{HoriuchiMillerShimomoto-IU14}, it was proved that if $R/xR$ is
$F$-injective and $H^j_{\frakm}(A/(x^\ell) )\to H^j_{\frakm}(A/(x))$
is surjective for all $\ell>1$ and $j$ then $A$ is $F$-injective. Ma
and de Stefani established deformation when the local cohomology
modules $H^\bullet_\frakm(A)$ have secondary decompositions that are
preserved by the Frobenius \cite{MaStefani-2009.09038}.
\end{rmk}

In \cite{EnescuHochster-ANT08} it is proved that  face rings of finite
simplicial complexes are \FH-finite. They showed further that an $F$-pure and
quasi-Gorenstein local ring is $\FH$-finite, and raised the question 
whether all $F$-pure and Cohen--Macaulay local rings are $\FH$-finite. Ma answered this question in the affirmative by proving the following result in \cite{Ma-IMRN14}.

\begin{thm}
Let $(A,\frakm)$ be a Noetherian local ring of characteristic $p$. If $A$ is $F$-pure, then $A$ and all power series rings over $A$ are $FH$-finite.
 \end{thm}
In the paper he also proved that if $A$ is $F$-pure (even just on the
punctured spectrum) then $H^\bullet_\frakm(A)$ is a finite length
$A\{f\}$-module, and he also established that the finite length
property is stable under localization. With Quy, he introduced more recently in
\cite{MaQuy-Nagoya18} the notions \emph{$F$-full} (when the Frobenius
action is surjective) and \emph{$F$-anti-nilpotent} (when the action
is injective on every $A\{f\}$-subquotient of local cohomology). They
established that $F$-anti-nilpotence implies $F$-fullness and equals
\FH-finiteness of \cite{EnescuHochster-ANT08}.
Inspired by ideas
from \cite{HoriuchiMillerShimomoto-IU14},
they prove the interesting fact
that both $F$-anti-nilpotence and $F$-fullness do deform.

The action of the Frobenius also ties in naturally with the action of
the Frobenius on the cohomology of projective varieties via the
identification \eqref{eqn-cd-Serre}. For example, the Segre product of
a smooth elliptic curve $E$ with $\PP^1_\KK$ has $F$-injective
coordinate ring (recall Definition \ref{dfn-F-pure}) if and only if
the curve is ordinary (the group $H^1(E;\calO_E)$ is the degree zero
part of $H^2_\frakm(A)$ and the Frobenius action is the induced one;
here $A$ is the coordinate ring of $E$). Compare Example
\ref{exa-ara-not-lcd}.

Hartshorne and Speiser in \cite{HartshorneSpeiser-Annals77}, and
Fedder and Watanabe in \cite{FedderWatanabe-CA89} studied 
$F$-actions on local cohomology with regards to vanishing of
cohomology on projective varieties, and with regards to singularity
types of local rings respectively.

According to \cite{SrinivasTakagi-Adv17}, a local ring $(A,\frakm)$ is
\emph{$F$-nilpotent} if the Frobenius action is nilpotent on
$H^{<\dim(A)}_\frakm(A)$ and $0^*_{H^{\dim(A)}_\frakm(A)}$ (the tight closure of the zero submodule of $H^{\dim(A)}_\frakm(A)$), and Srinivas and Takagi show that $A$ is
$F$-injective and $F$-nilpotent if and only it is $F$-rational. In
\cite{PolstraQuy-JA19}, Polstra and Quy characterize 
$F$-nilpotence as (under mild hypotheses) being equivalent
to the equality of tight and Frobenius closure for all parameter
ideals. This work extends the finite length case discussed in
\cite{Ma-MathAnn15} and is somewhat surprising since the complementary
notion of $F$-injectivity is not equivalent to the
Frobenius-closedness of all parameter ideals,
\cite{QuyShimomoto-Adv17}, but only implied by it.

Ma also shows in \cite{Ma-MathAnn15}, in his setting of finite length
lower local cohomology, that $F$-injectivity implies the ring being
Buchsbaum (a generalization of Cohen--Macaulay,
\cite{StueckradVogel}), and that the analogous statement in
characteristic zero is true in the sense that, if $A$ is a normal
standard graded $\KK$-algebra with $\KK\supseteq \QQ$ that is Du Bois
and has finite length lower local cohomology, then $A$ is Buchsbaum. (A
singularity $X$ embedded inside a smooth scheme over the complex
numbers is du Bois, following Schwede's paper
\cite{Schwede-Compositio07}, if and only if an embedded resolution
$\pi\colon Y\to X$ of $X=\Spec(A)$ with reduced total transform $E$
leads to an isomorphism $\calO_X=R\pi_*(\calO_E)$. Initially, Du Bois
singularities arose from Hodge-theoretic filtrations of the de Rham
complex in \cite{DuBois-BSMF81}; they include normal crossings and
quotient singularities). Du Bois singularities are closely related to (and conjecturally equivalent to) singularities of dense $F$-injective type. Recall that, a finite $\ZZ$-algebra $A_\ZZ$ is of \emph{dense $F$-injective type} if its
reductions $A_p$ modulo $p$ are $F$-injective for infinitely many primes $p\in\ZZ$. Schwede proved in \cite{Schwede-AJM09} that if a finite $\ZZ$-algebra $A_\ZZ$ is of dense $F$-injective type then the complex model $A_\CC=A_\ZZ\otimes_\ZZ \CC$ is Du Bois. The other implication remains an open problem and was proved to be equivalent to the Weak Ordinarity Conjecture (see \cite{BhattSchwedeTakagi} for details). 

We close this section with a brief discussion on the very interesting
topic of the interaction of the Frobenius with Hodge theory, crossing
characteristics.  Suppose $A$ is a finitely generated graded
$\CC$-algebra, and set $X:=\Proj(A)$.  It is known that certain
aspects of the Hodge theory of $X$ are encoded in the combinatorics of
the resolution of singularities of $X$,
\cite{Deligne-HodgeII,Deligne-HodgeIII,DonuParsaJarek-MathAnn13}. In
this context, Srinivas and Takagi proposed and studied in
\cite{SrinivasTakagi-Adv17} the following local conjecture.
\begin{conj}\label{conj-Fnil}
  If $\pointx$ is a normal isolated singularity on the $n$-dimensional
  $\CC$-scheme $X$ then the local ring at $\pointx$ is of
  $F$-nilpotent type if and only if for all $i<\dim(X)$, the zeroth
  graded piece $\Gr^0_F(H^i_{\{\pointx\}}(X^\an,\CC))$ of the Hodge
  filtration is zero.
\end{conj}
How much is still unknown in this fascinating area between
characteristics can be seen from the fact that the following
conjectural statement
is still open: let $V$ be an $(n-1)$-dimensional projective simple
normal crossings variety in characteristic zero; then the Frobenius
action on $H^i(V_p,\calO_{V_p})$ is not nilpotent for an infinite set
of reductions $V_p$ modulo $p$ of $V$. Srinivas and Takagi
\cite{SrinivasTakagi-Adv17} prove the
case $n-1=2$ of this and derive from it the case
$n=3$ for the conjecture above.

\subsubsection{The Lyubeznik functor $\calH_{R,A}$}

Assume that $A$ is a homomorphic image of a Noetherian regular ring $R$. The approach of building $H^j_I(R)$ using a finitely generated $R$-module results in a very useful functor $\calH_{R,A}$ from the category of cofinite $A\{f\}$-modules to the category of $F_R$-finite $F_R$-modules.

\begin{rmk}
Let $R=\kk[[x_1,\dots,x_n]]$ and $E=H^n_{(x_1,\dots,x_n)}(R)$. Denote
as before the Matlis dual functor $\Hom_R(-,E)$ by $D(-)$. Then there is a functorial $R$-module isomorphism
\[\tau:D(F_R(M))\cong F_R(D(M))\]
for all Artinian $R$-modules $M$. 

Let $A$ be a homomorphic image of $R$. Let $M$ be an $A\{f\}$-module. One can check that
\begin{gather}\label{eqn-alpha}
  \alpha: F_R(M)\xrightarrow{r\otimes m\mapsto rf(m)}M
\end{gather}
is an $R$-module homomorphism. Now, assume that $M$ is a cofinite $A\{f\}$-module. Taking the Matlis dual of $\alpha$, we have an $R$-module homomorphism
\[\beta=\tau\circ D(\alpha): D(M)\to F_R(D(M)),\]
and hence we have a direct system of Noetherian $R$-modules:
\[D(M)\xrightarrow{\beta} F_R(D(M))\xrightarrow{F_R(\beta)}F^2_R(D(M))\to \cdots\]

Analogously, let $R=\kk[x_1,\dots,x_n]$ and denote the graded Matlis dual functor $\sideset{^*}{_R}\Hom(-,E)$ by $D^*(-)$. There is a functorial graded $R$-module isomorphism
\[\tau:D^*(F_R(M))\cong F_R(D^*(M))\]
for all Artinian graded $R$-modules $M$.

Assume $R=\kk[x_1,\dots,x_n]$ and $A$ is a graded homomorphic image of $R$. By a graded $A\{f\}$-module, we mean a graded $A$-module $M$ such that $f:M_\ell\to M_{p\ell}$ for all integers $\ell$. One can check that then \eqref{eqn-alpha}
is a graded $R$-module homomorphism. Now, assume that $M$ is a cofinite graded $A\{f\}$-module. Taking the graded Matlis dual of $\alpha$, we have a graded $R$-module homomorphism
\[\beta=\tau\circ D^*(\alpha): D^*(M)\to F_R(D^*(M)),\]
and hence we have a direct system of graded Noetherian $R$-modules:
\[D^*(M)\xrightarrow{\beta} F_R(D^*(M))\xrightarrow{F_R(\beta)}F^2_R(D^*(M))\to \cdots\]
\end{rmk}

\begin{dfn}\label{dfn-HRA}
Let $R$ be a complete regular local ring $R$ of characteristic $p$ and let $A$ be a homomorphic image of $R$. For each cofinite $A\{f\}$-module $M$, we define
\[\calH_{R,A}(M):=\varinjlim(D(M)\xrightarrow{\beta} F_R(D(M))\xrightarrow{F_R(\beta)}F^2_R(D(M))\to \cdots)\]

The graded version $\calH^*_{R,A}$ is defined analogously on
homogeneous input.
\end{dfn}

\begin{exa}
\label{H functor on lc}
Let $R=\kk[[x_1,\dots,x_n]]$ (or $R=\kk[x_1,\dots,x_n]$ respectively) and let
$I$ be an ideal of $R$ (homogeneous, if  $R=\kk[x_1,\dots,x_n]$). Set $A=R/I$. Hence $H^{j}_{\frakm}(A)$ is an $A\{f\}$-module according to Remark \ref{characterization of Frobenius modules}. Since $H^{j}_{\frakm}(A)$ is Artinian, it is a cofinite $A\{f\}$-module. Then one can check that
\[\calH_{R,A}(H^j_{\frakm}(A))\cong H^{n-j}_I(R)\]
(which reads \[\calH^*_{R,A}(H^j_{\frakm}(A))\cong H^{n-j}_I(R)\]
  when $R=\kk[x_1,\dots,x_n]$).
\end{exa}

\begin{rmk}
The functor $\calH_{R,A}$ (resp.\ $\calH^*_{R,A}$) from the category of cofinite (graded) $A\{f\}$-modules to the category of (graded) $F$-finite $F$-modules is contravariant, additive, and exact.

Given a cofinite (graded) $A\{f\}$-module $M$, $\calH_{R,A}(M)=0$ (or $\calH^*_{R,A}(M)=0$ respectively) if and only if the additive map $\varphi:M\to M$ in Remark \ref{characterization of Frobenius modules} is nilpotent.

Now Lyubeznik's vanishing theorem in characteristic $p$ follows from Example \ref{H functor on lc}: $H^{n-j}_I(R)=0$ if and only if the natural Frobenius (induced by the Frobenius on $R$) on $H^{j}_{\frakm}(A)$ is nilpotent.
\end{rmk}

The nilpotence of the action of Frobenius on $H^{j}_{\frakm}(A)$ prompts the following definition ({\it cf.} \cite[Definition 4.1]{Lyubeznik-Compositio06}).

\begin{dfn}
Let $(A,\frakm)$ be a local ring of characteristic $p$. The $F$-depth of $A$ is the smallest $i$ such that $H^i_{\frakm}(A)$ is not $f$-nilpotent, where $f$ is the natural action of Frobenius on $H^i_{\frakm}(A)$ induced by the Frobenius endomorphism on $A$.
\end{dfn}

\begin{rmk}\label{rmk-Fdepth}
One can show that ({\it cf.} \cite[\S4]{Lyubeznik-Compositio06})
\begin{enumerate}
\item $\depth(A)\leq F\text{-depth}(A)\leq \dim(A)$,
\item $F\text{-depth}(A)=F\text{-depth}(\hat{A})$,
\item $F\text{-depth}(A)=F\text{-depth}(A_\red)$ where $A_\red=A/\sqrt{(0)}$.
\end{enumerate}
In terms of $F$-depth, the vanishing theorem via Frobenius in characteristic $p$ can be restated as follows: {\it let $(R,\frakm)$ be a regular local ring of characteristic $p$ and $I$ be an ideal. Then}
\[\lcd_R(I)=\dim(R)-F\text{-depth}(R/I).\]
(Compare also the corresponding statement in characteristic zero,
Theorem \ref{thm-Ogus-dRdepth}). 
\end{rmk}

In general, $F\text{-depth}(A)$ can be different from $\depth(A)$ as shown in the following example ({\it cf.} \cite[\S5]{Lyubeznik-Compositio06}).

\begin{exa}
Let $\kk$ be a perfect field of characteristic $p$ and let
$\calC\subseteq \PP^2_{\kk}$ denote the Fermat curve defined by
$x^3+y^3+z^3$. Let $R=\kk[x_0,\dots,x_5]$ and $I\subseteq R$ be the
defining ideal of $\calC\times \PP^1_{\kk}\subseteq \PP^5_{\kk}$. Set
$A=(R/I)_{\frakm}$ where $\frakm=(x_0,\dots,x_5)$.

If $3\mid (p-2)$, then
\[F\text{-depth}(A)=3>2=\depth(A).\]
See also Example \ref{exa-ara-not-lcd}.
\end{exa} 


Since $F$-finite $F$-modules have finite length in the category of
$F$-modules, it is natural to ask whether one can compute the length,
especially for local cohomology modules. It turns out that $F$-module
length of local cohomology modules is closely related to singularities
defined by the Frobenius, and
Lyubeznik's functor $\calH_{R,A}$ is a useful tool for studying this length.
To illustrate this, let $R$ be a regular local ring of characteristic
$p$. That $\calH_{R,A}$ sets up  a link between the length of
$H^{\hight(I)}_I(R)$ and the singularities of $A=R/I$ was first
discovered in \cite{BlickleIntersectionHomology2004}; this was later
extended and strengthened in \cite{KatzmanMaSmirnovZhang2018} as
follows, see also \cite{Bitoun-IMRN20}

\begin{thm}
\label{F-length KMSZ}
Let $R=\kk[[x_1,\dots,x_n]]$ (or $R=\kk[x_1,\dots,x_n]$), and set
$\frakm=\ideal{x_1,\dots,x_n}$. Let $A=R/I$ be a reduced and
equidimensional (and graded, if $R=\kk[x_1,\dots,x_n]$) ring of dimension $d\geq 1$. Let $c$ denote the number of minimal primes of $A$.

\begin{enumerate}
\item If $A$ has an isolated non-$F$-rational point at $\frakm$ and $\kk$ is separably closed, then 
\[\ell_{\calF_R}(H^{n-d}_I(R))=\dim_{\kk}(H^d_{\frakm}(A)_\st)+c.\]

\item If the non-$F$-rational locus of $A$ has dimension $\leq 1$ and $\kk$ is separably closed, then 
\[\ell_{\calF_R}(H^{n-d}_I(R))\leq \sum_{\dim(A/\frakp)=1}\dim_{\kappa(\frakp)}(H^{d-1}_{\frakp\widehat{A_{\frakp}}}(A_{\frakp})_\st)+\dim_{\kk}(H^d_{\frakm}(A)_\st)+c,\]

\item If $A$ is $F$-pure, then $\ell_{\calF_R}(H^{n-d}_I(R))$ is at least the number of special primes of $H^d_{\frakm}(A)$. Moreover, if $A$ is $F$-pure and quasi-Gorenstein, then $\ell_{\calF_R}(H^{n-d}_I(R))$ is precisely the number of special primes of $H^d_{\frakm}(A)$.
\end{enumerate}
\end{thm}

It remains an open problem whether one can extend Theorem \ref{F-length KMSZ} to the case of a higher dimensional non-$F$-rational locus.

Recently, in \cite{MontanerBoixZarzuelaIMRN2020}, \`Alvarez Montaner, Boix and Zarzuela computed $\ell_{\calF}(H^j_I(R))$ and $\ell_{\calD}(H^j_I(R))$ when $R$ is a polynomial ring over a field and $I$ is generated by square-free monomials and pure binomials ({\it i.e.} $I$ is a toric face ideal).

\subsection{Interaction between $D$-modules and $F$-modules}
In characteristic $p$, the theories of $D$-modules and $F$-modules are entwined; it has been fruitful to consider local cohomology modules from both perspectives.

\begin{rmk}
\label{F-modules are D-modules}
Let $\kk$ be a field of characteristic $p$ and let
$R=\kk[x_1,\dots,x_n]$ or $ R=\kk[[x_1,\dots,x_n]]$. It is clear from
the definition that, if $(M,\theta)$ is an $F$-module, the map
\[
\alpha_e:M\xrightarrow{\theta}F_R(M)\xrightarrow{F_R(\theta)}F^2_R(M)\xrightarrow{F^2_R(\theta)}
\cdots \to F^e_R(M)
\]
induced by $\theta$ is also an isomorphism.

This induces a $\calD=\calD(R,\kk)$-module structure on $M$. To
specify the induced $\calD$-module structure, it suffices to specify
how $\partial^{[i_1]}_1\cdots \partial^{[i_n]}_n$ acts on $M$. Choose
$e$ such that $p^e\geq (i_1+\cdots +i_n)+1$. Given $z\in M$, we
consider $\alpha_e(z)$ and we will write it as $\sum r_j\otimes z_j$
with $r_j\in {F^e_*R}$ and $z_j\in M$. Then define
\[
\partial^{[i_1]}_1\cdots \partial^{[i_n]}_n \cdot
z:=\alpha_e^{-1}(\sum \partial^{[i_1]}_1\cdots \partial^{[i_n]}_n
r_j\otimes z_j);
\]
that this is legal is due to a simplification in the product rule in
characteristic $p$: $(x^pg)'=x^p(g)'$.

Therefore, every $F$-module is also a $\calD$-module.
\end{rmk}

When $R=\kk[x_1,\dots,x_n]$ with its standard grading, the
$\calD$-module structure on each graded $F$-module as in Remark
\ref{F-modules are D-modules} is also graded. Moreover,
\cite{MaZhangEulerian} proves the following:
\begin{thm}
Let $R=\kk[x_1,\dots,x_n]$. Every graded $F$-module is an Eulerian graded $\calD$-module.
\end{thm}

Since every $F$-module is a $\calD$-module, given an $F$-finite $F$-module $M$, one may compare $\ell_{\calF}(M)$ and $\ell_{\calD}(M)$. A quick observation is that, since each filtration of $M$ in $\calF$ is also a filtration in $\calD$, one always has
\[\ell_{\calF}(M)\leq \ell_{\calD}(M).\]

It turns out that this inequality can be strict.

\begin{thm}[Proposition 7.5 in \cite{KatzmanMaSmirnovZhang2018}]
Let $p$ be a prime number such that $7\ |\ (p-4)$. Let $R=\overline{\FF}_p[x,y,z,t]$ and $f=tx^7+ty^7+z^7$. Then
\[\ell_{\calF}(H^1_{(f)}(R))=3<7= \ell_{\calD}(H^1_{(f)}(R)).\]
\end{thm}

On the other hand, the equality holds when hypotheses are added:
\begin{thm}
Let $R,I,A$ be as in Theorem \ref{F-length KMSZ}. If $A$ has an isolated non-$F$-rational point at $\frakm$ and $\kk$ is separably closed, then
\[\ell_{\calF}(H^{n-d}_I(R))=\ell_{\calD}(H^{n-d}_I(R)).\]
\end{thm}

$F$-modules and $\calD$-modules are deeply connected via a generating property. The following is a special case of \cite[Corollary 4.4]{AlvarezBlickleLyubeznik}.

\begin{thm}
\label{generators of root generate D-module}
Let $\kk$ be a field of characteristic $p$ such that $[\kk:\kk^p]<\infty$ and let $R=\kk[x_1,\dots,x_n]$ or $R=\kk[[x_1,\dots,x_n]]$. Let $M$ be an $F$-finite $F$-module. If $z_1,\dots,z_t\in M$ generate a root of $M$, then $z_1,\dots,z_t\in M$ generate $M$ as a $\calD$-module.
\end{thm}

Theorem \ref{generators of root generate D-module} plays a crucial
role in proving that $1/g$ generates $R_g$ as a $\calD$-module in
\cite{AlvarezBlickleLyubeznik}, and also in proving the finiteness of associated primes of local cohomology of smooth $\ZZ$-algebras in \cite{BBLSZ-2014}.

\section{Local cohomology and topology}
\label{sec-top}

In this section we discuss the interaction of local cohomology with
various themes of topological flavor. The interactions can typically be seen as a failure of flatness in
some family witnessed by specific elements of certain local cohomology.

We start with a classical discussion of the number of defining
equations for a variety, then elaborate on the more recent
developments that originate from this basic question. We survey
interactions with topology in characteristic zero, and with the
Frobenius map in positive characteristic. We discuss a collection of
applications of local cohomology to various areas: hypergeometric
functions, the theory of Milnor fibers, the Bockstein morphism from
topology. We close with a discussion on a set of numerical invariants
based on local cohomology modules introduced by Lyubeznik.

\subsection{Arithmetic rank}\label{subsec-ara}

The main object of interest here is described in our first definition. 
\begin{dfn}\label{dfn-ara}
  The \emph{arithmetic rank} $\ara_A(I)$ of the $A$-ideal $I$ is the
  minimum number of generators for an ideal with the same radical as $I$:
  \[
  \ara_A(I)=\min\{\ell\in\NN \mid \exists x_1,\ldots,x_\ell\in A,\
  \sqrt{I}=\sqrt{\ideal{x_1,\ldots,x_\ell}}\}.
  \]
  Here, $\sqrt{-}$ denotes the radical of the given ideal.
\end{dfn}

The arithmetic rank of an ideal has been of interest to algebraists
for as long as they have looked at ideals.  In a polynomial ring over
an algebraically closed field it answers the question by how many
hypersurfaces the affine 
variety defined by $I$ is cut out. The problem of finding this number
has a long history that is detailed excellently in
\cite{L-ara-survey,L-lc-survey}.  Some ground-breaking contributions
before the turn of the millennium included
\cite{Hartshorne-CDAV,EisenbudEvans-Inv73,Ogus-LCDAV,
  PeskineSzpiro-IHES73,Hartshorne-Bulletin74,HartshorneSpeiser-Annals77,
  Speiser-TAMS78,BrunsSchwanzl-BLMS90,Lyubeznik-AJM92}, and
\cite{Kunz-Intro} contains a gentle introduction to the problem.

\subsubsection{Some examples and conjectures}

Local cohomology is sensitive to arithmetic rank and relative
dimension. Indeed, it follows from the \v Cech complex point of view
that
\[
\max\{k\in\NN\mid H^k_I(A)\neq 0\}=\lcd_A(I)\le\ara_A(I),
\]
while a standard theorem in local cohomology asserts that
\[
\min\{k\in\NN\mid H^k_I(A)\neq 0\}=\depth_A(I,A),
\]
where $\depth_A(I,M)$ is the length of the longest $M$-regular sequence
in $I$. If $A$ is a Cohen--Macaulay ring, $\depth_A(I,A)$ is the height
of the ideal.

There are examples where the arithmetic rank exceeds the local
cohomological dimension, but it is often not easy to verify this since
the determination of $\lcd_A(I)$ and $\ara_A(I)$ is tricky.
\begin{exa}\label{exa-ara-not-lcd}
  \begin{asparaenum}
  \item Let $E$ be an elliptic curve over any field of characteristic $p>0$,
    and consider the Segre embedding 
    $E\times\PP^1_\KK\into\PP^5_\KK$. The curve $E$ is \emph{supersingular} if
    the Frobenius acts as zero on the one-dimensional space
    $H^1(E,\calO_E)$. It is known that if $E$ is defined over the integers
    then there are infinitely many $p$ for which the reduction $E_p$ is
    supersingular \cite{Elkies}, and infinitely many primes for which
    it is ordinary. For example, for $E=\Var(x^3+y^3+z^3)$,
    supersingularity is equivalent to $p-2$ being a multiple of $3$.  By
    \cite[Ex.~3]{HartshorneSpeiser-Annals77}, the local cohomological
    dimension of the ideal defining $E\times \PP^1_\KK$ in $\PP^5_\KK$
    equals three if and only if $E$ is supersingular (and it is $4$
    otherwise). However, by \cite{SinghWalther-Segre}, the arithmetic rank
    is always four, independently of supersingularity (and even in
    characteristic zero).
  \item Let $I\subseteq
    R=\CC[x_{11},x_{12},x_{13},x_{21},x_{22},x_{23}]$ be the ideal
    describing the image of the second Veronese map from $\PP^2$ to
    $\PP^5$ over the complex numbers. Then $\lcd_R(I)=3=\hight(I)$. On
    the other hand, as will be discussed in Example
    \ref{exa-Veronese}, the arithmetic rank of $I$ is $4$, not $3$.
    The underlying method, de Rham cohomology, is the topic of
    Subsection \ref{subsec-top-charzero}. Replacing de Rham arguments
    with \'etale cohomology, similar results hold in prime
    characteristic, Example \ref{exa-Veronese-p}. This is an example
    where the \'etale cohomological dimension of the projective
    complement $U$ surpasses the sum of the dimension and the
    cohomological dimension, $\ecd(U)=6>2+3=\cd(U)+\dim(U)$, compare
    also the discussion around Conjecture \ref{conj-L-ecd}.
  \end{asparaenum}
\end{exa}

Finding the arithmetic rank in concrete cases can be extremely
difficult; some of the long-standing open problems in this area
include general questions about ``large'' ambient spaces, but also about
concrete curves:
\begin{itemize}
\item Hartshorne's conjecture (\cite{Hartshorne-Bulletin74}: If
  $Y=\Proj(R/I)$ is a smooth $s$-dimensional subvariety of
  $\PP^n_\CC$, and $s>\frac{2n}{3}$, is then $Y$ a global complete
  intersection (\emph{i.e.}, is $Y$ the zero set of $\codim(Y)$ many
  projective hypersurfaces that, at each point of $Y$, are smooth and
  meet transversally)? 
  \item Is the \emph{Macaulay curve} in $\PP^3_\KK$, parameterized as
    $\{(s^4,s^3t,st^3,t^4)\}_{s/t\in\PP^1_\KK}$, a set-theoretic
    complete intersection (\emph{i.e.}, does the defining ideal have
    arithmetic rank 2, realized by homogeneous generators)?\\ This
    question is specific to characteristic zero, as in prime
    characteristic $p$, Hartshorne proved in \cite{Hartshorne-AJM79}
    that the Macaulay curve is a set-theoretic complete intersection
    for each $p$.
\end{itemize}
The (degree 5) Pl\"ucker embedding of the (6-dimensional) Grassmann
variety $\Gr_\CC(2,5)$ of affine $\CC$-planes in $\CC^5$ into
$\PP^9_\KK$ is not contained in a hyperplane, so Bezout's Theorem
indicates that we are not looking at a complete intersection. Thus,
the factor $2/3$ in the Hartshorne Conjecture is, in a weak sense,
optimal. Asymptotically, the coefficient must be at least $1/2$, but
Hartshorne writes in \cite{Hartshorne-Bulletin74}: ``I do not know any
infinite sequences of examples of noncomplete intersections which
would justify the fraction of the conjecture as
$n\rightarrow\infty$''. On the other hand, even less is known when
$\dim(Y)$ is small. For example, scores of articles have been devoted
to the study of monomial curves in $\PP^3_\CC$; in larger ambient
spaces \cite{Vogel-Monatsberichte71} contains a criterion for
estimating arithmetic rank in terms of ideal transforms, the functors
$\dlim_\ell\Ext^\bullet_A(I^\ell,-)$.

Some
of the major vanishing theorems in local cohomology came out of an
unsuccessful attempt to use local cohomology in order to show that
certain curves in $\PP^3_\KK$ cannot be defined set-theoretically by two
equations. For example, let $I\subseteq R=\KK[x_1,\ldots,x_4]$ define an
irreducible  projective
curve. In order for the arithmetic rank of $I$ to be $2$, $H^3_I(R)$
and $H^4_I(R)$ should both be zero. That $H^4_I(R)$ vanishes follows from
the Hartshorne--Lichtenbaum theorem. That $H^3_I(R)$ is also always zero
is the \emph{Second Vanishing Theorem} discussed in Subsection
\ref{subsec-vanishing}, in its incarnations due to  Ogus, Peskine--Szpiro,
Hartshorne, and Huneke--Lyubeznik.
In particular,
the desired obstruction to $\ara_A(I)=2$ cannot not materialize, but the
attempt led to the discovery of the Second Vanishing theorem.

On the positive side, Moh proved in \cite{Moh-MR0784166} that in
\emph{positive} characteristic every monomial curve in $\PP^3_\KK$ is
defined set-theoretically by two binomials; compare also
\cite{CowsikNori,Ferrand-MR0555692,Hartshorne-AJM79,BresinskyStueckradRenschuch,RoloffStueckrad}. The
construction of the two binomials uses heavily the Frobenius and, as
one might expect, the equations that work in one characteristic do not
work in another \cite{BarileMorales-MR1621700}. In characteristic
zero, Kneser proved that a curve in $\PP^3_\KK$ is cut out by three
equations if it has a $\KK$-rational point, and monomial space curves
are cut out by three binomials by \cite{BarileMorales-MR1621700}, but
nothing better is known at this point.  

\smallskip

There is recent progress on arithmetic rank and local
cohomological dimension in
toric and monomial situations.

In \cite{Varbaro-TAMS12}, Varbaro shows that if $X$ is a general
smooth hypersurface of projective $n$-space of degree less than $2n$
then the arithmetic rank of the natural embedding of the Segre product
of $X$ with a projective line is at most $2n$. This generalizes an
observation that appeared in \cite{SinghWalther-Segre} where $X$ is an
elliptic curve.  Moreover, Varbaro continues, if $X$ is a smooth conic
in the projective plane then its Segre product with projective
$m$-space has arithmetic rank exactly $3m$, as long as the characteristic is
not $2$.

Toric varieties, by which we
mean here the spectra of semigroup rings $\KK[S]$ where $S\subseteq
\ZZ^d$ is a finitely generated semigroup,
provide a standard testing ground for theories and conjectures.
Note that, for example, the
Macaulay curve falls into this category.

Barile and her coauthors have studied the question whether a toric
variety is a complete or almost complete intersection in
\cite{Barile-MR2222158,
  Barile-MR2244384,BarileMoralesThoma-MR1752767,
  BarileMoralesThoma-MR1896020}. Building on this,
\cite{Barile-MR2441609} shows that certain toric ideals of codimension
two are not complete intersections, and that their arithmetic rank is
equal to $3$.  The combinatorial condition with arithmetic flavor of
being $p$-glued has been shown to be pertinent here.
A semigroup can be $p$-glued for exactly one prime $p$,
\cite{BarileLyubeznik-PAMS05}. That such examples might be possible is
explained in part by the fact that the depth of the semigroup ring may
depend on the chosen characteristic: Hochster's theorem from
\cite{Hochster-Annals72} indicates for example how Cohen--Macaulayness
can toggle with $p$.

Monomial ideals and their local cohomology have been studied by
\`Alvarez-Montaner and his collaborators, see \cite{Montaner-LNM13} for notes to
a lecture series. At the heart of this work stands the
Galligo--Granger--Maisonobe correspondence between perverse sheaves
and hypercubes detailed in \cite{GalligoGrangerMaisonobe}.  Morally,
this is similar to the quiver encoding from Subsection
\ref{subsubsec-inv} and will receive a second look in Section
\ref{subsec-lambda}; compare specifically \cite{Alvarez-JSC05} on the
catogory of regular holonomic $\calD$-modules with support on a normal
crossing divisor and variation zero, and \cite{AlvarezGarciaZarzuela}.

In \cite{SchmittVogel} a technique is given how to find generators (up
to radical) for ideals that are intersections of ideals with given
generators. Application to monomial ideals 
relates to systematic search for the arithmetic rank of certain
intersections, compare \cite{Barile-MR2104196}.  Goresky and MacPherson noted in
\cite{GoreskyMacPherson-SMT,Jewell-Top94} a formula on the singular
cohomology of the complement of a complex subspace arrangement. The
article \cite{Yan-JPAA00} generalizes the formula to subspace
arrangements over any separably closed field using \'etale cohomology
and sheaf theory. These results are then applied to determine the
arithmetic rank of monomial ideals. In \cite{Yan-JA99}, Yan studies a
question of Lyubeznik on the arithmetic rank of certain resultant
systems and again uses \'etale cohomology to get some lower
bounds. More recently, Kimura and her collaborators have produced a
wealth of new information on arithmetic rank of monomial
ideals, \emph{cf.} \cite{KimuraMontero-JCA17} and its bibliography tree.

\subsubsection{Endomorphisms of local cohomology}

As always, $(A,\frakm,\kk)$ is a Noetherian local ring and $\fraka$ an
ideal of $A$.  In this subsection we discuss some 
challenges that have arisen in the last two decades, connecting the
question of finding the arithmetic rank to problems about
$\calD$-modules, with the focus on the question of determining whether
a given ideal be a complete intersection.

We recall that the local cohomologcal dimension $\lcd_A(\fraka)$ is a
lower bound for the arithmetic rank $\ara_A(\fraka)$ and that the two
invariants may not be equal, Example \ref{exa-ara-not-lcd}. Nonetheless, as work primarily
of Hellus and St\"uckrad shows, local cohomology modules contain
information that can lead to the determination of arithmetic
rank. However, decoding it successfully is at this point a serious challenge.

The story starts with a result of Hellus from
\cite{Hellus-CiA05}. Denote $E=E_A(\kk)$ the injective hull of the
residue field. Suppose $f_1,\ldots, f_c$ are elements of $\fraka$, and
write for simplicity $\frakb_i$ for the $A$-ideal generated by
$f_1,\ldots,f_i$. Assuming that $\lcd_A(\fraka)=c$, Hellus showed that
these elements
generate $\fraka$ up to radical if and only if $f_i$ operates
surjectively on $H^{c+1-i}_\fraka(A/\frakb_{i-1})$ for $1\le i\le
c$. This has the following corollary pertaining  to 
set-theoretic complete intersections: if $f_1,\ldots,f_c$ is an
$A$-regular sequence (in our situation this means that
$H^i_\fraka(A)=0$ unless $i=c$), then the sequence generates $\fraka$ up to
radical if and only if they form a regular sequence on
$D(H^c_\fraka(A))$ where, as before,
\[
D(M):=\Hom_A(M,E)
\]
is the Matlis dual. This is discussed from a new angle in \cite{HartshornePolini-20}

This motivates (when only one $H^i_\fraka(A)$ is nonzero) the study of
the multiplication operators $f_i\colon D(H^c_\fraka(A))\to
D(H^c_\fraka(A))$, and in particular the associated primes of
$D(H^c_\fraka(A))$. In fact, Hellus offers the following conjecture:
if $(A,\frakm,\kk)$ is local Noetherian, 
\begin{eqnarray}\label{eqn-Hellus-cnj}
\text{Is }  \Ass_A(D(H^i_{\frakb_i}(A)))=\{\frakp\in\Spec A\mid
  H^i_{\frakb_i}(A/\frakp)\neq 0\} \,\,?
\end{eqnarray}
(One always has the inclusion $\subseteq$ above, and in the
equi-characteristic case, the set $\{\frakp\in\Spec(A)\mid
f_1,\ldots,f_i\text{ is part of an s.o.p.\ for }A/\frakp\}$ is contained in
$\Ass_A(D(H^i_{\frakb_i}(A)))$---but this may not be an
equality. In mixed characteristic, a similar statement can be
made). Hellus proceeds to show that this conjecture is equivalent to
$\Ass_A(D(H^i_{\frakb_i}(A)))$ being stable under
\emph{generalization}, and also gives 
the following reformulation:
\begin{prb}\label{prb-Hellus}
  For all Noetherian local \emph{domains}
  $(A,\frakm,\kk)$ and for all $f_1,\ldots,f_c\in A$, show that the
  nonvanishing of $H^i_{\ideal{f_1,\ldots,f_i}}(A)$ implies that the
  zero ideal is associated to $D(H^i_{\ideal{f_1,\ldots,f_i}}(A))$.
\end{prb}

\begin{rmk}
A significant part of Problem \ref{prb-Hellus} was resolved positively
in \cite{LYildirim-PAMS18}. Namely, if $R$ is a regular Noetherian
local ring of prime characteristic, then $\Ass_R(D(H^i_I(R)))$ contains
$\{0\}$, as long as $H^i_I(R)$ is nonzero. In fact, it is shown for
all $F$-finite $F$-modules $M$ that $\{0\}$ has to be associated to at
least one of $M, D(M)$. The proof is an explicit construction of an
element that is not torsion.

{Motivated by their work in prime characteristic, they conjectured in \cite[Conjecture 1]{LYildirim-PAMS18} that, if $(R,\frakm)$ is a regular local ring and $I$ is an ideal such that $H^i_I(R)\neq 0$, then $(0)\in \Ass_R(D(H^i_I(R)))$.}
\end{rmk}

\begin{rmk}
Let $R=\ZZ_2[[x_0,\dots,x_5]]$ and let $I$ be the monomial ideal as in
Example \ref{Reisner example}. It follows from \cite[Remark
  5.3]{DattaSwitalaZhang} that the arithmetic rank of $I$ is 4;
equivalently there are $f_1,\dots,f_4\in R$ such that
$H^4_I(R)=H^4_{\ideal{f_1,\dots,f_4}}(R)$. By \cite[Proposition 5.5]{DattaSwitalaZhang}, $H^4_I(R)\cong E_{\bar{R}}(\bar{R}/\frakm))$, where $\bar{R}=R/\ideal{2}$ and $\frakm=\ideal{2,x_0,\dots,x_5}$. Hence 
\[D(H^4_{\ideal{f_1,\dots,f_4}}(R))=D(H^4_I(R))\cong \bar{R}.\] 
Consequently the zero ideal is not associated to
$D(H^4_{(f_1,\dots,f_4)}(R))$. This answers Hellus' question in Problem \ref{prb-Hellus} in the negative for unramified regular local rings of mixed characteristic, and provides a counterexample to the conjecture of Lyubeznik and Yildirim in mixed characteristic.
\end{rmk}

In \cite{Hellus-CiA07Matlis}, an example is given
where arithmetic rank and local cohomological dimension differ. What
is special here is that $\lcd_A(\fraka)=1$; Hellus gives a
criterion for the arithmetic rank to be one, based on
the \emph{prime avoidance property} of $\Ass_A(D(H^1_\fraka(A)))$. In the
same year and journal \cite{Hellus-CiA07dual}, he shows for
Cohen--Macaulay rings the curious identity
$H^c_\fraka(D(H^c_\fraka(R)))=D(R)$, provided that $c=\lcd_A(\fraka)$
is also the grade of $\fraka$. This was subsequently generalized in
\cite{Khashyarmanesh-Basel07}.

In \cite{HellusStueckrad-PAMS08duals}, Hellus and St\"uckrad continue their
study of associated primes of, and regular sequences on,
$D(H^c_\fraka(A))$. They show that $H^m_{\ideal{f_1,\ldots,f_m}}(A)$
always surjects onto
$H^{m+n}_{\ideal{f_1,\ldots,f_m,g_1,\ldots,g_n}}(A)$ for $m>0$ and
derive from this some insights about the inclusion \eqref{eqn-Hellus-cnj},
and about Problem \ref{prb-Hellus} when $A$ is a complete domain and
$\fraka$ a 1-dimensional prime.  In \cite{HellusStueckrad-PAMS08endo}
the same authors show that
in a complete local ring, when $\fraka$ has the local cohomological
behavior of a complete interection (\emph{i.e.}, $H^i_\fraka(A)=0$ unless $i=c$),
then the natural map $A\to \End_A(H^c_\fraka(A))$ is an
isomorphism. (In general, this map is not surjective and has a kernel). 
In particular, no element of $A$ annihilates $H^c_\fraka(A)$. 
By results mentioned above, this means that if $\fraka$ behaves local
cohomologically like a complete intersection and if $f_1,\ldots,f_c$ is
an $A$-regular sequence in $\fraka$, then
$D(H^c_\fraka(D(H^c_\fraka(A))))$ is an ideal of $A$ which, if computable,
predicts whether $\fraka$ is a complete intersection. {For more on $\End_A(H^c_\fraka(A))$, see \cite{SchenzelPAMS2009, Khashyarmanesh-Canadian10, SchenzelArchMath2010, SchenzelJA2011}.}

In \cite{HellusSchenzel-JA08} it is investigated which ideals behave
like a complete intersection from the point of local cohomology, by
establishing relations to iterated local cohomology functors which
then lead to Lyubeznik numbers (see Section \ref{subsec-lambda}). For
example, if $\fraka=\ideal{f_1,\ldots,f_c}$ is an ideal of dimension
$d$ in a local Gorenstein ring, and if $\fraka$ is a complete
intersection outside the maximal ideal, then [$H^i_\fraka(A)=0$ unless $i=c$]
precisely when [$H^d_\frakm(H^c_\fraka(A))=E_A(\kk)$ and
  $H^i_\frakm(H^c_\fraka(A))=0$ for $i\neq d$]. In particular, the
complete intersection property of $\fraka$ is then completely detectable
from $H^c_\fraka(A)$ alone. A new version of some of these ideas is given in a
recent work of Hartshorne and Polini, who introduce and investigate
coregular sequences and codepth in \cite{HartshornePolini-20}.

\subsection{Relation with de Rham and \'etale cohomology}
\label{subsec-top-charzero}

\subsubsection{The \v Cech--de Rham complex}
\label{subsubsec-top-deRham}

Suppose $I\subseteq R_\KK:=\KK[x_1,\ldots,x_n]$ is generated by
$f_1,\ldots,f_m$  and assume that $\KK$ is a field containing $\QQ$. The finitely many
coefficients of $f_1,\ldots,f_m$ all lie in some finite extension field $\kk$ of
$\QQ$, and because of flatness one has
$H^i_I(\KK[x_1,\ldots,x_n])=H^i_{I_\kk}(\kk[x_1,\ldots,x_n])\otimes_{\kk}\KK$,
where $I_\kk=\ideal{f_1,\ldots,f_m}R_\kk=I\cap R_\kk$ with 
$R_\kk=\kk[x_1,\ldots,x_n]$.

The finite extension $\kk$ can be embedded into $\CC$ and then, by
flatness again,
$H^i_I(\CC[x_1,\ldots,x_n])=H^i_{I_\kk}(\kk[x_1,\ldots,x_n])\otimes_{\kk}\CC$. It
follows that most aspects of the behavior of local cohomology in
characteristic zero can be studied over the complex numbers.
\begin{cnv}
  In this subsection, $\kk=\CC$ and $I$ is an ideal of
  $R=\CC[x_1,\ldots,x_n]$.  The advantage of working over $\CC$ is
  that one has access to topological notions and tools.
\end{cnv}

The arithmetic rank of the ideal $I$
is the
smallest number of principal open affine sets $U_f$ that cover the
complement $U_I=\CC^n\setminus \Var(I)$. Any $U_f$ arises also as the
closed affine set defined by $f\cdot x_0=1$ inside $\CC^n\times\CC^1_{x_0}$.

Complex affine space as well as all its Zariski closed subsets are
\emph{Stein spaces}. 
This is a
complex analytic condition that includes separatedness by holomorphic
functions, and a convexity condition about compact sets under
holomorphic functions. It implies, among other things, that a Stein
space of complex dimension $n$ has the homotopy type of an
$n$-dimensional CW-complex.
In
particular, a Stein space $S$ of complex dimension $n$ cannot have
singular cohomology $H^i_\Sing(S;-)$ beyond degree $n$.
That complex affine varieties have this
latter property is the \emph{Andreotti--Frankel Theorem}.  (For example, a
Riemann surface is Stein exactly when it is not compact).
In the ``spirit of GAGA'', \cite{Serre-FAC,Serre-GAGA}, Stein spaces
are the notion that corresponds to affine varieties.

Now consider the complement $U_I=U_{f_1}\cup\cdots\cup U_{f_m}$ of the
variety $\Var(I)$. It follows from the
Mayer--Vietoris principle that $H^i_\Sing(U;-)=0$ for all
$i>n+m-1$ and all coefficients. Being Stein is not a local property:
\begin{exa}
  Let $I=\ideal{x,y}\subseteq \CC[x,y]$. Then $U_I$ is homotopy
  equivalent to the $3$-sphere and in particular cannot be Stein.
\end{exa}

\medskip

The \v Cech complex on a set of generators for $I$ is always a
complex in the category of $\calD$-modules.
Let $\varphi\colon X\to Y$
be a morphism of smooth algebraic varieties. We refer to \cite{HTT}
for background and details on the following continuation of
the discussion on functors on $D$-modules in Section \ref{sec-DF}.

There are (both regular and exceptional) direct and inverse image
functors between the categories of bounded complexes of $\calD$-modules on
$X$ and $Y$. These functors preserve the categories of complexes with
holonomic cohomology. In particular, one can apply them to the
structure sheaf, or to local cohomology modules and \v Cech
complexes.

If $\iota\colon U\into X$ is an open embedding and $M$ a
$\calD_U$-module, then the direct image of $M$ under $\iota$ as
$\calD$-module agrees with the direct image as $\calO$-module.
For example, in both categories there is an exact triangle
\[
\RGamma_{X\setminus U}(-)\to \id\to \iota_*((-)|_U)\stackrel{+1}{\to}.
\]
Let $X=\CC^n$ and choose $\varphi\colon X\to Y$ be the projection to a
point $Y$. Write
\[
\omega_X=\calD_X/\ideal{\del_1,\ldots\del_n}\cdot\calD_X;
\]
this gives the canonical sheaf of the manifold $X$ a right
$\calD_X$-structure in a functorial way.  Then under $\iota\colon
U\into X$, $\calO_U$ turns into a complex of sheaves that is
represented on global sections by the \v Cech complex on generators of
the ideal $I=\ideal{f_1,\ldots,f_m}$ describing $X\setminus U$. The
$\calD$-module direct image under $\varphi$ corresponds to the functor
$\omega_X\otimes^L_{\calD_X}(-)$ whose output is a complex of vector
spaces.  Applying this functor to the \v Cech complex for $I$ invites
the inspection of a \emph{\v Cech--de Rham spectral sequence} starting
with $\Tor_\bullet^{\calD_X}(\omega_X,H^\bullet_I(\calO_X))$.  With
$R=\Gamma(X,\calO_X)$, $\omega_R=\Gamma(X,\omega_X)$, and
$D=\Gamma(X,\calD_X)$, the Grothendieck Comparison Theorem
\cite{Grothendieck-deRham} asserts that on global sections, the
abutment of the spectral sequence is the reduced de Rham cohomology of
$U$,
\begin{gather}\label{eqn-Cech-deRham}
  E_2^{i,j}=\Tor_{n-j}^{D}(\omega_R,H^i_I(R))\Rightarrow
  \tilde H^{i+j-1}_\dR(U;\CC).
\end{gather}
We note in passing that there are algorithmic methods that can compute
the pages of this spectral sequence as vector spaces over $\CC$, see
\cite{OT1,W-algdR00,W-alg-det01,OT2}.  In the sequence
\eqref{eqn-Cech-deRham}, the Tor-groups involved vanish for the index
exceeding $\dim X$, and so the spectral sequence operates clearly
inside the rectangle $0\le i\le \lcd_R(I)$, $0\le j\le n$. However, it
is actually limited to a much smaller, triangular region, compare
\cite{RWZ-lambda}.

This now opens the door to direct comparisons between local cohomology
groups of high index and singular cohomology groups of high index; the
de Rham type arguments in the following example are written down in
\cite{LSW,HartshornePolini-simple}, but are folklore and were known to
the authors of \cite{Ogus-LCDAV} and \cite{Hartshorne-DRCAV}. For
example, Theorem 2.8 in \cite{Ogus-LCDAV} shows that in a regular
local ring $R$ over $\QQ$ with closed point $\frakp$, the vanishing of
local cohomology $H^j_I(R)$ for all $j>r$ implies the vanishing of the
local de Rham cohomology groups $H^i_\frakp(\Spec(R/I))$ for all
$i<\dim(R)-r$ (and is in fact equivalent to it if one already knows
that the support of $H^j_I(R)$ is inside $\frakp$ for $j>r$).

\begin{exa}\label{exa-2x3-betti}
  We continue Example \ref{exa-2x3} with $\KK=\CC$. The open set
  $U=\CC^6\setminus\Var(I)$ consists of  the set of
  $2\times 3$ complex matrices of rank two. The closed set $V=\Var(I)$
  is smooth outside the origin, as one sees from the
  $GL(2,\CC)$-action. Since $\dim(R/I)=4$, the height of $I$ is $2$
  and so $H^{2+1}_I(R)$ must be, if nonzero, supported at the origin
  only, by Remark \ref{rmk-other-powers}. 

  Since $H^3_I(R)$ is also a holonomic $D$-module, $D$ being the
  ring of $\CC$-linear differential operators on $R$, Kashiwara
  equivalence (\cite[\S1.1.6]{HTT}) asserts that $H^3_I(R)$ is a
  finite direct sum of $\lambda$ copies of $E$, the $R$-injective hull
  of the residue field at the origin. The number $\lambda$ can be
  evaluated as follows.

  Since $I$ is $3$-generated, $H^{>3}_I(R)=0$ and the \v Cech--de Rham
  spectral sequence shows that $H^i(U;\CC)$ vanishes for
  $i>6+3-1=8$. Moreover, an easy exercise shows that
  $\Tor_{n-j}^D(\omega_R,E)=0$ unless $j=0$, and in that case
  returns one copy of $\CC$ so that the only possibly nonzero
  $E_2$-entry in the spectral sequence \eqref{eqn-Cech-deRham} in
  column $3$ is the entry $E_2^{3,6}=\CC^\lambda$. The workings of the
  spectral sequence make it clear that all differentials into and out
  of position $(3,6)$ on all pages numbered $2$ and up vanish. So,
  $\CC^\lambda=E^{3,6}_2=E^{3,6}_\infty=H^8(U;\CC)$. We now compute
  this group explicitly via the following argument taken from Mel
  Hochster's unpublished notes on local cohomology.

  Let $A$ be a point of $U$, representing a rank two $2\times 3$
  matrix. Consider the deformation that scales the top row to length
  $1$, followed by the deformation (based on gradual row reduction)
  that makes the bottom row perpendicular to the top row and then
  scales it to length $1$ as well. Then the top row varies in the
  $5$-sphere, and for each fixed top row the bottom row varies in a
  $3$-sphere. Let $M$ be this retract of $U$ and note that, projecting
  to the top row, it is the total space of an $S^3$-bundle over
  $S^5$. Both base and fiber are orientable, and the base is simply
  connected. Thus, $M$ is an orientable compact manifold of dimension
  $8$ which forces $1=\dim_\CC H^8(M;\CC)=\dim_\CC
  H^8(U;\CC)=\lambda$.
\end{exa}  

\begin{rmk}\label{rmk-2x3}
  Already Ogus proved in  \cite{Ogus-LCDAV} results that relate the local
  cohomology module $H^3_I(R)$ of Example \ref{exa-2x3-betti} to
  topological information. We discuss this in and after Theorem
  \ref{thm-Ogus-dRdepth} below. In brief, the non-vanishing of
  $H^3_I(R)$ is ``to be blamed'' on the failure of the restriction map
  $H^2_\dR(\PP^5_\CC)\to H^2_\dR(Y)$ to be surjective. Here, $Y$ is
  the image of the Segre map and
  $\dim_\CC(H^2_\dR(Y))=\dim_\CC(H^2_\dR(\PP^1_\KK\times
  \PP^2_\KK))=2$ by the K\"unneth theorem.
\end{rmk}

\begin{exa}[{Compare \cite[Exa.~4.6]{Ogus-LCDAV}}]\label{exa-Veronese}
  Let $\iota\colon \PP^2_\CC\into \PP^5_\CC$ be the second Veronese
  morphism, denote the target by $X$, the image by $Z$ and write
  $U:=X\setminus Z$. There is a long exact sequence of singular (local)
  cohomology
  \[
  H^p_Z(X;-)\to H^p(X;-)\to H^p(U;-)\stackrel{+1}{\to}
  \]
  and a natural identification $H^p_Z(X;-)\cong (H^{2\dim
    X-p}_c(Z;-))^\vee$ with compactly supported cohomology, for any
  coefficient field, compare \cite[\S6.6]{Iversen}. Via Poincar\'e
  duality, this allows to identify the map $H^p_Z(X;-)\to H^p(X;-)$ as
  the dual to $H_{2\dim X-p}(Z;-)\to H_{2\dim X-p}(X;-)$. Now take
  $\ZZ/2\ZZ$ as coefficients. Then, since $\iota^*$ sends the generator
  of $H^2(X;\ZZ)$ to twice the generator of $H^2(Z;\ZZ)$, the
  long exact sequence shows that $H^8(U;\ZZ/2\ZZ)$ is nonzero. Thus,
  $U$ cannot be covered by three affine sets and $\ara_A(I)\geq 4$. (In
  fact, $\ara_A(I)=4$ as one finds easily from experiments).
\end{exa}

\subsubsection{Algebraic de Rham cohomology}
\label{subsubsec-alg-dR}

In \cite{Hartshorne-DRCAV}, Hartshorne defines and develops for
(possibly singular) schemes over a field of characteristic zero a
purely algebraic (co)homology theory that he connects to singular
cohomology via comparison theorems. In a nutshell, the de Rham
cohomology $H^q_{\dR}(Y)$ of $Y$ embedded into a smooth scheme $X$ is
the $q$-th hypercohomology on $X$ of the de Rham complex on $X$,
completed along $Y$. Similarly, the de Rham homology $H_q(Y)$ of $Y$
is the $(2\dim(X)-q)$-th local hypercohomology group with support in
$Y$ of the de Rham complex on $X$. (We add here a pointer to Remark
\ref{rmk-GreenleesMay}). Hartshorne develops many tools of
singular (co)homology: Mayer--Vietoris sequences, Thom--Gysin sequences,
Poincar\'e duality, and a local (relative) version. With it, he shows
foundational finiteness as well as Lefschetz type theorems.

One of the most remarkable applications of his theory as it relates to
local cohomology is worked out in
the thesis of Ogus, and based on the following definition.
\begin{dfn}[{\cite[Dfn.~2.12]{Ogus-LCDAV}}]
 Let $Y$ be a scheme over a field of characteristic zero. The \emph{de
   Rham depth} $\dRdepth(Y)$ of $Y$ is the greatest integer $d$ such
 that for every point $\fraky\in Y$ (closed or not) one has
 \[
 H^i_\fraky(Y)=0\qquad\text{ for } i<d-\dim(\overline{ \{y\}}).
 \]
\end{dfn}
This number never exceeds the dimension of $Y$ as one sees by looking
at a closed point $\fraky$. Ogus uses it in the following fundamental
result; we point here at Remark \ref{rmk-Fdepth} for the corresponding result in
positive characteristic and note the formal similarities both of de Rham and
$F$-depth, and the corresponding results on local cohomological
dimension.

\begin{thm}[{\cite[Thm.~2.13]{Ogus-LCDAV}}]\label{thm-Ogus-dRdepth}
  If $Y$ is a closed subset of a smooth Noetherian scheme $X$ of
  dimension $n$ over a field $\kk$ of characteristic zero, then for each
  $d\in \NN$ one has
  \[
    [\lcd(X,Y)\le n-d] \Leftrightarrow [\dRdepth(Y)\geq d].
    \]
    In particular, if $Y=\Spec(R/I)$ for some regular $\kk$-algebra $R$ then
    $n-\lcd_R(I)=\dRdepth(Y)$ is intrinsic to $Y$ and does not depend
    on $X$.
\end{thm}

Now let $Y$ be a projective variety over the field $\kk$ of
characteristic zero, embedded into $\PP^n_\kk$. Let $R$ be the
coordinate ring of $\PP^n_\kk$ and $I$ the defining ideal of $Y$; 
of course, these are not determined by $Y$.
Then Ogus obtains in
\cite[Thm.~4.1]{Ogus-LCDAV} the equivalences
\begin{eqnarray*}
  [\lcd(\PP^n_\kk,Y)\le r]&\Leftrightarrow&
  [\Supp_R(H^i_I(R))\subseteq \frakm\text{ for }i>r]\\
  &\Leftrightarrow&
  [\dRdepth(Y)\geq n-r].
\end{eqnarray*}
In particular, for any such embedding, the smallest integer $r$ such
that $H^{>r}_I(R)$ is Artinian is intrinsic to $Y$.

One might wonder whether a similar result holds for $\lcd(R,I)$
itself.
With the same notations as in the previous theorem, Ogus proves in
\cite[Thm.~4.4]{Ogus-LCDAV}: 
\[
  [\cd(\PP^n_\kk\smallsetminus Y)< r]
    \]
  (that is, $\lcd(R,I)\le r$) 
is equivalent to
  \[
  \left[ [\dRdepth(Y)\geq n-r]\text{ and }
  [H^i_{\dR}(\PP^n_\kk)\onto H^i_{\dR}(Y)\text{ for }i<n-r]\right].
\]
Note that these restriction maps are always injective, and
surjectivity is preserved under Veronese maps.

\subsubsection{Lefschetz and Barth Theorems}
\label{subsubsec-top-Barth}

Let $X\subseteq \PP^n_\CC$ be a projective variety and $H\subseteq
\PP^n_\CC$ a hyperplane. Setting $Y=X\cap H$, the \emph{Lefschetz
  hyperplane theorem} states that under suitable hypotheses the
natural restriction map
\begin{gather}
\label{eqn-rho}
\rho_{X,Y}^i\colon H^i(X;\CC)\to H^i(Y;\CC)
\end{gather}
is an isomorphism for $i<\dim(Y)$ and injective for $i=\dim(Y)$. In
the original formulation by Lefschetz, $X$ is supposed to be smooth
and $H$ should be generic (which then entails $Y$ being
smooth). Inspection showed that the relevant condition is that
the affine scheme $X\setminus Y$ be smooth, since then the relative groups
$H^i_\Sing(X,Y;\CC)$ are zero in the required range.

It is clear that one can iterate this procedure and derive similar
connections between the cohomology of $X$ and the cohomology of
complete intersections on $X$ that are well-positioned with respect to
the singularities of $X$. (Recall that any hypersurface section can be
cast as a hyperplane section via a suitable Veronese embedding of
$X$).

A rather more difficult problem is to establish connections when $Y$
is not a complete intersection. At the heart of the problem is the
issue that in general $X\setminus Y$ will not be affine and thus might
allow
more complicated cohomology. 

In \cite{Barth-AJM70}, Barth developed theorems that connect, for
$Y\subseteq \PP^n_\CC$ smooth (and of small codimension), the
surjectivity of $\rho_{\PP^n_\CC,Y}^i$ to the surjectivity of corresponding
restrictions $\rho_{\PP^n_\CC,Y}^i(\calF)$ of coherent sheaves $\calF$
and hence to the cohomological dimension of $\PP^n_\CC\setminus Y$ and
the arithmetic rank of the defining ideal of $Y$. More
precisely, he proved that surjectivity of
$\rho_{\PP^n_\CC,Y}^i(\calF)$ occurs for $i\le 2\dim(Y)-n$ and proved
for $\calF=\calO_{\PP^n_\CC}$ that surjectivity of
$\rho_{\PP^n_\CC,Y}^i(\calF)$ is equivalent to surjectivity of
$\rho_{\PP^n_\CC,Y}^i$ in the sense of Equation \eqref{eqn-rho}
above.
As a corollary, he obtained a more general
form of the Lefschetz Hyperplane Theorem: if $\PP^n_\CC\supseteq X, Y$
are smooth with $\dim(X)=a$, $\dim(Y)=b$ then
\[
\rho_{Y,X\cap Y}^i\colon H^i(Y;\CC)\to H^i(X\cap Y;\CC)
\]
is an isomorphism for $i\le \min(2b-n,a+b-n-1)$.
It is worth looking at the special case when $X$ is the ambient
projective space.
For $i=0$
the theorem then 
generalizes the fundamental fact that a smooth subvariety of pure
dimension $a$ is connected whenever $2a\geq n$. But it also gives
obstructions for embedding varieties into projective spaces of given
dimension, since it forces the singular cohomology groups $H^i(Y;\CC)$ to
agree with those of $\PP^n_\CC$ in the range $i\le 2\dim(Y)-n$.
For example, an Abelian variety $Y$ of dimension $b$ cannot be
embedded in $Y=\PP^{2b-1}_\CC$ since with such embedding the map
$H^1(\PP^{2b-1};\CC)\to H^1(Y;\CC)$ should be surjective.

Barth uses the special unitary group action on $\PP^n_\CC$ to
``spread'' the classes on $Y$ to classes on $\PP^n_\CC$ near $Y$. In
order to glue them, he then needs a suitable cohomological triviality
of the complement of $Y$.  In \cite{Ogus-LCDAV}, Ogus gives an
algebraic version of Barth's transplanting technique, and succeeds (in
his section 4) in
proving
various statements that connect the isomorphy of the restriction maps
of de Rham cohomology of two schemes $X\subseteq Y\subseteq \PP^n_\kk$
to the de Rham depths of $X$, $Y$ and $X\smallsetminus Y$.

In \cite{Speiser-TAMS78}, Speiser studies in
varying characteristics the
cohomological dimension of the complement $C_Y$ of the diagonal in
$Y\times Y$. As a stepping stone he studies $C_{\PP^n_\KK}$ for
arbitrary fields. In any characteristic, the diagonal
scheme is the set-theoretic intersection of $2n-1$ very ample
divisors. However, a big difference appears for cohomological dimension:
$\cd(C_{\PP^n_\KK})=2n-2$ when $\QQ\subseteq \KK$, but 
$\cd(C_{\PP^n_\KK})=n-1$ in positive characteristic. The discrepancy
is due to the Peskine--Szpiro Vanishing
since the diagonal comes with a Cohen--Macaulay
coordinate ring.

In characteristic zero, Speiser's results imply that the diagonal of
projective space is cut out set-theoretically by $2n-1$ and no fewer
hypersurfaces. More generally, for Cohen--Macaulay $Y$, he shows
in \cite[Thm.~3.3.1]{Speiser-TAMS78} a similar vanishing result about
$C_Y$ in positive characteristic over algebraically closed fields: the
cohomological dimension of $Y\times Y\setminus \Delta$ is bounded by
$2n-2$ whenever $Y\subseteq \PP^n$ is a Cohen--Macaulay scheme of
dimension $s\geq(n+1)/2$.

\subsubsection{Results via \'etale cohomology}
\label{subsubsec-top-etale}

Suppose $U$ is an open subset of affine space $X=\CC^n$ whose closed
complement $V=X\setminus U$ is defined by an ideal $I$ in the appropriate
polynomial ring $R$.  We have seen in \eqref{eqn-Cech-deRham}
that the local cohomological dimension $\lcd_R(I)$ is related to the de
Rham cohomology via the vanishing
\begin{gather}\label{eqn-lcd-dR}
  [H^i_\dR(U;\CC)=0] \text{ whenever }
  [i\geq\lcd(I)+n-1=\cd(U)+n].
\end{gather}
We mention here a variant of this in arbitrary
characteristic, involving \'etale cohomology. This is a cohomology
theory that interweaves topological data with arithmetic
information. We refer to \cite{Milne,Milne-notes} for guidance on
\'etale cohomology.

One significant difference to the de Rham case is that the basic
version of \'etale cohomology involves coefficients that are torsion
(\emph{i.e.}, sheaves with stalk $\ZZ/\ell\ZZ$) of order not divisible
by $p=\ch(\kk)$. 

In many aspects, over a separably closed field $\kk$, \'etale
cohomology behaves quite similar to de Rham or singular cohomology
over the complex numbers. For example, on non-singular projective
varieties there is a version of Poincar\'e duality, there is a
K\"unneth theorem, and if a variety is defined over $\ZZ$ then its
model over $\CC$ has singular cohomology group ranks equal to the
corresponding \'etale cohomology ranks of the reductions modulo $p$
for most primes $p$.

The \'etale cohomology groups on a scheme $X$ vanish beyond $2\dim X$,
and even beyond $\dim(X)$ if $X$ is affine, similar to the
Andreotti--Frankel Theorem. So, it makes sense to talk of \'etale
cohomological dimension $\ecd(-)$, the largest index of a
non-vanishing \'etale cohomology group. The Mayer--Vietoris principle
implies that if $V$ is a variety inside affine $n$-space $X\neq V$
over the algebraically closed field $\kk$, cut out by the ideal $I$,
then with $U=X\setminus V$ one has
\begin{gather}\label{eqn-ara-and-ecd}
  \ecd(U)\le n+\ara_A(I)-1.
\end{gather}
Note that $\ara_A(I)\geq\lcd_R(I)=\cd(U)+1$. 

In \cite{L-lc-survey}, Lyubeznik formulates the following conjecture.
\begin{conj}\label{conj-L-ecd}
  Over a separably closed field $\kk$,
\[
\ecd(U)\geq  \dim(U)+\cd(U).
\]
\end{conj}
In this conjecture, $U$ need not be the complement of an affine
variety or even smooth.  Comparison with \eqref{eqn-lcd-dR} shows that
(for complements of varieties in affine or projective spaces) the
conjecture can be interpreted to say that \'etale cohomology always
provides a better lower bound for arithmetic rank than local
cohomological dimension does.  At present, this conjecture seems wide
open.  Varbaro shows in \cite{Varbaro-TAMS12} that it holds over $\CC$
in the case that $U$ is the complement in projective space
$\PP^n_\CC\setminus V$ of a smooth variety $V$
with $\cd(\PP^n\setminus V) > \codim_{\PP^n}(V)-1$.

\begin{exa}\label{exa-Veronese-p}
  We continue Example \ref{exa-Veronese}. For $\KK=\CC$ and
  all other field coefficients of characteristic not equal to 2, one has
  $H^8(U;\KK)=0$. Thus, we cannot conclude that
  $\lcd_R(I)\geq 4$ in the way we
  concluded in Example \ref{exa-2x3-betti}.
  In fact, as Ogus \cite[Exa.\ 4.6]{Ogus-LCDAV}
  proved, $\cd(U)=2$ (and so $\lcd_R(I)=3$) and, in particular,
  $\ecd(U)>\dim(U)+\cd(U)$. 

  In finite characteristic different from 2, if one replaces
  ``singular'' by ``\'etale'', the same formal arguments as in Example
  \ref{exa-Veronese} show that the arithmetic rank of the defining
  ideal of the Pl\"ucker embedding is $4$ while (since the coordinate
  ring is Cohen--Macaulay) $\lcd_R(I)=3$.

  In characteristic $2$, the arithmetic rank drops to $3$ and the
  ideal is generated up to radical by
  $\{t_{xx}t_{yy}-t_{xy}t_{xy},t_{xx}t_{zz}-t_{xz}t_{xz},
  t_{yy}t_{zz}-t_{yz}t_{yz}\}$ since, for example,
  $t_{xx}t_{zz}(t_{xx}t_{yy}-t_{xy}t_{xy})+
  t_{xy}t_{xy}(t_{xx}t_{zz}-t_{xz}t_{xz})+t_{xx}t_{xx}(t_{yy}t_{zz}-
  t_{yz}t_{yz})=(t_{xx}t_{yz}-t_{xy}t_{xz})^2$ in characteristic $2$.
\end{exa}
\begin{rmk}
  In \cite[Rmk.~2.13]{Varbaro-TAMS12}, Varbaro points out that Example
  \ref{exa-Veronese-p} shows that the \'etale cohomological dimension
  of the complement of an embedding of $\PP^2_\kk$ into $\PP^5_\kk$
  depends on the embedding: for a subspace embedding it is at most
  $3+4$ since the subspace is covered by three affine spaces of
  dimension $5$, but for the Veronese it is $8$ (compare also
  \cite{Barile-JA95} for arithmetic rank consequences that highlight
  variable behavior in varying characteristic).  This contrasts with
  his Theorem 2.4, which states that the quasi-coherent cohomological
  dimension is independent of the embedding (intrinsic to the given
  smooth projective subvariety).
  \item Ogus proved in \cite[Ex.~4.6]{Ogus-LCDAV} for any Veronese map of a
    projective space in characteristic zero that the local cohomological
    dimension agrees with the height of the defining ideal.  In
    positive characteristic, the same follows from Peskine--Szpiro
    \cite[Prop.~III.4.1]{PeskineSzpiro-IHES73}. In \cite{Pandey-2005.03250},
    Pandey shows that this is even true over the integers, and by
    extension then over every commutative Noetherian ring. 
\end{rmk}

Now, recall Speiser's result from Subsection
\ref{subsubsec-top-Barth}, on the arithmetic rank $2n-1$ of the
diagonal of $\PP^n_\KK\times \PP^n_\KK$.  In \cite{Varbaro-TAMS12}
Varbaro shows that it remains true in every characteristic as long as
$\KK$ is separably closed; note, however, that the cohomological
dimension of the complement is much smaller in finite characteristic,
always equal to $n-1$.  The main ingredient comes from K\"unneth
theorems on \'etale cohomology.

There are Lefschetz and Barth type results for \'etale cohomology. For
example, in \cite[Prop 9.1]{L-ecd}, Lyubeznik proves the following:
assume $\KK$ to be separably closed, of any characteristic, and pick
two varieties $Y\subseteq X$ with $X\setminus Y$ smooth. If
$\ecd(U)<2\dim(X)-r$ then $H^i_{\et}(X,\ZZ/\ell\ZZ)\to
H^i_{\et}(Y,\ZZ/\ell\ZZ)$ is an isomorphism for $i<r$ and
injective for $i=r$.

In the \cite{Varbaro-TAMS12}, Varbaro also investigates
the interaction of \'etale cohomological dimension with intersections:
let $\KK$ be an algebraic closed field of arbitrary characteristic and
let $X$ and $Y$ be two smooth projective varieties of dimension at
least 1. Set $Z = X \times Y \subseteq \PP^N_\KK$ (any embedding) and
$U = P^N \setminus Z$. Then $\ecd(U) \geq 2N - 3$. In particular, if
$\dim Z \geq 3$ then $Z$ cannot be a set-theoretic complete
intersection by (\ref{eqn-ara-and-ecd}).

\subsection{Other applications of local cohomology to geometry}

\subsubsection{Bockstein morphisms}
In this subsection we discuss a construction that originates (to our
knowledge) in topology but can, in principle, be used as a tool to
study any linear functor in prime characteristic. 

For this we need the following concept. A collection of functors
$\{\funcF^\bullet\}$ is a \emph{covariant $\delta$-functor} (in the
sense of Grothendieck) if for each short exact
sequence of $A$-modules $0\to M'\to M\to M''\to 0$ one obtains a
functorial long exact sequence
\[
\ldots\to \funcF^i(M')\to \funcF^i(M)\to \funcF^i(M'')\to
\funcF^{i+1}(M')\to\cdots
\]

Now suppose that for some $A$-module $M$, multiplication by 
$f\in A$ induces an injection $0\to
M\stackrel{f}{\to}M\stackrel{\pi}{\to} M/fM\to 0$. If 
$\funcF^\bullet$ is a covariant $\delta$-functor that is $A$-linear
(\emph{i.e.}, each $\funcF^i$ is additive, and
$\funcF^i(M\stackrel{a\cdot
  h}{\to}N)=\funcF^i(M)\stackrel{a\cdot\funcF^i(h)}{\to}\funcF^i(N)$ for all
$a\in A$ and all $h\in\Hom_A(M,N)$)
then there is an induced long exact
sequence
\[
\cdots
\to
\funcF^i(M/fM)\stackrel{\delta^{\funcF,i}_f}{\to}
\funcF^{i+1}(M)\stackrel{f\cdot}{\to}
\funcF^{i+1}(M)\stackrel{\pi^{\funcF,i+1}_f}{\to}
\funcF^{i+1}(M/fM)\stackrel{\delta^{\funcF,i+1}_f}{\to}\cdots
\]
Now one can define a sequence of \emph{Bockstein} morphisms
\[
\beta^{\funcF,i}_f\colon \funcF^i(M/fM)\to
\funcF^{i+1}(M/fM)
\]
as the composition
\[
\beta^{\funcF,i}_f=\pi^{\funcF,i+1}_f\circ \delta^{\funcF,i}_f.
\]
\begin{rmk}\label{rmk-Bockstein}
  \begin{asparaenum}
  \item Clearly, $f,i,\funcF$ and the $A/fA$-module $M/fM$ are
    ingredients of a Bockstein morphism. However, while the notation
    does not indicate this, is also depends on $A$ and the avatar
    $M\stackrel{f}{\to} M$ for $M/fM$ (or at least an infinitesimal
    avatar $0\to M/fM\to M/f^2M\to M/fM\to0$). Bocksteins are not intrinsic
    but arise from a specialization.
  \item It is
    possible to modify the constructions to include contravariant
    functors, or $A$-modules $N$ on which $f$ acts surjectively.
  \end{asparaenum}
\end{rmk}

The orginal version of a Bockstein morphism appeared in topology,
where $A=\ZZ$, $f$ is a prime number, $M$ is an Abelian  group without
$p$-torsion, and
$\funcF^\bullet$ is singular homology (or cohomology) with coefficients $M$ on a
fixed space $X$. Generally, in this context there is a Bockstein
spectral sequence that arises from the short exact sequence of
singular chains on $X$ with coefficients in $M,M$ and $M/pM$
respectively. It starts with $E^1_{i,j}=H_{i+j}(X;M/pM)$, the
differential on the $E^1$-page is the Bockstein morphism, and it converges to the tensor product of $\ZZ/p\ZZ$ with the free part of $H_{i+j}(X;M)$.

In \cite{SinghWalther-Bockstein}, Bockstein maps were introduced and
studied in local cohomology. So, $A$ is a Noetherian $\ZZ$-algebra,
$I=\ideal{\bsg=g_1,\ldots,g_m}\subseteq A$ is an ideal, and $M$ is a
$p$-torsion free $A$-module. In this setup there are several
$\delta$-functors $\funcF^\bullet$ that arise naturally: the local
cohomology functor $\funcF^i=H^i_I(-)$ with support in $I$, the
extension functors $\funcF^i=\Ext^i_A(A/I^\ell,-)$, the Koszul
cohomology functors $\funcF^i=H^i(-;\bsg)$. It is shown in
\cite{SinghWalther-Bockstein} that in the same way that these three
$\delta$-functors allow natural transformations, the three families of
Bockstein morphisms are compatible. Several examples are given, based
(for example) on the arithmetic of elliptic curves and on subspace
arrangements.

One result of \cite{SinghWalther-Bockstein} states that when $A$ is a
polynomial ring over $\ZZ$ containing the ideal $I$, then the
Bockstein on $H^\bullet_I(R/pR)$ is zero except for a finite set of
primes $p$. On a more topological note, the same article investigates
the interplay between Bocksteins on local cohomology and those on
singular homology in the context of Stanley--Reisner rings. More
precisely, let $R=\ZZ[\bsx]$ be the $\ZZ^n$-graded polynomial ring on
the vertices of the simplicial complex $\Delta$ on $n$ vertices, and let
$\frakm=\ideal{\bsx}$ be the graded maximal ideal.  
Hochster linked the multi-graded components of the local cohomology
$H^\bullet_\frakm(M\otimes_\ZZ R/I)$ with the singular cohomology with
coefficients in $M$ of a certain simplicial subcomplex of $\Delta$
determined by the chosen multi-degree,
\cite{Hochster-Norman-notes}. Then \cite{SinghWalther-Bockstein} shows
that the topological Bocksteins on these links are compatible with the
local cohomology Bocksteins via Hochster's identification, and that it
behaves well with respect to local duality.

It follows easily from the definitions that the composition of
Bocksteins $\beta^{\funcF,i+1}_f\circ \beta^{\funcF,i}_f$ is zero;
this is the origin of the Bockstein spectral sequence mentioned
above. Its ingredients are the \emph{Bockstein cohomology modules}
$\ker(\beta^{\funcF,i+1}_f)/\image(\beta^{\funcF,i}_f)$.  In
\cite{Puthenpurakal-Vietnam19}, this notion is used to study the
extended Rees ring $A[It,t^{-1}]$ of an $\frakm$-primary ideal in the
local ring $(A,\frakm)$ as $M$, using $t$ for $f$ and $\funcF$ is the
local cohomology with support in $\frakm$. The accomplishment consists
in vanishing theorems for local cohomology of the associated graded
ring $\gr_I(A)$, extending earlier such results of Narita, and
Huckaba--Huneke \cite{Narita-PCPS63,HuckabaHuneke-Crelle99}.

\subsubsection{Variation of Hodge structures and GKZ-systems}

Here we give a brief motivation of $A$-hypergeometric
systems and explain how local cohomology of toric varieties enters the
picture. We recommend \cite{SST00,Stienstra,RSSW} for more
detailed information and literature sources.

\begin{ntn}\label{ntn-A}
  Let $\{\bolda_1,\ldots,\bolda_n\}=A\subseteq \ZZ^{d\times n}$ satisfy
the following properties:
\begin{asparaenum}
\item \label{ntn-A-1} the cone $C_A:=\RR_{\geq 0}A$ spanned by the columns of $A$
  inside $\RR^d$ is $d$-dimensional and its lineality (the
  dimension of the largest real vector space that it contains) is
  zero;
\item \label{ntn-A-2} there exists a $\ZZ$-linear functional
  \[
  h\colon \ZZ^d\to\ZZ
  \]
  such that $h(\bolda_j)=1$ for $1\le j\le n$;
\item \label{ntn-A-3} the semigroup $\NN A:=\sum_{j=1}^n \NN \bolda_j$
  agrees with the intersection $\ZZ^d\cap C_A$.
\end{asparaenum}
\end{ntn}

The graded (via $h$) semigroup ring
\[
S_A:=\CC[\NN A]
\]
gives rise to a
projective toric variety $Y_A\subseteq \PP^{n-1}_\CC$ of dimension
$d-1$ and its cone $X_A=\Spec(S_A)\subseteq \CC^n$.
They can be viewed as (partial) compactifications of the $(d-1)$-torus
\[
\TT:=\underbrace{\Hom_\ZZ(\ZZ^d,\CC^*)}_{=:\tilde\TT}/\Hom_\ZZ(\ZZ \bolda_0,\CC^*),
\]
and $\tilde \TT$ respectively, 
where $\bolda_0$ is a suitable element of $\ZZ^d\cap C_A$ that induces
$h$ in the sense that $h(\bolda_j)$ is the dot product
$\langle\bolda_0,\bolda_j\rangle$. 

A global section $F_{A,\pointx}\in \Gamma(Y_A,\calO_{Y_A}(1))$ is an
element $\sum \pointx_jt^{\bolda_j}$ of the Laurent polynomial ring
$\CC[t_1^\pm,\ldots,t_d^\pm]$ that is equivariant under the action of
$\Hom_\ZZ(\ZZ \bolda_0,\CC^*)$. Its vanishing defines a hypersurface
$Z_\pointx$ inside $\TT$ with complement $U_\pointx=\TT\setminus
Z_\pointx$.  Batyrev initiated the study of the Hodge theory of these
objects in his search for mirror symmetry on toric varieties and their
hypersurfaces \cite{Batyrev-Duke93}. As is explained in Stienstra's
article \cite{Stienstra}, for understanding the weight filtration on
the cohomology of $Z_\pointx$ it is useful to study Hodge aspects of the
cohomology $H^\bullet(\tilde\TT,\tilde Z_\pointx;\CC)$ relative to the affine
cone
\[
\tilde Z_\pointx:=\tilde\TT\cap\Var(F_{A,\pointx}-1).
\]
A powerful tool in this endeavor is the idea of letting the section vary and
studying these cohomology groups as a family, viewing the coefficients
of the Laurent polynomial as parameters. For this, read the parameter
$\pointx$
of $F_{A,\pointx}$ as a
point in $\CC^n$. For any face
$\tau$ of the cone $C_A$ let $F^\tau_{A,\pointx}$ be the subsum of
$F_{A,\pointx}$ of terms with support on $\tau$. Then $\pointx$ is
\emph{non-degenerate} if the singular locus of $F^\tau_{A,\pointx}$
does not meet $\tilde \TT$ for any $\tau$, including the case
$\tau=A$.

For non-degenerate $\pointx$, 
$H^i(\tilde\TT,\tilde Z_{\pointx};\CC)$ is nonzero only when $i=d$ and there is a
natural identification of $H^d(\tilde\TT,\tilde Z_{\pointx};\CC)$ with
the stalk of the solutions of a certain natural $D$-module that we
describe next.

For what is to follow, we assume that $A$ satisfies condition
\ref{ntn-A}\eqref{ntn-A-1} but not necessarily
\ref{ntn-A}\eqref{ntn-A-2} and \ref{ntn-A}\eqref{ntn-A-3}, unless
indicated expressly.

Let $D_A$ be the Weyl algebra
$\CC[x_1,\ldots,x_n]\langle\del_1,\ldots,\del_n\rangle$ with subring $O_A=\CC[x_1,\ldots,x_n]$, and let $L_A$ be
the $\ZZ$-kernel of $A$. Define two types of operators
\begin{align*}
 && &E_i:=\sum_{j=1}^n a_{i,j}x_j\del_j&\text{(Euler operators);}\\
 && &\Box_\boldu:=\del^{\boldu_+}-\del^{\boldu_-}&\text{(box operators).}
\end{align*}
Here, $1\le i\le d$ and $\boldu\in L_A$ with
$(\boldu_+)_i=\max\{\boldu_i,0\}$ and 
$(\boldu_-)_i=\max\{-\boldu_i,0\}$.
Then choose a parameter vector
$\beta\in\CC^d$ and define the \emph{hypergeometric ideal}
\[
H_A(\beta)=D_A\cdot\{E_i-\beta_i\}_{i=1}^d+D_A\cdot\{\Box_\boldu\}_{\boldu\in
  L_A}
\]
and the \emph{hypergeometric module}
\[
M_A(\beta):=D_A/H_A(\beta)
\]
to $A$ and $\beta$. These modules were defined by Gelfand, Graev,
Kapranov and Zelevinsky in a string of articles including
\cite{GGZ87,GKZ89} during their investigations of Aomoto type
integrals. 
The modules are always holonomic \cite{Adolphson-Duke94}, and they are
regular holonomic if and only if $A$ satisfies Condition
\ref{ntn-A}.\eqref{ntn-A-2}, \cite{SchulzeWalther-slopes}.
We refer to \cite{SST00} for extensive background
on hypergeometric functions, their associated differential equations,
and how they relate to hypergeometric modules $M_A(\beta)$  via a
dehomogenization technique investigated in \cite{BMW-Adv19}. The
article \cite{RSSW} is a gentle introduction to hypergeometric
$D$-modules, combined with a survey on recent applications to Hodge theory.

Let $R_A=\CC[\del_1,\ldots,\del_n]$; while this is a subring of
operators of $D_A$, one can also view it as a polynomial ring in its
own right. The ideal
\[
I_A:=R_A\cdot\{\Box_\boldu\}_{\boldu\in L_A}
\]
that forms part of the defining equations for $H_A(\beta)$ is called
the \emph{toric ideal}; its variety in $\hat \CC^n=\Spec R_A$ is the
toric variety $X_A$. We use here the ``hat'' to
distinguish the copy of complex $n$-space that arises as $\Spec R_A$
from that which arises as $\Spec O_A$. The two are domain and target
of the Fourier--Laplace transform $\FL(-)$ which, on elements of $D_A$, amounts
to $x_j\mapsto \del_j, \del_j\mapsto -x_j$. 

Local cohomology arises in two ways in the study of $M_A(\beta)$: in
connection with the dimension of the space of solutions, and in
the limitations of a functorial description of $M_A(\beta)$ via a
$D$-module theoretic pushforward.

For any holonomic $D_A$-module $M$ there is a Zariski open set
of $\CC^n$ on which $M$ is a connection; we call the rank of
this connection the \emph{rank} of $M$. For $M_A(\beta)$, this open
set is determined by the non-vanishing of the
\emph{$A$-discriminant}, a generalization of the discriminant of a
polynomial. In particular, it does not depend on $\beta$; we denote it
$U_A$.
If $A$ satisfies Condition \ref{ntn-A}.\eqref{ntn-A-3}
then the connection on $U_A$ has rank equal to the simplicial volume
in $\RR^d$ of the convex hull of the origin and the columns of $A$,
\cite{GGZ87,GKZ89,Adolphson-Duke94}. Indeed, the hypothesis implies that the
semigroup ring $S_A$ is Cohen--Macaulay by Hochster's theorem
\cite{Hochster-Annals72}, and this allows a certain spectral sequence to
degenerate, which determines the rank. In fact, one can even produce
the solutions often in explicit forms, by writing down suitable
hypergeometric series and proving convergence \cite{GGZ87,SST00}.

In the absence of Condition
\ref{ntn-A}.\eqref{ntn-A-3}, the situation can be more interesting since
then there may be choices of $\beta$  with the effect of changing the
rank \cite{MMW05}. That the possibility of changing
rank exists at all was discovered in \cite{SturmfelsTakayama98}.
A certain Koszul-like complex based
on the Euler operators $E_i$ that appeared in \cite{MMW05} can be used
to substitute for the (now not degenerating) spectral sequence.

A
natural question is which parameters $\beta$ will show a change in
rank. Because of basic principles, the rank at special $\beta$ can
only go up \cite{MMW05}. Since $S_A$ is $A$-graded via
$\deg_A(\del_j)=\bolda_j\in\ZZ A$, so are its local cohomology modules
$H^i_\bsdel(S_A)$ supported at the homogeneous maximal ideal. Set
\[
\calE_A:=\overline{\bigcup_{i=0}^{d-1}\deg_A(H^i_\bsdel(S_A))}^{\text
  Zariski},
\]
the Zariski closure of the union of all $A$-degrees of nonzero
elements in a local cohomology module with $i<d$. Note that the union of
these degrees can be seen as witnesses to the failure of $S_A$ being
Cohen--Macaulay: the union is empty if and only if $S_A$ has full
depth. In generalization of the implication of equal rank for all
$\beta$ in the Cohen--Macaulay case, it is shown in \cite{MMW05} that  
\[
[\rk(M_A(\beta)>\vol(A)]\Leftrightarrow [\beta\in\calE_A].
\]

Consider now the monomial map
\begin{eqnarray*}
\varphi=\varphi_A\colon \tilde \TT&\to&\hat\CC^n,\\
\pointt&\mapsto&(\pointt^{\bolda_1},\ldots,\pointt^{\bolda_n})
\end{eqnarray*}
induced by $A$. The map is an isomorphism onto the image by Condition
\ref{ntn-A}.\eqref{ntn-A-1}, and its closure is the toric variety
$X_A$. On $\tilde \TT$ one has for each $\beta$ the (regular)
connection
$\calL_\beta=D_\TT/D_\TT\cdot\{t_i\del_{t_i}+\beta_i\}_1^d$. In
\cite{GKZ90}, Gelfand, Kapranov and Zelevinsky proved that if $\beta$
is sufficiently generic then the Fourier--Laplace transform
$\FL(M_A(\beta))$ agrees with the $\calD$-module direct image
$\varphi_+(\calL_\beta)$, where the set of ``good'' $\beta$ forms the
complement of a countably infinite and locally finite hyperplane
arrangement called the \emph{resonant} parameters, and given by all
$L_A$-shifts of the bounding hyperplanes of the cone $C_A$. In
\cite{SchulzeWalther-ekdi} this result was refined and completed to an
equivalence
\[
[M_A(\beta)=\varphi_+(\calL_\beta)]\Leftrightarrow [\beta\text{ is not
    strongly resonant}].
\]
Here, following \cite{SchulzeWalther-ekdi}, $\beta\in\CC^d$ is
\emph{strongly resonant} if and only if there is a finitely generated
$R_A$-submodule of $\bigoplus_{j=1}^n H^1_{\del_j}(S_A)$ containing
$\beta$ in the Zariski closure of its $A$-degrees.  (Since the local
cohomology modules here are not coherent, being strongly resonant is more
special than being in the Zariski closure of the $A$-degrees of the
direct sum). Some further improvements have been made in
\cite{Steiner-JPAA,Steiner-JA}.

\begin{rmk}
  As it turns out, when Conditions \ref{ntn-A} are in force in full
strength, then certain $M_A(\beta)$, including the case
$\beta=\boldzero$,  are not just a regular $D_A$-module but in
fact carry a mixed Hodge module structure in the sense of Saito,
\cite{SaitoMHM}. The Hodge and weight filtrations of hypergeometric
systems
have been studied in \cite{ReiSe-Hodge,Reich2,RW-weight}, showing
connections to intersection homology of toric varieties. See
\cite{RSSW} for a survey.
\end{rmk}

\subsubsection{Milnor fibers and torsion in the Jacobian ring}

Let $f$ be a non-unit in $R=\CC[x_1,\ldots,x_n]$  and put
$X:=\CC^n=\Spec(R)$. By the ideal $J_f$ we mean the ideal
generated by the partial derivatives $\frac{\del f}{\del
  x_1},\ldots,\frac{\del f}{\del x_n}$; this ideal varies with the
choice of coordinate system in which we calculate. In contrast, the
Jacobian ideal $\Jac(f)=J_f+(f)$ is independent.

If
$\pointx\in\Var(f)$, let $B(\pointx,\eps)$ denote the $\eps$-ball around
$\pointx\in\Var(f)\subseteq \CC^n$.  Milnor \cite{Milnor} proved that
the diffeomorphism type of the open real manifold
\[
M_{f,\pointx,t,\eps}=B(\pointx,\eps)\cap \Var(f-t)
\]
is independent of $\eps,t$ as long as $0<|t|\ll\eps\ll 1$. Abusing
language, for $0<t\ll\eps\ll 1$ denote by $M_{f,\pointx}$ the fiber
of the bundle
\[
B(\pointx,\eps)\cap \{\pointy\in\CC^n\mid 0<|f(\pointy)|<t\}\to
f(\pointy).
\]

%

If $f$ has an isolated singularity at $\pointx$ then the Milnor fiber
$M_{f,\pointx}$ is a bouquet of $(n-1)$-spheres, and
$H^{n-1}(M_{f,\pointx};\CC)$ can be identified non-canonically with
the Jacobian ring $R/\Jac(f)$ as vector space; in particular, the
Jacobian ring ``knows'' the number of spheres in the
bouquet. 

We call $f$ \emph{quasi-homogeneous} under the weight
$(w_1,\ldots,w_n)\in\QQ^n$ if $\sum_{i=1}^n w_i\frac{\del}{\del
  x_i}(f)=f$.  In this case, the Jacobian ring acquires a
$\QQ[s]$-module structure where $s$ acts via the Euler homogeneity,
compare \cite{Malgrange-isolee}. This is actually true for general
isolated singularities, not just in the presence of homogeneity, and
the eigenvalues of the action of $s$ turn out to be the non-trivial
roots of the local Bernstein--Sato polynomial of $f$ at $\pointx$. The
$s$-action comes then from the Gau\ss--Manin connection. Compare also
\cite{Saito-ASPM09,Saito-1609.04801,Malgrange-evan,
  Kashiwara-LNM1016,Walther-survey}. Compare \cite{Steenbrink-survey}
for details on the Hodge structure on the cohomology of the Milnor fiber.

For non-isolated singularities, most of this must break down, since
$R/\Jac(f)$ is not Artinian in that case. Suppose from now on that $f$
is homogeneous, and that $\pointx$ is the origin. Note
that now $\Jac(f)=J_f$; we abbreviate $M_{f,\pointx}$ to $M_f=\Var(f-1)$. The
\emph{Jacobian module}
\[
H^0_\frakm(R/J_f)=\{g+J_f\mid \exists k\in\NN, \forall i, x_i^kg\in J_f\}
\]
has been studied in \cite{Pellikaan,vanStratenWarmt} for various symmetry
properties and connections with geometry. Note that this finite length
module agrees
with the Jacobian ring in the case of an isolated singularity, it can
hence be considered a generalization of it in more general settings.

If
\[
\eta=\sum_ix_i\de x_1\wedge\cdots\wedge\widehat{\de
  x_i}\wedge\cdots\wedge \de x_n
\]
denotes the canonical $(n-1)$-form on $X$, then (via residues) every
class in $H^{n-1}(M_f;\CC)$ is of the form $g\eta$ for suitable $g\in
R$, and 
if $g\in R$ is the smallest degree homogeneous
polynomial such that $g\eta$ represents a chosen class in
$H^{n-1}(M_f;\CC)$ then $-\deg(g\eta)/\deg(f)$ is a root of the
Bernstein--Sato polynomial of $f$, \cite{Walther-Bernstein}.  Suppose
the singular locus of $f$ is (at most) $1$-dimensional. Then
by \cite{Saito-BS-arr,Walther-zieg,Saito-1703.05741}, with $1\le k\le
d$ and $\lambda=\exp(2\pi \sqrt{-1}k/d)$, the following holds:
\[
\dim_\CC [H^0_\frakm(R_n/\Jac(f))]_{d-n+k}\le \dim_\CC \gr^{{\rm
    Hodge}}_{n-2}(H^{n-1}(M_f;\CC)_\lambda),
\]
where the right hand side indicates the $\lambda$-eigenspace of the
associated graded object to the Hodge filtration on
$H^{n-1}(M_{f};\CC)$. Dimca and Sticlaru have used this inequality to
study nearly free divisors and pole order filtrations,
\cite{DimcaSticlaru-pole,DimcaSticlaru-W}. It would be interesting to
find more general inequalities of this type. The above estimate is
based on local cohomology of logarithmic forms introduced in
\cite{KSaito-logforms}; such
modules have been calculated in \cite{DSSWW} for generic hyperplane
arrangements. See \cite{Walther-survey}
for more connections to monodromy and zeta-functions.

\subsection{Lyubeznik numbers}
\label{subsec-lambda}

Let $(R,\frakm,\kk)$ be a commutative regular local Noetherian ring of
dimension $n$ that contains its residue field.
For any
ideal $I$ of $R$, Lyubeznik proved in \cite{L-Dmods} that the
$\kk$-dimension
\[
  \lambda_{i,j}(R,I):=\dim_\kk(\Ext^i_R(\kk,H^{n-j}_I(R)))
\]
is for each $i,j\in\NN$ only a function of $R/I$ and so does not depend
on the presentation of $R/I$ as a quotient of a regular local ring.

In his seminal paper, Lyubeznik also showed that $\lambda_{i,j}(R,I)$
agrees with the socle dimension in $H^i_\frakm(H^{n-j}_I(R))$, and
hence with the $i$-th Bass number of $H^{n-j}_I(R)$ with respect to
$\frakm$. In fact,
 $H^i_\frakm(H^{n-j}_I(R))$  is the direct sum of
$\lambda_{i,j}(R/I)$ many copies of $E_R(\kk)$, the injective hull of
$\kk$ when viewed as $R$-module.

It follows from the local cohomology interpretation that
$\lambda_{i,j}(R,I)=\lambda_{i,j}(\hat R,I\hat R)$ is invariant under
completion. By the Cohen structure theorems, every complete local
Noetherian ring containing its residue field is the quotient of a
complete regular Noetherian local ring containing its residue
field. One can thus define for every local Noetherian ring $A$ the
\emph{$(i,j)$-Lyubeznik number}
\[
  \lambda_{i,j}(A):=\lambda_{i,j}(R,I)
\]
via any surjection $R\onto R/I=\hat A$ from a complete regular ring $R$
onto the completion of $A$.

\begin{ntn}
Throughout this subsection, $(R,\frakm,\kk)$ is a regular local ring containing its
residue field, $\hat R$ its completion along $\frakm$,
and $I$ an ideal of $R$ such that $A=R/I$. Set $d:=\dim(A)$.
Field extensions $R\leadsto \KK\otimes_\kk R$ have no impact on the
Lyubeznik numbers, so that one can always assume $\kk$ to be
algebraically or separably closed if necessary. Moreover, since
$\Gamma_I(M)=\Gamma_{\sqrt{I}}(M)$, one may assume that $A$ is
reduced.
\end{ntn}

By Grothendieck's vanishing theorem, $\lambda_{i,j}(A)$ is zero if
$j<0$, and by the depth sensitivity of local cohomology,
$\lambda_{i,j}(A)=0$ if $j>\dim(A)$, \cite{24h}. By construction, the
dimension of the support of $H^{n-j}_I(R)$ is contained in the variety
of $I$, so that $\lambda_{i,j}(A)=0$ for all $i>d$.

We can thus write $\Lambda(A)$ for the \emph{Lyubeznik table}
\[
  \Lambda(A):=\left(
    \begin{array}{ccc}
      \lambda_{0,0}(A)&\dots&\lambda_{0,d}(A)\\
      \vdots&&\vdots\\
      \lambda_{d,0}(A)&\dots&\lambda_{d,d}(A)
    \end{array}
  \right)
\]
It has been shown in \cite{HunekeSharp} in the case $\charac(R)>0$, and then in
\cite{L-Dmods} when $\QQ\subseteq \kk$ that the injective dimension of
$H^k_I(R)$ is always bounded above by the dimension of its
support. However, it is standard that the support of $H^{n-j}_I(R)$ is
contained in a variety of dimension at most $j$.
This implies that the nonzero entries of $\Lambda(A)$ are on
or above the main diagonal of $\Lambda(A)$. 

\medskip

There is a Grothendieck spectral sequence
\begin{eqnarray}\label{eqn-Gss}
  H^i_\frakm(H^j_I(R))\Longrightarrow H^{i+j}_\frakm(R).
\end{eqnarray}
It follows directly from this spectral sequence that
\begin{itemize}
\item the alternating sum $\sum_{i,j}(-1)^{i+j}\lambda_{i,j}(A)$
  equals $1$;
\item $\lambda_{0,d}(A)=\lambda_{1,d}(A)=0$ for all $A$ unless
  $\dim(A)\le 1$;
\item if $R/I$ is a complete intersection, then $\lambda_{i,j}(A)$
  vanishes unless $i=j=d$. (We say that \emph{the Lyubeznik
    table is trivial}).
\item Moreover, following \cite{AlvarezYanagawa18} let
  $\rho_j(A):=-\delta_{0,j}+\sum_{i=0}^{d-j}\lambda_{i,i+j}(A)$ be the
  reduced sum along the $j$-th super-diagonal in $\Lambda(A)$, where
  $\delta$ denotes the Kronecker-$\delta$.  Then $\rho_d(A)$ is
  always zero, and) non-vanishing of $\rho_j(A)$ implies that of
  either $\rho_{j-1}(A)$ or $\rho_{j+1}(A)$, compare
  \cite{BetancourtSpiroffWitt}.
\end{itemize}
In characteristic $p>0$, the (iterated) Frobenius functor sends
a free resolution of the ideal $I$ to a free resolution of the
Frobenius power $I^{[p^e]}$. As the Frobenius powers of $I$ are
cofinal with the usual powers, $H^k_I(R)=0$ whenever $k$ exceeds the
projective dimension of $R/I$. In particular, if $I$ is perfect
(\emph{i.e.}, $R/I$ is Cohen--Macaulay), the Lyubeznik table of $R/I$
is trivial in positive characteristic. In characteristic zero, this is
not so; for example, the Lyubeznik table for the (perfect) ideal of
the $2\times 2$ minors of a $2\times 3$ matrix of indeterminates over
$\kk\supseteq \QQ$ is
\[
  \left(
    \begin{array}{ccccc}
      0&0&0&1&0\\
      \cdot&0&0&0&0\\
      \cdot&\cdot&0&0&1\\
      \cdot&\cdot&\cdot&0&0\\
      \cdot&\cdot&\cdot&\cdot&1
    \end{array}\right)
\]
as one sees from the fact that $I$ is $3$-generated and
$H^3_I(R)=E_R(\kk)$, compare Example \ref{exa-2x3}, or the
computations in \cite{W-lcD}. (Single dots indicate a zero entry.)

\begin{dfn}
The \emph{highest Lyubeznik number} of $A$ is $\lambda_{d,d}(A)$.\end{dfn}

It follows directly from the spectral sequence that for $d\le 1$, only
$\lambda_{d,d}(A)$ is nonzero (and thus equal to $1$).

Lyubeznik proved in \cite{L-Dmods} that $\lambda_{d,d}(A)$ is always
positive. For 2-dimensional complete local rings, with separably
closed residue field, it was shown in \cite{W-lambda,Kawasaki00} that
the Lyubeznik table is independent of the $1$-dimensional components of
$I$. Indeed, one has:
\[
  \Lambda(A)=\left(\begin{array}{ccc}
                     0&t-1&0\\
                     \cdot&0&0\\
                     \cdot&\cdot&t
                   \end{array}
  \right)
\]
where $t$ is the number of components of the punctured spectrum of
$A$.
In any
dimension $d$, the
number $\lambda_{d,d}(A)$ is $1$ if $A$ is analytically normal \cite{L-Dmods} or has
Serre's condition $S_2$ \cite{Kawasaki02}. On the other hand,
$\lambda_{d,d}(A)$ can be $1$ without $A$ 
being Cohen--Macaulay or even $S_2$, \cite{Kawasaki02}.

More generally, consider the \emph{Hochster--Huneke graph} of $A$: the
vertices of $\HH(A)$ are the $d$-dimensional primes of $A$ and an edge
links two such primes if the height of their sum is $1$. Then Zhang,
generalizing the case $d\le 2$ from \cite{W-lambda} and the case
$\charac(\kk)>0$ from \cite{L-lcInvs}, proved (in a
characteristic-independent way) in
\cite{Zhang-highest} that $\lambda_{d,d}(A)$ agrees
with the number of connected components of $\HH(A)$. The main result in \cite{Zhang-highest} has been extended to mixed characteristic in \cite{ZhangVanishingLCMixedChar}. See
also \cite{SchenzelMahmood,Schenzel-OnConnectedness} for more on the
relationship between connectedness and the structure of local cohomology.

\subsubsection{Combinatorial cases and topology}

If $I$ is a monomial ideal, then Alvarez, Vahidi and Yanagawa
\cite{Yanagawa03,AlvarezVahidi14,Alvarez15,AlvarezYanagawa18} have
obtained the following results:
\begin{itemize}
\item Lyubeznik numbers of monomial ideals relate to linear strands
  of the minimal free resolution of their Alexander duals;
\item If $A$ is sequentially Cohen--Macaulay (\emph{i.e.}, every
  $\Ext^i_R(A,R)$ is zero or Cohen--Macaulay of dimension $i$) then
  both in characteristic $p>0$ and also if $I$ is monomial then the
  Lyubeznik table is trivial. 
\item there are Thom--Sebastiani type results for Lyubeznik tables of
  monomial ideals in disjoint sets of variables.
\item Lyubeznik numbers of Stanley--Reisner rings are topological
  invariants attached to the underlying simplicial complex.
\end{itemize}

In a different direction, consider the case when $I_{r,s,t}$ is the ideal
generated by the $(t+1)\times (t+1)$ minors of an $r\times s$ matrix of indeterminates over
the field $\kk$. In positive characteristic, the Cohen--Macaulayness
of $R/I$ implies triviality of the Lyubeznik table. In characteristic
zero, however, these numbers carry interesting combinatorial
information related to representations of the general linear
group. L\"orincz and Raicu proved in \cite{LoerinczRaicu-lambda} the following.
Write the Lyubeznik numbers into a bivariate generating function 
\[
  L_{r,s,t}(q,w):=\sum_{i,j\geq 0}\lambda_{i,j}(A_{r,s,t})\cdot
  q^i\cdot w^j
\]
with $A_{r,s,t}=\CC[\{x_{i,j}|1\le i\le r,1\le j\le
s\}]/I_{r,s,t}$, with $r>s>t$. Then
\[
  L_{r,s,t}(q,w)=\sum_{i=0}^t q^{i^2+i(r-s)}\cdot {s\choose
    i}_{q^2}\cdot w^{t^2+2t+i(r+s-2t-2)}\cdot {s-1-i\choose t-i}_{w^2}.
\]
Here, the subscripts to the binomial coefficient indicate the Gaussian
$q$-binomial expression
${a\choose b}_c
=\frac{(1-c^a)(1-c^{a-1})\cdot\ldots\cdot(1-c^{a-b+1})}{(1-c^b)(1-c^{b-1})\cdot\ldots\cdot(1-c)}$.

There is a similar formula for the case $r=s>t$. 

\medskip

We now turn to topological interpretations on Lyubeznik tables.
The earliest such results  were formulated
by Garc\'ia L\'opez and Sabbah. Suppose $A$ has an isolated
singularity at $\frakm$. Then
$H^{n-j}_I(R)$ is $\frakm$-torsion for $n-j\neq d$. Hence, by the
spectral sequence \eqref{eqn-Gss}, $\Lambda(A)$ is concentrated in the
top row and the rightmost column, and there are equalities
$\lambda_{0,j}+\delta_{j+1,d}=\lambda_{j+1,d}$, using again
Kronecker notation. It is shown in \cite{GarciaSabbah} that, if the
coefficient field is $\CC$, then
$\lambda_{0,j}$ equals the $\CC$-dimension of the topological local
cohomology group of the analytic space $\Spec(V)$ with support in the
vertex $\frakm$. 

This result was then generalized by Blickle and Bondu as
follows. Suppose (over $\CC$) that the constant sheaf on the spectrum
of $A$ is self-dual in
the sense of Verdier outside the vertex $\frakm$. This is 
the case when $\frakm$ is an isolated singularity, but it also occurs
in more general cases. For example, on a hypersurface $f=0$ this condition
is equivalent to the Bernstein--Sato polynomial of $f$ having no other
integral root but $-1$, and $-1$ occurring with multiplicity one,
\cite{Torrelli-intHom}. Blickle and Bondu prove in \cite{BlickleBondu} that
in this situation the same interpretation of $\Lambda(A)$ can be made
as in the article by Garc\'ia L\'opez and Sabbah. In parallel, they
also show that if the field has finite characteristic, a corresponding
interpretation can be made in terms of local \'etale cohomology with
supports at the vertex.

\medskip

Lyubeznik numbers also contain information on connectedness of
algebraic varieties. For example, 
as mentioned before, for $\dim(A)=2$ over a separably closed field, the
Lyubeznik table is entirely characterized by the number of connected
components of the $2$-dimensional part of the punctured spectrum.

Suppose $A$ is equidimensional, with separably closed coefficient
field $\kk$. Denote by $\kappa(A)$ the \emph{connectedness dimension
  of $A$}, the smallest dimension $t$ of a subvariety $Y$ in
$\Spec(A)$ whose removal leads to a
disconnection. N\'u\~nez-Betancourt, Spiroff and Witt discuss in
\cite{BetancourtSpiroffWitt} the relationship between the number
$\kappa(A)$ and the vanishing of certain Lyubeznik numbers. Their
results generalize a consequence of the Second Vanishing Theorem
that can be phrased as: $H^{n-1}_I(R)=0 $ if and only if
$\kappa(A)\neq 0$. To be precise, they show for an equidimensional
ring $A$:
\begin{itemize}
\item $[\kappa(A)\geq 1]\Longleftrightarrow [\lambda_{0,1}(A)=0]$;
\item $[\kappa(A)\geq 2]\Longleftrightarrow
  [\lambda_{0,1}(A)=\lambda_{1,2}(A)=0]$;
\item for $i< \dim(A)$, $[\kappa(A)\geq i]\Longleftarrow [\lambda_{0,1}(A)=\cdots=\lambda_{i-1,i}(A)=0]$.
\end{itemize}
Earlier, Dao and Takagi, inspired by remarks of Varbaro, showed that
over any field, Serre's condition $S_3$ implies that
$\lambda_{d-1,d}=0$, \cite{DaoTakagi-Compositio16},
while in increasing generality it was shown in
\cite{W-lambda,L-lcInvs,Zhang-highest} that $[\kappa(A)\geq \dim(A)-1]\Leftrightarrow
[\lambda_{d,d}(A)=1]$. In \cite{RWZ-lambda} are some other results on the effect
of Serre's conditions $(S_t)$ on $\Lambda(A)$.

\subsubsection{Projective Lyubeznik numbers}

Suppose $X=X_\kk$ is a projective variety of dimension $d-1$, with
embedding $\iota\colon X\into\PP:=\PP^{n-1}_\kk$ via sections of the
line bundle $\calL=\iota^*(\calO_\PP(1))$.  With this embedding comes
a global coordinate ring $\Gamma_*(\PP)$ of $\PP$ and a homogeneous
ideal defining the cone $C(X)$ over $X$ in the corresponding affine
space. Let $R$ be the localization of $\Gamma_*(\PP)$ at the vertex,
and let $I$ be the ideal defining the germ of $C(X)$ in $R$. A natural
question is to ask:
\begin{prb}
  To what extent are the Lyubeznik numbers of $R/I$ dependent on the
  embedding $\iota$?
\end{prb}
Certainly, if two such cones $(R,I)$ and $(R',I')$ arise from one
another by an automorphism of $\PP$, then the attached Lyubeznik
tables are equal. It is less clear from the definitions
whether two embeddings that produce the same sheaf $\calL$ on $X$, or at
least the same element in the Picard group,
should give the same Lyubeznik tables. And even more difficult is the
question whether $\iota,\iota'$ should give rise to equal Lyubeznik
tables when $\calL_\iota\neq\calL_{\iota'}$ in the Picard group.

We say that \emph{$\Lambda(X)$ (\emph{or just $\lambda_{i,j}(X)$}) is
  projective} if each cone derived from a projective embedding of $X$
produces the same $\Lambda$-table (or at least the same
$\lambda_{i,j}$). 
Positive known results include the following:
\begin{itemize}
\item If $\dim(X)\le 1$ then $\Lambda(X)$ is projective by
  \cite{W-lambda}, since then each cone ring is at most
  $2$-dimensional, and connectedness of the punctured $d$-dimensional
  spectrum of $R/I$ is equivalent to connectedness of the
  $(d-1)$-dimensional part of $X$.
\item If $X$ is smooth and $\kk=\CC$, then each cone has an isolated
  singularity, so that the Lyubeznik numbers can be expressed in terms
  of topological local cohomology as in \cite{GarciaSabbah}. Switala proves in 
  \cite{Switala} that these data are actually intrinsic to $X$,
  appearing as cokernels of the cup product with the Chern  class of
  the embedding on singular cohomology of $X$. By independence
  of Lyubeznik numbers under field extensions, this also works when
  just $\QQ\subseteq \kk$.
\item Since $\lambda_{0,1}(A)=0$ is equivalent to $H^{n-1}_I(R)=0$,
  which in turn is equivalent to connectedness of the punctured
  spectrum of $A$, $\lambda_{0,1}(X)$ is projective.
\item Similarly, the simultaneous vanishing of
  $\lambda_{0,1},\lambda_{1,2},\ldots,\lambda_{i-1,i}$
    is projective since it measures by \cite{BetancourtSpiroffWitt} the
  connectedness dimension of the cone, which corresponds to
  connectedness dimension of $X$ itself.
\end{itemize}
Consider the module
$\calE_{i,j}(\iota):=\Ext^{n-i}_R(\Ext^{n-j}_R(R/I,\Omega_R),\Omega_R)$
where $\Omega_R$ is the canonical module of $R$.  In
\cite{Zhang-lambda}, Zhang proves that in finite characteristic, the
degree zero part of $\calE_{i,j}(\iota)$ supports a natural action of
Frobenius, whose stable part is independent of $\iota$ and has
$\kk$-dimension $\lambda_{i,j}(R/I)$. In particular, $\Lambda$ is
projective in positive characteristic.

In characteristic zero, after base change to $\CC$, the modules
$H^i_\frakm(H^{n-j}_I(R))$ have a natural structure as mixed Hodge
modules. This has been exploited in \cite{RSW-lambda} to prove that in
this setting, on the level of constructible sheaves via the
Riemann--Hilbert correspondence,
\[
  \lambda_{i,j}(R/I)=\dim_\CC H^i \tau^!{}^p\calH^{-j}(\DD\QQ_C).
\]
Here, $\QQ_C$ is the constant sheaf on the cone $C=C(X)$ under any
embedding of $X$, $\DD$ is Verdier duality (corresponding to holonomic
duality), ${}^p\calH$ is taking perverse cohomology (corresponding to
usual cohomology for $D$-modules via the Riemann--Hilbert
correspondence) and $\tau^!$ is the exceptional inverse image for
constructible sheaves under
the embedding $\tau$ of the vertex into the cone. One can then recast
this as the dimension of the cohomology of a certain related sheaf on
the punctured cone, and this cohomology is the middle term in an exact
sequence whose other terms are kernels and cokerels of the Chern class
of $\calL_\iota$ on certain sheaves on $X$. These sheaves are
relatives of, but not always equal to, intersection cohomology of
$X$. This difference is then exploited to construct examples of
(reducible) varieties whose Lyubeznik numbers are not projective. In
\cite{Wang-lambda}, the construction was modified to yield irreducible
ones with non-projective $\Lambda$-table.

The construction of \cite{RSW-lambda} starts with a variety whose
Picard number is greater than one, and from it constructs a suitable
$X$.  In \cite{RWZ-lambda} it is shown that if the rational Picard
group of $X$ is $\QQ$ then almost all Lyubeznik numbers of $X$ are
projective. In particular, this applies to determinantal varieties so
that the L\"orincz--Raicu computation in \cite{LoerinczRaicu-lambda}
determines the vast majority of the entries of the Lyubeznik tables for
such varieties under all embeddings. 

\begin{rmk}
  A similar set (to Lyubeznik numbers) of invariants is introduced in
  \cite{Bridgland} (but see also \cite{Switala-dR-complete}). It is
  shown that if $I$ is an ideal in a polynomial ring over the
  complex numbers then the \v Cech-to-de Rham spectral sequence whose
  abutment is the reduced singular cohomology of the complement of the variety
  of $I$ has terms on page two that do not depend on the embedding of
  the variety of $I$ into an affine space, at least when suitably
  re-indexed. Using algebraic de Rham cohomology, this is actually
  shown over all fields of characteristic zero. These \emph{\v
    Cech--de Rham numbers} are further investigated in
  \cite{RWZ-lambda} from the viewpoint of projectivity since, if $I$ is
  homogeneous, one can ask to what extent these numbers are defined by
  the associated projective variety (rather than the affine cone). \cite{RWZ-lambda} studies their behavior under Veronese maps, and
  the degeneration of the spectral sequence.
\end{rmk}

\newcommand{\etalchar}[1]{$^{#1}$}
\def\cprime{$'$}
\providecommand{\bysame}{\leavevmode\hbox to3em{\hrulefill}\thinspace}
\providecommand{\MR}{\relax\ifhmode\unskip\space\fi MR }
\providecommand{\MRhref}[2]{%
  \href{http://www.ams.org/mathscinet-getitem?mr=#1}{#2}
}
\providecommand{\href}[2]{#2}


\begin{thebibliography}{HNnBPW19}

\bibitem[ABW13]{DonuParsaJarek-MathAnn13}
Donu Arapura, Parsa Bakhtary, and Jaros{\l}aw W{\l}odarczyk, \emph{Weights on
  cohomology, invariants of singularities, and dual complexes}, Math. Ann.
  \textbf{357} (2013), no.~2, 513--550. \MR{3096516}

\bibitem[Ado94]{Adolphson-Duke94}
Alan Adolphson, \emph{Hypergeometric functions and rings generated by
  monomials}, Duke Math. J. \textbf{73} (1994), no.~2, 269--290.

\bibitem[AKM98]{AberbachKatzmanMcCrimmon-JA98}
Ian Aberbach, Mordechai Katzman, and Brian MacCrimmon, \emph{Weak
  {F}-regularity deforms in {${\bf Q}$}-{G}orenstein rings}, J. Algebra
  \textbf{204} (1998), no.~1, 281--285. \MR{1623973}

\bibitem[AM69]{A+M}
M.~F. Atiyah and I.~G. Macdonald, \emph{Introduction to commutative algebra},
  Addison-Wesley Publishing Co., Reading, Mass.-London-Don Mills, Ont., 1969.
  \MR{0242802}

\bibitem[AM00]{MontanerJPPA2000}
Josep \`Alvarez~Montaner, \emph{Characteristic cycles of local cohomology
  modules of monomial ideals}, J. Pure Appl. Algebra \textbf{150} (2000),
  no.~1, 1--25. \MR{1762917}

\bibitem[AM04]{MontanerJPPA2004}
\bysame, \emph{Characteristic cycles of local cohomology modules of monomial
  ideals. {II}}, J. Pure Appl. Algebra \textbf{192} (2004), no.~1-3, 1--20.
  \MR{2067186}

\bibitem[AM05]{Alvarez-JSC05}
\bysame, \emph{Operations with regular holonomic {$\mathscr D$}-modules with
  support a normal crossing}, J. Symbolic Comput. \textbf{40} (2005), no.~2,
  999--1012. \MR{2167680}

\bibitem[AM15]{Alvarez15}
\bysame, \emph{Lyubeznik table of sequentially {C}ohen-{M}acaulay rings}, Comm.
  Algebra \textbf{43} (2015), no.~9, 3695--3704. \MR{3360843}

\bibitem[AMBL05]{AlvarezBlickleLyubeznik}
Josep Alvarez-Montaner, Manuel Blickle, and Gennady Lyubeznik, \emph{Generators
  of {$D$}-modules in positive characteristic}, Math. Res. Lett. \textbf{12}
  (2005), no.~4, 459--473. \MR{2155224}

\bibitem[AMBZ20]{MontanerBoixZarzuelaIMRN2020}
Josep \`Alvarez~Montaner, Alberto~F. Boix, and Santiago Zarzuela, \emph{On some
  local cohomology spectral sequences}, Int. Math. Res. Not. IMRN (2020),
  no.~19, 6197--6293. \MR{4165477}

\bibitem[AMGLZA03]{AlvarezGarciaZarzuela}
Josep \`Alvarez~Montaner, Ricardo Garc\'{\i}a~L\'{o}pez, and Santiago
  Zarzuela~Armengou, \emph{Local cohomology, arrangements of subspaces and
  monomial ideals}, Adv. Math. \textbf{174} (2003), no.~1, 35--56. \MR{1959890}

\bibitem[AMHNnB17]{AlvarezMontanerHunekeBetancourt}
Josep \`Alvarez~Montaner, Craig Huneke, and Luis N\'{u}\~{n}ez Betancourt,
  \emph{{$D$}-modules, {B}ernstein-{S}ato polynomials and {$F$}-invariants of
  direct summands}, Adv. Math. \textbf{321} (2017), 298--325. \MR{3715713}

\bibitem[AML06]{LeykinMontaner2006}
Josep \`Alvarez~Montaner and Anton Leykin, \emph{Computing the support of local
  cohomology modules}, J. Symbolic Comput. \textbf{41} (2006), no.~12,
  1328--1344. \MR{2271328}

\bibitem[AMV14]{AlvarezVahidi14}
Josep \`Alvarez~Montaner and Alireza Vahidi, \emph{Lyubeznik numbers of
  monomial ideals}, Trans. Amer. Math. Soc. \textbf{366} (2014), no.~4,
  1829--1855. \MR{3152714}

\bibitem[AMY18]{AlvarezYanagawa18}
Josep \`Alvarez~Montaner and Kohji Yanagawa, \emph{Lyubeznik numbers of local
  rings and linear strands of graded ideals}, Nagoya Math. J. \textbf{231}
  (2018), 23--54. \MR{3845587}

\bibitem[ATJLL99]{TarrioJeremiasLipman}
Leovigildo Alonso~Tarr\'{\i}o, Ana Jerem\'{\i}as~L\'{o}pez, and Joseph Lipman,
  \emph{Studies in duality on {N}oetherian formal schemes and non-{N}oetherian
  ordinary schemes}, Contemporary Mathematics, vol. 244, American Mathematical
  Society, Providence, RI, 1999. \MR{1716706}

\bibitem[Bah17]{Bahmanpour-CiA2017}
Kamal Bahmanpour, \emph{A note on {L}ynch's conjecture}, Comm. Algebra
  \textbf{45} (2017), no.~6, 2738--2745. \MR{3594553}

\bibitem[Bar70]{Barth-AJM70}
W.~Barth, \emph{Transplanting cohomology classes in complex-projective space},
  Amer. J. Math. \textbf{92} (1970), 951--967. \MR{287032}

\bibitem[Bar95]{Barile-JA95}
Margherita Barile, \emph{Arithmetical ranks of ideals associated to symmetric
  and alternating matrices}, J. Algebra \textbf{176} (1995), no.~1, 59--82.
  \MR{1345294}

\bibitem[Bar04]{Barile-MR2104196}
\bysame, \emph{On the computation of arithmetical ranks}, Int. J. Pure Appl.
  Math. \textbf{17} (2004), no.~2, 143--161. \MR{2104196}

\bibitem[Bar06a]{Barile-MR2222158}
\bysame, \emph{On toric varieties of high arithmetical rank}, Yokohama Math. J.
  \textbf{52} (2006), no.~2, 125--130. \MR{2222158}

\bibitem[Bar06b]{Barile-MR2244384}
\bysame, \emph{On toric varieties which are almost set-theoretic complete
  intersections}, J. Pure Appl. Algebra \textbf{207} (2006), no.~1, 109--118.
  \MR{2244384}

\bibitem[Bar07]{Barile-MR2441609}
\bysame, \emph{On simplicial toric varieties of codimension 2}, Rend. Istit.
  Mat. Univ. Trieste \textbf{39} (2007), 9--42. \MR{2441609}

\bibitem[Bat93]{Batyrev-Duke93}
Victor~V. Batyrev, \emph{Variations of the mixed {H}odge structure of affine
  hypersurfaces in algebraic tori}, Duke Math. J. \textbf{69} (1993), no.~2,
  349--409.

\bibitem[BB05]{BlickleBondu}
Manuel Blickle and Raphael Bondu, \emph{Local cohomology multiplicities in
  terms of \'{e}tale cohomology}, Ann. Inst. Fourier (Grenoble) \textbf{55}
  (2005), no.~7, 2239--2256. \MR{2207383}

\bibitem[BB11]{BlickleBoeckle-Crelle11}
Manuel Blickle and Gebhard B\"{o}ckle, \emph{Cartier modules: finiteness
  results}, J. Reine Angew. Math. \textbf{661} (2011), 85--123. \MR{2863904}

\bibitem[BBL{\etalchar{+}}14]{BBLSZ-2014}
Bhargav Bhatt, Manuel Blickle, Gennady Lyubeznik, Anurag~K. Singh, and Wenliang
  Zhang, \emph{Local cohomology modules of a smooth {$\mathbb{Z}$}-algebra have
  finitely many associated primes}, Invent. Math. \textbf{197} (2014), no.~3,
  509--519. \MR{3251828}

\bibitem[BBL{\etalchar{+}}19]{BBLSZ2}
\bysame, \emph{Stabilization of the cohomology of thickenings}, Amer. J. Math.
  \textbf{141} (2019), no.~2, 531--561. \MR{3928045}

\bibitem[BBL{\etalchar{+}}21]{BBLSZ3}
\bysame, \emph{An asymptotic vanishing theorem for the cohomology of
  thickenings}, Math. Ann. \textbf{380} (2021), no.~1-2, 161--173. \MR{4263681}

\bibitem[BE18]{BoixEghbali-Annihilators}
Alberto~F. Boix and Majid Eghbali, \emph{Annihilators of local cohomology
  modules and simplicity of rings of differential operators}, Beitr. Algebra
  Geom. \textbf{59} (2018), no.~4, 665--684. \MR{3871100}

\bibitem[Ber72]{Bernstein-bfu}
I.~N. Bern\v{s}te\u{\i}n, \emph{Analytic continuation of generalized functions
  with respect to a parameter}, Funkcional. Anal. i Prilo\v{z}en. \textbf{6}
  (1972), no.~4, 26--40. \MR{0320735}

\bibitem[BGfGf72]{BGG-Ell}
I.~N. Bern\v{s}te\u{\i}n, I.~M. Gel\cprime~fand, and S.~I. Gel\cprime~fand,
  \emph{Differential operators on a cubic cone}, Uspehi Mat. Nauk \textbf{27}
  (1972), no.~1(163), 185--190. \MR{0385159}

\bibitem[BH93]{BrunsHerzogCMBook}
Winfried Bruns and J\"{u}rgen Herzog, \emph{Cohen-{M}acaulay rings}, Cambridge
  Studies in Advanced Mathematics, vol.~39, Cambridge University Press,
  Cambridge, 1993. \MR{1251956}

\bibitem[Bha12]{Bhatt-ANT2012}
Bhargav Bhatt, \emph{Annihilating the cohomology of group schemes}, Algebra
  Number Theory \textbf{6} (2012), no.~7, 1561--1577. \MR{3007159}

\bibitem[Bha20]{BhattRPlusMixedChar}
\bysame, \emph{Cohen-{M}acaulayness of absolute integral closures},
  arXiv:2008.08070, 2020.

\bibitem[Bit20]{Bitoun-IMRN20}
Thomas Bitoun, \emph{Length of local cohomology in positive characteristic and
  ordinarity}, Int. Math. Res. Not. IMRN (2020), no.~7, 1921--1932.
  \MR{4089437}

\bibitem[Bj{\"{o}}79]{BjorkBook}
J.-E. Bj{\"{o}}rk, \emph{Rings of differential operators}, North-Holland
  Mathematical Library, vol.~21, North-Holland Publishing Co., Amsterdam-New
  York, 1979. \MR{549189}

\bibitem[BL05]{BarileLyubeznik-PAMS05}
Margherita Barile and Gennady Lyubeznik, \emph{Set-theoretic complete
  intersections in characteristic {$p$}}, Proc. Amer. Math. Soc. \textbf{133}
  (2005), no.~11, 3199--3209. \MR{2160181}

\bibitem[BL10]{BerkeschLeykin}
Christine Berkesch and Anton Leykin, \emph{Algorithms for {B}ernstein-{S}ato
  polynomials and multiplier ideals}, I{SSAC} 2010---{P}roceedings of the 2010
  {I}nternational {S}ymposium on {S}ymbolic and {A}lgebraic {C}omputation, ACM,
  New York, 2010, pp.~99--106. \MR{2920542}

\bibitem[Bli04]{BlickleIntersectionHomology2004}
Manuel Blickle, \emph{The intersection homology {$D$}-module in finite
  characteristic}, Math. Ann. \textbf{328} (2004), no.~3, 425--450.
  \MR{2036330}

\bibitem[BM98]{BarileMorales-MR1621700}
Margherita Barile and Marcel Morales, \emph{On the equations defining
  projective monomial curves}, Comm. Algebra \textbf{26} (1998), no.~6,
  1907--1912. \MR{1621700}

\bibitem[BM19]{BarileMacchia-MR3957102}
Margherita Barile and Antonio Macchia, \emph{On determinantal ideals and
  algebraic dependence}, Comm. Algebra \textbf{47} (2019), no.~6, 2357--2366.
  \MR{3957102}

\bibitem[BMT00]{BarileMoralesThoma-MR1752767}
Margherita Barile, Marcel Morales, and Apostolos Thoma, \emph{On simplicial
  toric varieties which are set-theoretic complete intersections}, J. Algebra
  \textbf{226} (2000), no.~2, 880--892. \MR{1752767}

\bibitem[BMT02]{BarileMoralesThoma-MR1896020}
\bysame, \emph{Set-theoretic complete intersections on binomials}, Proc. Amer.
  Math. Soc. \textbf{130} (2002), no.~7, 1893--1903. \MR{1896020}

\bibitem[BMW19]{BMW-Adv19}
Christine Berkesch, Laura~Felicia Matusevich, and Uli Walther, \emph{Torus
  equivariant {$D$}-modules and hypergeometric systems}, Adv. Math.
  \textbf{350} (2019), 1226--1266. \MR{3949610}

\bibitem[Bre05]{BrennerLinearBound-MMJ2005}
Holger Brenner, \emph{A linear bound for {F}robenius powers and an inclusion
  bound for tight closure}, Michigan Math. J. \textbf{53} (2005), no.~3,
  585--596. \MR{2207210}

\bibitem[Bri20]{Bridgland}
Nicole Bridgland, \emph{On the de {R}ham homology of affine varieties in
  characteristic 0}, arXiv:2006.01334.

\bibitem[Bru89]{Bruns-Berkeley87}
Winfried Bruns, \emph{Additions to the theory of algebras with straightening
  law}, Commutative algebra ({B}erkeley, {CA}, 1987), Math. Sci. Res. Inst.
  Publ., vol.~15, Springer, New York, 1989, pp.~111--138. \MR{1015515}

\bibitem[BS90]{BrunsSchwanzl-BLMS90}
Winfried Bruns and Roland Schw\"{a}nzl, \emph{The number of equations defining
  a determinantal variety}, Bull. London Math. Soc. \textbf{22} (1990), no.~5,
  439--445. \MR{1082012}

\bibitem[BS98]{BrodmannSharp}
M.~P. Brodmann and R.~Y. Sharp, \emph{Local cohomology: an algebraic
  introduction with geometric applications}, Cambridge Studies in Advanced
  Mathematics, vol.~60, Cambridge University Press, Cambridge, 1998.
  \MR{1613627}

\bibitem[BSR81]{BresinskyStueckradRenschuch}
Henrik Bresinsky, J\"{u}rgen St\"{u}ckrad, and Bodo Renschuch,
  \emph{Mengentheoretisch vollst\"{a}ndige {D}urchschnitte verschiedener
  rationaler {R}aumkurven im {${\bf P}^{3}$} \"{u}ber {K}\"{o}rpern von
  {P}rimzahlcharakteristik}, Math. Nachr. \textbf{104} (1981), 147--169.
  \MR{657889}

\bibitem[BST15]{BlickleSchwedeTucker-AJM2015}
Manuel Blickle, Karl Schwede, and Kevin Tucker, \emph{{$F$}-singularities via
  alterations}, Amer. J. Math. \textbf{137} (2015), no.~1, 61--109.
  \MR{3318087}

\bibitem[BST17]{BhattSchwedeTakagi}
Bhargav Bhatt, Karl Schwede, and Shunsuke Takagi, \emph{The weak ordinarity
  conjecture and {$F$}-singularities}, Higher dimensional algebraic
  geometry---in honour of {P}rofessor {Y}ujiro {K}awamata's sixtieth birthday,
  Adv. Stud. Pure Math., vol.~74, Math. Soc. Japan, Tokyo, 2017, pp.~11--39.
  \MR{3791207}

\bibitem[BSTZ10]{BlickleSchwedeTakagiZhang-MathAnn2010}
Manuel Blickle, Karl Schwede, Shunsuke Takagi, and Wenliang Zhang,
  \emph{Discreteness and rationality of {$F$}-jumping numbers on singular
  varieties}, Math. Ann. \textbf{347} (2010), no.~4, 917--949. \MR{2658149}

\bibitem[Cha04]{ChardinRegularity-2004}
Marc Chardin, \emph{Regularity of ideals and their powers}, 2004,
  Pr{\'e}publication 364. Institut de math{\'e}matiques de Jussieu.

\bibitem[CN78]{CowsikNori}
R.~C. Cowsik and M.~V. Nori, \emph{Affine curves in characteristic {$p$} are
  set theoretic complete intersections}, Invent. Math. \textbf{45} (1978),
  no.~2, 111--114. \MR{472835}

\bibitem[Cou95]{CoutinhoBook}
S.~C. Coutinho, \emph{A primer of algebraic {$D$}-modules}, London Mathematical
  Society Student Texts, vol.~33, Cambridge University Press, Cambridge, 1995.
  \MR{1356713}

\bibitem[DB81]{DuBois-BSMF81}
Philippe Du~Bois, \emph{Complexe de de {R}ham filtr\'{e} d'une vari\'{e}t\'{e}
  singuli\`ere}, Bull. Soc. Math. France \textbf{109} (1981), no.~1, 41--81.
  \MR{613848}

\bibitem[Del71]{Deligne-HodgeII}
Pierre Deligne, \emph{Th\'{e}orie de {H}odge. {II}}, Inst. Hautes \'{E}tudes
  Sci. Publ. Math. (1971), no.~40, 5--57. \MR{498551}

\bibitem[Del74]{Deligne-HodgeIII}
\bysame, \emph{Th\'{e}orie de {H}odge. {III}}, Inst. Hautes \'{E}tudes Sci.
  Publ. Math. (1974), no.~44, 5--77. \MR{498552}

\bibitem[DS19a]{DimcaSticlaru-pole}
Alexandru Dimca and Gabriel Sticlaru, \emph{Computing the monodromy and pole
  order filtration on {M}ilnor fiber cohomology of plane curves}, J. Symbolic
  Comput. \textbf{91} (2019), 98--115. \MR{3860886}

\bibitem[DS19b]{DimcaSticlaru-W}
\bysame, \emph{Line and rational curve arrangements, and {W}alther's
  inequality}, Atti Accad. Naz. Lincei Rend. Lincei Mat. Appl. \textbf{30}
  (2019), no.~3, 615--633. \MR{4002214}

\bibitem[DSM20]{MaStefani-2009.09038}
Alessandro De~Stefani and Linquan Ma, \emph{{F}-stable secondary
  representations and deformation of {F}-injectivity}, Preprint
  arXiv:2009.09038, 2020.

\bibitem[DSS{\etalchar{+}}13]{DSSWW}
G.~Denham, H.~Schenck, M.~Schulze, M.~Wakefield, and U.~Walther, \emph{Local
  cohomology of logarithmic forms}, Ann. Inst. Fourier (Grenoble) \textbf{63}
  (2013), no.~3, 1177--1203. \MR{3137483}

\bibitem[DSZ19]{DattaSwitalaZhang}
Rankeya Datta, Nicholas Switala, and Wenliang Zhang, \emph{Annihilators of
  {$D$}-modules in mixed characteristic}, arXiv:1907.09948, Math. Res. Lett., to appear.

\bibitem[DT16]{DaoTakagi-Compositio16}
Hailong Dao and Shunsuke Takagi, \emph{On the relationship between depth and
  cohomological dimension}, Compos. Math. \textbf{152} (2016), no.~4, 876--888.
  \MR{3484116}

\bibitem[EE73]{EisenbudEvans-Inv73}
David Eisenbud and E.~Graham Evans, Jr., \emph{Every algebraic set in
  {$n$}-space is the intersection of {$n$} hypersurfaces}, Invent. Math.
  \textbf{19} (1973), 107--112. \MR{327783}

\bibitem[EH08]{EnescuHochster-ANT08}
Florian Enescu and Melvin Hochster, \emph{The {F}robenius structure of local
  cohomology}, Algebra Number Theory \textbf{2} (2008), no.~7, 721--754.
  \MR{2460693}

\bibitem[Eis95]{Eisenbud}
David Eisenbud, \emph{Commutative algebra}, Graduate Texts in Mathematics, vol.
  150, Springer-Verlag, New York, 1995, With a view toward algebraic geometry.
  \MR{1322960}

\bibitem[EK04]{Emerton-Kisin}
Matthew Emerton and Mark Kisin, \emph{The {R}iemann-{H}ilbert correspondence
  for unit {$F$}-crystals}, Ast\'{e}risque (2004), no.~293, vi+257.
  \MR{2071510}

\bibitem[Elk87]{Elkies}
Noam~D. Elkies, \emph{The existence of infinitely many supersingular primes for
  every elliptic curve over {${\bf Q}$}}, Invent. Math. \textbf{89} (1987),
  no.~3, 561--567. \MR{903384}

\bibitem[EMS00]{EisenbudMustataStillman}
David Eisenbud, Mircea Musta\c{t}\v{a}, and Mike Stillman, \emph{Cohomology on
  toric varieties and local cohomology with monomial supports}, J. Symbolic
  Comput. \textbf{29} (2000), no.~4-5, 583--600, Symbolic computation in
  algebra, analysis, and geometry (Berkeley, CA, 1998). \MR{1769656}

\bibitem[Ene03]{Enescu-PAMS03}
Florian Enescu, \emph{F-injective rings and {F}-stable primes}, Proc. Amer.
  Math. Soc. \textbf{131} (2003), no.~11, 3379--3386. \MR{1990626}

\bibitem[Ene09]{Enescu-JA09}
\bysame, \emph{Local cohomology and {F}-stability}, J. Algebra \textbf{322}
  (2009), no.~9, 3063--3077. \MR{2567410}

\bibitem[Fal80]{FaltingsCrelle1980}
Gerd Faltings, \emph{\"{U}ber lokale {K}ohomologiegruppen hoher {O}rdnung}, J.
  Reine Angew. Math. \textbf{313} (1980), 43--51. \MR{552461}

\bibitem[Fed83]{Fedder-TAMS83}
Richard Fedder, \emph{{$F$}-purity and rational singularity}, Trans. Amer.
  Math. Soc. \textbf{278} (1983), no.~2, 461--480. \MR{701505}

\bibitem[Fer79]{Ferrand-MR0555692}
Daniel Ferrand, \emph{Set-theoretical complete intersections in characteristic
  {$p>0$}}, Algebraic geometry ({P}roc. {S}ummer {M}eeting, {U}niv.
  {C}openhagen, {C}openhagen, 1978), Lecture Notes in Math., vol. 732,
  Springer, Berlin, 1979, pp.~82--89. \MR{555692}

\bibitem[FW89]{FedderWatanabe-CA89}
Richard Fedder and Keiichi Watanabe, \emph{A characterization of
  {$F$}-regularity in terms of {$F$}-purity}, Commutative algebra ({B}erkeley,
  {CA}, 1987), Math. Sci. Res. Inst. Publ., vol.~15, Springer, New York, 1989,
  pp.~227--245. \MR{1015520}

\bibitem[Gab81]{GabberIntegrability1981}
Ofer Gabber, \emph{The integrability of the characteristic variety}, Amer. J.
  Math. \textbf{103} (1981), no.~3, 445--468. \MR{618321}

\bibitem[Gab04]{Gabber-04}
\bysame, \emph{Notes on some {$t$}-structures}, Geometric aspects of {D}work
  theory. {V}ol. {I}, {II}, Walter de Gruyter, Berlin, 2004, pp.~711--734.
  \MR{2099084}

\bibitem[Gal85]{Galligo-Linz85}
Andr\'{e} Galligo, \emph{Some algorithmic questions on ideals of differential
  operators}, E{UROCAL} '85, {V}ol. 2 ({L}inz, 1985), Lecture Notes in Comput.
  Sci., vol. 204, Springer, Berlin, 1985, pp.~413--421. \MR{826576}

\bibitem[GGM85]{GalligoGrangerMaisonobe}
A.~Galligo, M.~Granger, and Ph. Maisonobe, \emph{{${\mathscr D}$}-modules et
  faisceaux pervers dont le support singulier est un croisement normal}, Ann.
  Inst. Fourier (Grenoble) \textbf{35} (1985), no.~1, 1--48. \MR{781776}

\bibitem[GGZ87]{GGZ87}
I.~M. Gel{\cprime}fand, M.~I. Graev, and A.~V. Zelevinsky, \emph{Holonomic
  systems of equations and series of hypergeometric type}, Dokl. Akad. Nauk
  SSSR \textbf{295} (1987), no.~1, 14--19. \MR{902936}

\bibitem[GKZ90]{GKZ90}
Israel~M. Gel{\cprime}fand, Mikhail~M. Kapranov, and Andrei~V. Zelevinsky,
  \emph{Generalized {E}uler integrals and {$A$}-hypergeometric functions}, Adv.
  Math. \textbf{84} (1990), no.~2, 255--271.

\bibitem[GLS98]{GarciaSabbah}
R.~Garc\'{\i}a~L\'{o}pez and C.~Sabbah, \emph{Topological computation of local
  cohomology multiplicities}, Collect. Math. \textbf{49} (1998), no.~2-3,
  317--324, Dedicated to the memory of Fernando Serrano. \MR{1677136}

\bibitem[GM88]{GoreskyMacPherson-SMT}
Mark Goresky and Robert MacPherson, \emph{Stratified {M}orse theory},
  Ergebnisse der Mathematik und ihrer Grenzgebiete (3) [Results in Mathematics
  and Related Areas (3)], vol.~14, Springer-Verlag, Berlin, 1988. \MR{932724}

\bibitem[GM92]{GreenleesMay}
J.~P.~C. Greenlees and J.~P. May, \emph{Derived functors of {$I$}-adic
  completion and local homology}, J. Algebra \textbf{149} (1992), no.~2,
  438--453. \MR{1172439}

\bibitem[GQS70]{GuilleminQuillenSternberg1970}
Victor~W. Guillemin, Daniel Quillen, and Shlomo Sternberg, \emph{The
  integrability of characteristics}, Comm. Pure Appl. Math. \textbf{23} (1970),
  no.~1, 39--77. \MR{461597}

\bibitem[Gro66]{Grothendieck-deRham}
A.~Grothendieck, \emph{On the de {R}ham cohomology of algebraic varieties},
  Inst. Hautes \'{E}tudes Sci. Publ. Math. (1966), no.~29, 95--103. \MR{199194}

\bibitem[Gro67]{EGA4-4}
\bysame, \emph{\'{E}l\'{e}ments de g\'{e}om\'{e}trie alg\'{e}brique. {IV}.
  \'{E}tude locale des sch\'{e}mas et des morphismes de sch\'{e}mas {IV}},
  Inst. Hautes \'{E}tudes Sci. Publ. Math. (1967), no.~32, 361. \MR{238860}

\bibitem[Gro68]{SGA2}
Alexander Grothendieck, \emph{Cohomologie locale des faisceaux coh\'{e}rents et
  th\'{e}or\`emes de {L}efschetz locaux et globaux {$(SGA$} {$2)$}},
  North-Holland Publishing Co., Amsterdam; Masson \& Cie, \'{E}diteur, Paris,
  1968, Augment\'{e} d'un expos\'{e} par Mich\`ele Raynaud, S\'{e}minaire de
  G\'{e}om\'{e}trie Alg\'{e}brique du Bois-Marie, 1962, Advanced Studies in
  Pure Mathematics, Vol. 2. \MR{0476737}

\bibitem[GS]{Macaulay2}
Daniel~R. Grayson and Michael~E. Stillman, \emph{Macaulay2, a software system
  for research in algebraic geometry}.

\bibitem[GZK89]{GKZ89}
I.~M. Gel{\cprime}fand, A.~V. Zelevinsky, and M.~M. Kapranov,
  \emph{Hypergeometric functions and toric varieties}, Funktsional. Anal. i
  Prilozhen. \textbf{23} (1989), no.~2, 12--26. \MR{1011353}

\bibitem[Har66]{Hartshorne-RD}
Robin Hartshorne, \emph{Residues and duality}, Lecture notes of a seminar on
  the work of A. Grothendieck, given at Harvard 1963/64. With an appendix by P.
  Deligne. Lecture Notes in Mathematics, No. 20, Springer-Verlag, Berlin-New
  York, 1966. \MR{0222093}

\bibitem[Har67]{Hartshorne-lc-notes}
\bysame, \emph{Local cohomology}, A seminar given by A. Grothendieck, Harvard
  University, Fall, vol. 1961, Springer-Verlag, Berlin-New York, 1967.
  \MR{0224620}

\bibitem[Har68]{Hartshorne-CDAV}
\bysame, \emph{Cohomological dimension of algebraic varieties}, Ann. of Math.
  (2) \textbf{88} (1968), 403--450. \MR{232780}

\bibitem[Har74]{Hartshorne-Bulletin74}
\bysame, \emph{Varieties of small codimension in projective space}, Bull. Amer.
  Math. Soc. \textbf{80} (1974), 1017--1032. \MR{384816}

\bibitem[Har75]{Hartshorne-DRCAV}
\bysame, \emph{On the {D}e {R}ham cohomology of algebraic varieties}, Inst.
  Hautes \'{E}tudes Sci. Publ. Math. (1975), no.~45, 5--99. \MR{432647}

\bibitem[Har77]{Hartshorne-book}
\bysame, \emph{Algebraic geometry}, Springer-Verlag, New York-Heidelberg, 1977,
  Graduate Texts in Mathematics, No. 52. \MR{0463157}

\bibitem[Har79]{Hartshorne-AJM79}
\bysame, \emph{Complete intersections in characteristic {$p>0$}}, Amer. J.
  Math. \textbf{101} (1979), no.~2, 380--383. \MR{527998}

\bibitem[Har98]{Hara-AJM98}
Nobuo Hara, \emph{A characterization of rational singularities in terms of
  injectivity of {F}robenius maps}, Amer. J. Math. \textbf{120} (1998), no.~5,
  981--996. \MR{1646049}

\bibitem[Har70]{Hartshorne-AffineDuality}
Robin Hartshorne, \emph{Affine duality and cofiniteness}, Invent. Math.
  \textbf{9} (1969/70), 145--164. \MR{257096}

\bibitem[Hel05]{Hellus-CiA05}
M.~Hellus, \emph{On the associated primes of {M}atlis duals of top local
  cohomology modules}, Comm. Algebra \textbf{33} (2005), no.~11, 3997--4009.
  \MR{2183976}

\bibitem[Hel07a]{Hellus-CiA07dual}
\bysame, \emph{Finiteness properties of duals of local cohomology modules},
  Comm. Algebra \textbf{35} (2007), no.~11, 3590--3602. \MR{2362672}

\bibitem[Hel07b]{Hellus-CiA07Matlis}
Michael Hellus, \emph{Matlis duals of top local cohomology modules and the
  arithmetic rank of an ideal}, Comm. Algebra \textbf{35} (2007), no.~4,
  1421--1432. \MR{2313677}

\bibitem[HH90]{HochsterHuneke-JAMS1990}
Melvin Hochster and Craig Huneke, \emph{Tight closure, invariant theory, and
  the {B}rian\c{c}on-{S}koda theorem}, J. Amer. Math. Soc. \textbf{3} (1990),
  no.~1, 31--116. \MR{1017784}

\bibitem[HH92]{HochsterHuneke-AnnMath1992}
\bysame, \emph{Infinite integral extensions and big {C}ohen-{M}acaulay
  algebras}, Ann. of Math. (2) \textbf{135} (1992), no.~1, 53--89. \MR{1147957}

\bibitem[HH94]{HH-TAMS94}
\bysame, \emph{{$F$}-regularity, test elements, and smooth base change}, Trans.
  Amer. Math. Soc. \textbf{346} (1994), no.~1, 1--62. \MR{1273534}

\bibitem[HH99]{HuckabaHuneke-Crelle99}
Sam Huckaba and Craig Huneke, \emph{Normal ideals in regular rings}, J. Reine
  Angew. Math. \textbf{510} (1999), 63--82. \MR{1696091}

\bibitem[HJ20]{HochsterJeffries-MRL2019}
Melvin Hochster and Jack Jeffries, \emph{Faithfulness of top local cohomology
  modules in domains}, Math. Res. Lett. \textbf{27} (2020), no.~6, 1755--1765.
  \MR{4216603}

\bibitem[HK91]{HunekeKoh-cofiniteness}
Craig Huneke and Jee Koh, \emph{Cofiniteness and vanishing of local cohomology
  modules}, Math. Proc. Cambridge Philos. Soc. \textbf{110} (1991), no.~3,
  421--429. \MR{1120477}

\bibitem[HKM09]{HunekeKatzMarley}
Craig Huneke, Daniel Katz, and Thomas Marley, \emph{On the support of local
  cohomology}, J. Algebra \textbf{322} (2009), no.~9, 3194--3211. \MR{2567416}

\bibitem[HL90]{HunekeLyubeznik}
C.~Huneke and G.~Lyubeznik, \emph{On the vanishing of local cohomology
  modules}, Invent. Math. \textbf{102} (1990), no.~1, 73--93. \MR{1069240}

\bibitem[HL07]{HunekeLyubeznik-AdvMath2007}
Craig Huneke and Gennady Lyubeznik, \emph{Absolute integral closure in positive
  characteristic}, Adv. Math. \textbf{210} (2007), no.~2, 498--504.
  \MR{2303230}

\bibitem[HMS14]{HoriuchiMillerShimomoto-IU14}
Jun Horiuchi, Lance~Edward Miller, and Kazuma Shimomoto, \emph{Deformation of
  {$F$}-injectivity and local cohomology}, Indiana Univ. Math. J. \textbf{63}
  (2014), no.~4, 1139--1157, With an appendix by Karl Schwede and Anurag K.
  Singh. \MR{3263925}

\bibitem[HNnB17]{HochsterNunezSupportLC}
Melvin Hochster and Luis N\'{u}\~{n}ez Betancourt, \emph{Support of local
  cohomology modules over hypersurfaces and rings with {FFRT}}, Math. Res.
  Lett. \textbf{24} (2017), no.~2, 401--420. \MR{3685277}

\bibitem[HNnBPW19]{BetancourtHernandezPerezWitt-TAMS19}
Daniel~J. Hern\'{a}ndez, Luis N\'{u}\~{n}ez Betancourt, Felipe P\'{e}rez, and
  Emily~E. Witt, \emph{Lyubeznik numbers and injective dimension in mixed
  characteristic}, Trans. Amer. Math. Soc. \textbf{371} (2019), no.~11,
  7533--7557. \MR{3955527}

\bibitem[Hoc72]{Hochster-Annals72}
M.~Hochster, \emph{Rings of invariants of tori, {C}ohen-{M}acaulay rings
  generated by monomials, and polytopes}, Ann. of Math. (2) \textbf{96} (1972),
  318--337. \MR{304376}

\bibitem[Hoc77]{Hochster-Norman-notes}
Melvin Hochster, \emph{Cohen-{M}acaulay rings, combinatorics, and simplicial
  complexes}, Ring theory, {II} ({P}roc. {S}econd {C}onf., {U}niv. {O}klahoma,
  {N}orman, {O}kla., 1975), 1977, pp.~171--223. Lecture Notes in Pure and Appl.
  Math., Vol. 26. \MR{0441987}

\bibitem[Hoc19]{Hochster-SurveyLC-CiA2020}
\bysame, \emph{Finiteness properties and numerical behavior of local
  cohomology}, Comm. Algebra \textbf{47} (2019), no.~6, 1--11. \MR{3941632}

\bibitem[Hol04]{HolmDiffOperHyperplaneArr}
P\"{a}r Holm, \emph{Differential operators on hyperplane arrangements}, Comm.
  Algebra \textbf{32} (2004), no.~6, 2177--2201. \MR{2099582}

\bibitem[HP21a]{HartshornePolini-20}
Robin Hartshorne and Claudia Polini, \emph{Quasi-cyclic modules and coregular
  sequences}, Math. Z. \textbf{299} (2021), no. 1-2, 123--138. \MR{4311598}

\bibitem[HP21b]{HartshornePolini-simple}
\bysame, \emph{Simple {$\mathcal{D}$}-module components of local cohomology
  modules}, J. Algebra \textbf{571} (2021), 232--257. \MR{4200718}

\bibitem[HS77]{HartshorneSpeiser-Annals77}
Robin Hartshorne and Robert Speiser, \emph{Local cohomological dimension in
  characteristic {$p$}}, Ann. of Math. (2) \textbf{105} (1977), no.~1, 45--79.
  \MR{441962}

\bibitem[HS93]{HunekeSharp}
Craig~L. Huneke and Rodney~Y. Sharp, \emph{Bass numbers of local cohomology
  modules}, Trans. Amer. Math. Soc. \textbf{339} (1993), no.~2, 765--779.
  \MR{1124167}

\bibitem[HS97]{HunekeSmith-KodairaVanishing-1997}
Craig Huneke and Karen~E. Smith, \emph{Tight closure and the {K}odaira
  vanishing theorem}, J. Reine Angew. Math. \textbf{484} (1997), 127--152.
  \MR{1437301}

\bibitem[HS08a]{HellusStueckrad-PAMS08endo}
M.~Hellus and J.~St\"{u}ckrad, \emph{On endomorphism rings of local cohomology
  modules}, Proc. Amer. Math. Soc. \textbf{136} (2008), no.~7, 2333--2341.
  \MR{2390499}

\bibitem[HS08b]{HellusSchenzel-JA08}
Michael Hellus and Peter Schenzel, \emph{On cohomologically complete
  intersections}, J. Algebra \textbf{320} (2008), no.~10, 3733--3748.
  \MR{2457720}

\bibitem[HS08c]{HellusStueckrad-PAMS08duals}
Michael Hellus and J\"{u}rgen St\"{u}ckrad, \emph{Matlis duals of top local
  cohomology modules}, Proc. Amer. Math. Soc. \textbf{136} (2008), no.~2,
  489--498. \MR{2358488}

\bibitem[Hsi12]{Jenchieh-TAMS12}
Jen-Chieh Hsiao, \emph{{$D$}-module structure of local cohomology modules of
  toric algebras}, Trans. Amer. Math. Soc. \textbf{364} (2012), no.~5,
  2461--2478. \MR{2888215}

\bibitem[Hsi15]{Hsiao-big}
\bysame, \emph{A remark on bigness of the tangent bundle of a smooth projective
  variety and {$D$}-simplicity of its section rings}, J. Algebra Appl.
  \textbf{14} (2015), no.~7, 1550098, 10. \MR{3339397}

\bibitem[HTT08]{HTT}
Ryoshi Hotta, Kiyoshi Takeuchi, and Toshiyuki Tanisaki, \emph{{$D$}-modules,
  perverse sheaves, and representation theory}, Progress in Mathematics, vol.
  236, Birkh\"{a}user Boston, Inc., Boston, MA, 2008, Translated from the 1995
  Japanese edition by Takeuchi. \MR{2357361}

\bibitem[Hun92a]{HunekeProblemsLC}
Craig Huneke, \emph{Problems on local cohomology}, Free resolutions in
  commutative algebra and algebraic geometry ({S}undance, {UT}, 1990), Res.
  Notes Math., vol.~2, Jones and Bartlett, Boston, MA, 1992, pp.~93--108.
  \MR{1165320}

\bibitem[Hun92b]{Huneke-Invent1992}
\bysame, \emph{Uniform bounds in {N}oetherian rings}, Invent. Math.
  \textbf{107} (1992), no.~1, 203--223. \MR{1135470}

\bibitem[Hun00]{HunekeSaturation-CiA2000}
\bysame, \emph{The saturation of {F}robenius powers of ideals}, Comm. Algebra
  \textbf{28} (2000), no.~12, 5563--5572, Special issue in honor of Robin
  Hartshorne. \MR{1808589}

\bibitem[Hun07]{Huneke-lc-notes}
\bysame, \emph{Lectures on local cohomology}, Interactions between homotopy
  theory and algebra, Contemp. Math., vol. 436, Amer. Math. Soc., Providence,
  RI, 2007, Appendix 1 by Amelia Taylor, pp.~51--99. \MR{2355770}

\bibitem[HV16]{HenriquesVarbaro16}
In\^{e}s Bonacho Dos~Anjos Henriques and Matteo Varbaro, \emph{Test, multiplier
  and invariant ideals}, Adv. Math. \textbf{287} (2016), 704--732. \MR{3422690}

\bibitem[HZ18]{HochsterZhang-TAMS2018}
Melvin Hochster and Wenliang Zhang, \emph{Content of local cohomology,
  parameter ideals, and robust algebras}, Trans. Amer. Math. Soc. \textbf{370}
  (2018), no.~11, 7789--7814. \MR{3852449}

\bibitem[ILL{\etalchar{+}}07]{24h}
Srikanth~B. Iyengar, Graham~J. Leuschke, Anton Leykin, Claudia Miller, Ezra
  Miller, Anurag~K. Singh, and Uli Walther, \emph{Twenty-four hours of local
  cohomology}, Graduate Studies in Mathematics, vol.~87, American Mathematical
  Society, Providence, RI, 2007. \MR{2355715}

\bibitem[Ive86]{Iversen}
Birger Iversen, \emph{Cohomology of sheaves}, Universitext, Springer-Verlag,
  Berlin, 1986. \MR{842190}

\bibitem[Jew94]{Jewell-Top94}
Ken Jewell, \emph{Complements of sphere and subspace arrangements}, Topology
  Appl. \textbf{56} (1994), no.~3, 199--214. \MR{1269311}

\bibitem[Jon94]{Jones-toric}
A.~G. Jones, \emph{Rings of differential operators on toric varieties}, Proc.
  Edinburgh Math. Soc. (2) \textbf{37} (1994), no.~1, 143--160. \MR{1258039}

\bibitem[Kas83]{Kashiwara-LNM1016}
M.~Kashiwara, \emph{Vanishing cycle sheaves and holonomic systems of
  differential equations}, Algebraic geometry ({T}okyo/{K}yoto, 1982), Lecture
  Notes in Math., vol. 1016, Springer, Berlin, 1983, pp.~134--142. \MR{726425
  (85e:58137)}

\bibitem[Kas95]{Kashi-thesis}
Masaki Kashiwara, \emph{Algebraic study of systems of partial differential
  equations}, M\'{e}m. Soc. Math. France (N.S.) (1995), no.~63, xiv+72.
  \MR{1384226}

\bibitem[Kas03]{Kashi-book}
\bysame, \emph{{$D$}-modules and microlocal calculus}, Translations of
  Mathematical Monographs, vol. 217, American Mathematical Society, Providence,
  RI, 2003, Translated from the 2000 Japanese original by Mutsumi Saito,
  Iwanami Series in Modern Mathematics. \MR{1943036}

\bibitem[Kas77]{Kashiwara-bfu}
\bysame, \emph{{$B$}-functions and holonomic systems. {R}ationality of roots of
  {$B$}-functions}, Invent. Math. \textbf{38} (1976/77), no.~1, 33--53.
  \MR{430304}

\bibitem[Kat98]{KatzmanComplexity-JA1998}
Mordechai Katzman, \emph{The complexity of {F}robenius powers of ideals}, J.
  Algebra \textbf{203} (1998), no.~1, 211--225. \MR{1620654}

\bibitem[Kat02]{Katzman-infinite-ass-primes}
\bysame, \emph{An example of an infinite set of associated primes of a local
  cohomology module}, J. Algebra \textbf{252} (2002), no.~1, 161--166.
  \MR{1922391}

\bibitem[Kat06]{KatzmanSupportTopLC}
\bysame, \emph{The support of top graded local cohomology modules}, Commutative
  algebra, Lect. Notes Pure Appl. Math., vol. 244, Chapman \& Hall/CRC, Boca
  Raton, FL, 2006, pp.~165--174. \MR{2184796}

\bibitem[Kaw00]{Kawasaki00}
Ken-ichiroh Kawasaki, \emph{On the {L}yubeznik number of local cohomology
  modules}, Bull. Nara Univ. Ed. Natur. Sci. \textbf{49} (2000), no.~2, 5--7.
  \MR{1814657}

\bibitem[Kaw02]{Kawasaki02}
\bysame, \emph{On the highest {L}yubeznik number}, Math. Proc. Cambridge
  Philos. Soc. \textbf{132} (2002), no.~3, 409--417. \MR{1891679}

\bibitem[Kha07]{Khashyarmanesh-Basel07}
Kazem Khashyarmanesh, \emph{On the {M}atlis duals of local cohomology modules},
  Arch. Math. (Basel) \textbf{88} (2007), no.~5, 413--418. \MR{2316886}

\bibitem[Kha10]{Khashyarmanesh-Canadian10}
\bysame, \emph{On the endomorphism rings of local cohomology modules}, Canad.
  Math. Bull. \textbf{53} (2010), no.~4, 667--673. \MR{2761689}

\bibitem[KLZ09]{KatzmanLyubeznikZhang-JA2009}
Mordechai Katzman, Gennady Lyubeznik, and Wenliang Zhang, \emph{On the
  discreteness and rationality of {$F$}-jumping coefficients}, J. Algebra
  \textbf{322} (2009), no.~9, 3238--3247. \MR{2567418}

\bibitem[KLZ16]{KLZ-PAMS2016}
\bysame, \emph{An extension of a theorem of {H}artshorne}, Proc. Amer. Math.
  Soc. \textbf{144} (2016), no.~3, 955--962. \MR{3447649}

\bibitem[KM17]{KimuraMontero-JCA17}
Kyouko Kimura and Paolo Mantero, \emph{Arithmetical rank of strings and
  cycles}, J. Commut. Algebra \textbf{9} (2017), no.~1, 89--106. \MR{3631828}

\bibitem[KMSZ18]{KatzmanMaSmirnovZhang2018}
Mordechai Katzman, Linquan Ma, Ilya Smirnov, and Wenliang Zhang,
  \emph{{$D$}-module and {$F$}-module length of local cohomology modules},
  Trans. Amer. Math. Soc. \textbf{370} (2018), no.~12, 8551--8580. \MR{3864387}

\bibitem[Kol97]{Kollar-Proc97}
J{\'a}nos Koll{\'a}r, \emph{Singularities of pairs}, Algebraic
  geometry---{S}anta {C}ruz 1995, Proc. Sympos. Pure Math., vol.~62, Amer.
  Math. Soc., Providence, RI, 1997, pp.~221--287. \MR{1492525 (99m:14033)}

\bibitem[KS90]{KashiwaraShapira}
Masaki Kashiwara and Pierre Schapira, \emph{Sheaves on manifolds}, Grundlehren
  der Mathematischen Wissenschaften [Fundamental Principles of Mathematical
  Sciences], vol. 292, Springer-Verlag, Berlin, 1990, With a chapter in French
  by Christian Houzel. \MR{1074006}

\bibitem[KSSZ14]{KatzmanSchwedeSinghZhang}
Mordechai Katzman, Karl Schwede, Anurag~K. Singh, and Wenliang Zhang,
  \emph{Rings of {F}robenius operators}, Math. Proc. Cambridge Philos. Soc.
  \textbf{157} (2014), no.~1, 151--167. \MR{3211813}

\bibitem[Kun69]{KunzRegularRing1969}
Ernst Kunz, \emph{Characterizations of regular local rings of characteristic
  {$p$}}, Amer. J. Math. \textbf{91} (1969), 772--784. \MR{252389}

\bibitem[Kun85]{Kunz-Intro}
\bysame, \emph{Introduction to commutative algebra and algebraic geometry},
  Birkh\"{a}user Boston, Inc., Boston, MA, 1985, Translated from the German by
  Michael Ackerman, With a preface by David Mumford. \MR{789602}

\bibitem[KZ18]{KatzmanZhangSupportLC}
Mordechai Katzman and Wenliang Zhang, \emph{The support of local cohomology
  modules}, Int. Math. Res. Not. IMRN (2018), no.~23, 7137--7155. \MR{3920344}

\bibitem[Lip02]{Lipman-Guanajuato}
Joseph Lipman, \emph{Lectures on local cohomology and duality}, Local
  cohomology and its applications ({G}uanajuato, 1999), Lecture Notes in Pure
  and Appl. Math., vol. 226, Dekker, New York, 2002, pp.~39--89. \MR{1888195}

\bibitem[LR20]{LoerinczRaicu-lambda}
Andr\'{a}s~C. L\H{o}rincz and Claudiu Raicu, \emph{Iterated local cohomology
  groups and {L}yubeznik numbers for determinantal rings}, Algebra Number
  Theory \textbf{14} (2020), no.~9, 2533--2569. \MR{4172715}

\bibitem[LRW19]{LoerinczRaicuWeyman-CiA19}
Andr\'{a}s~C. L\H{o}rincz, Claudiu Raicu, and Jerzy Weyman, \emph{Equivariant
  {$\mathcal D$}-modules on binary cubic forms}, Comm. Algebra \textbf{47}
  (2019), no.~6, 2457--2487. \MR{3957110}

\bibitem[LS01]{LyubeznikSmith-TAMS01}
Gennady Lyubeznik and Karen~E. Smith, \emph{On the commutation of the test
  ideal with localization and completion}, Trans. Amer. Math. Soc. \textbf{353}
  (2001), no.~8, 3149--3180. \MR{1828602}

\bibitem[LSW16]{LSW}
Gennady Lyubeznik, Anurag~K. Singh, and Uli Walther, \emph{Local cohomology
  modules supported at determinantal ideals}, J. Eur. Math. Soc. (JEMS)
  \textbf{18} (2016), no.~11, 2545--2578. \MR{3562351}

\bibitem[LW19]{LorinczWalther-equi}
Andr\'{a}s~C. L\H{o}rincz and Uli Walther, \emph{On categories of equivariant
  {$\mathcal{D}$}-modules}, Adv. Math. \textbf{351} (2019), 429--478.
  \MR{3952575}

\bibitem[LY18]{LYildirim-PAMS18}
Gennady Lyubeznik and Tu\u{g}ba Yildirim, \emph{On the {M}atlis duals of local
  cohomology modules}, Proc. Amer. Math. Soc. \textbf{146} (2018), no.~9,
  3715--3720. \MR{3825827}

\bibitem[Lyn12]{Lynch-Annihilators}
Laura~R. Lynch, \emph{Annihilators of top local cohomology}, Comm. Algebra
  \textbf{40} (2012), no.~2, 542--551. \MR{2889480}

\bibitem[Lyu85]{Lyubeznik-PAMS1985}
Gennady Lyubeznik, \emph{Some algebraic sets of high local cohomological
  dimension in projective space}, Proc. Amer. Math. Soc. \textbf{95} (1985),
  no.~1, 9--10. \MR{796437}

\bibitem[Lyu89]{L-ara-survey}
\bysame, \emph{A survey of problems and results on the number of defining
  equations}, Commutative algebra ({B}erkeley, {CA}, 1987), Math. Sci. Res.
  Inst. Publ., vol.~15, Springer, New York, 1989, pp.~375--390. \MR{1015529}

\bibitem[Lyu92]{Lyubeznik-AJM92}
\bysame, \emph{The number of defining equations of affine algebraic sets},
  Amer. J. Math. \textbf{114} (1992), no.~2, 413--463. \MR{1156572}

\bibitem[Lyu93a]{L-ecd}
\bysame, \emph{\'{E}tale cohomological dimension and the topology of algebraic
  varieties}, Ann. of Math. (2) \textbf{137} (1993), no.~1, 71--128.
  \MR{1200077}

\bibitem[Lyu93b]{L-Dmods}
\bysame, \emph{Finiteness properties of local cohomology modules (an
  application of {$D$}-modules to commutative algebra)}, Invent. Math.
  \textbf{113} (1993), no.~1, 41--55. \MR{1223223}

\bibitem[Lyu97]{LyubeznikFModules}
\bysame, \emph{{$F$}-modules: applications to local cohomology and
  {$D$}-modules in characteristic {$p>0$}}, J. Reine Angew. Math. \textbf{491}
  (1997), 65--130. \MR{1476089}

\bibitem[Lyu00a]{LyubeznikFinitenessCharFree2000}
\bysame, \emph{Finiteness properties of local cohomology modules: a
  characteristic-free approach}, J. Pure Appl. Algebra \textbf{151} (2000),
  no.~1, 43--50. \MR{1770642}

\bibitem[Lyu00b]{LyubeznikUnramified}
\bysame, \emph{Finiteness properties of local cohomology modules for regular
  local rings of mixed characteristic: the unramified case}, Comm. Algebra
  \textbf{28} (2000), no.~12, 5867--5882, Special issue in honor of Robin
  Hartshorne. \MR{1808608}

\bibitem[Lyu00c]{LyubeznikInjDim2000}
\bysame, \emph{Injective dimension of {$D$}-modules: a characteristic-free
  approach}, J. Pure Appl. Algebra \textbf{149} (2000), no.~2, 205--212.
  \MR{1757731}

\bibitem[Lyu02]{L-lc-survey}
\bysame, \emph{A partial survey of local cohomology}, Local cohomology and its
  applications ({G}uanajuato, 1999), Lecture Notes in Pure and Appl. Math.,
  vol. 226, Dekker, New York, 2002, pp.~121--154. \MR{1888197}

\bibitem[Lyu06a]{L-lcInvs}
\bysame, \emph{On some local cohomology invariants of local rings}, Math. Z.
  \textbf{254} (2006), no.~3, 627--640. \MR{2244370}

\bibitem[Lyu06b]{Lyubeznik-Compositio06}
\bysame, \emph{On the vanishing of local cohomology in characteristic {$p>0$}},
  Compos. Math. \textbf{142} (2006), no.~1, 207--221. \MR{2197409}

\bibitem[Lyu07]{Lyubeznik-AdvMath2007}
\bysame, \emph{On some local cohomology modules}, Adv. Math. \textbf{213}
  (2007), no.~2, 621--643. \MR{2332604}

\bibitem[Lyu11]{LyubeznikJPPA2011}
\bysame, \emph{A characteristic-free proof of a basic result on {$\mathcal
  D$}-modules}, J. Pure Appl. Algebra \textbf{215} (2011), no.~8, 2019--2023.
  \MR{2776441}

\bibitem[Ma14]{Ma-IMRN14}
Linquan Ma, \emph{Finiteness properties of local cohomology for {$F$}-pure
  local rings}, Int. Math. Res. Not. IMRN (2014), no.~20, 5489--5509.
  \MR{3271179}

\bibitem[Ma15]{Ma-MathAnn15}
\bysame, \emph{{$F$}-injectivity and {B}uchsbaum singularities}, Math. Ann.
  \textbf{362} (2015), no.~1-2, 25--42. \MR{3343868}

\bibitem[Mal75]{Malgrange-isolee}
B.~Malgrange, \emph{Le polyn\^ome de {B}ernstein d'une singularit\'e isol\'ee},
  Fourier integral operators and partial differential equations ({C}olloq.
  {I}nternat., {U}niv. {N}ice, {N}ice, 1974), Springer, Berlin, 1975,
  pp.~98--119. Lecture Notes in Math., Vol. 459. \MR{0419827 (54 \#7845)}

\bibitem[Mal83]{Malgrange-evan}
\bysame, \emph{Polyn\^omes de {B}ernstein-{S}ato et cohomologie \'evanescente},
  Analysis and topology on singular spaces, {II}, {III} ({L}uminy, 1981),
  Ast\'erisque, vol. 101, Soc. Math. France, Paris, 1983, pp.~243--267.
  \MR{737934 (86f:58148)}

\bibitem[Mar01]{Marley-finiteness-ass-primes}
Thomas Marley, \emph{The associated primes of local cohomology modules over
  rings of small dimension}, Manuscripta Math. \textbf{104} (2001), no.~4,
  519--525. \MR{1836111}

\bibitem[Mil21]{Milne-notes}
J.S Milne, \emph{Lectures on etale cohomology}, Version 2.21, {\tt
  https://www.jmilne.org/math/CourseNotes/lec.html }.

\bibitem[Mil68]{Milnor}
John Milnor, \emph{Singular points of complex hypersurfaces}, Annals of
  Mathematics Studies, No. 61, Princeton University Press, Princeton, N.J.,
  1968. \MR{0239612 (39 \#969)}

\bibitem[Mil80]{Milne}
James~S. Milne, \emph{\'{E}tale cohomology}, Princeton Mathematical Series,
  vol.~33, Princeton University Press, Princeton, N.J., 1980. \MR{559531}

\bibitem[MMW05]{MMW05}
Laura~Felicia Matusevich, Ezra Miller, and Uli Walther, \emph{Homological
  methods for hypergeometric families}, J. Amer. Math. Soc. \textbf{18} (2005),
  no.~4, 919--941 (electronic).

\bibitem[Moh85]{Moh-MR0784166}
T.~T. Moh, \emph{Set-theoretic complete intersections}, Proc. Amer. Math. Soc.
  \textbf{94} (1985), no.~2, 217--220. \MR{784166}

\bibitem[Mon13]{Montaner-LNM13}
Josep~\`Alvarez Montaner, \emph{Local cohomology modules supported on monomial
  ideals}, Monomial ideals, computations and applications, Lecture Notes in
  Math., vol. 2083, Springer, Heidelberg, 2013, pp.~109--178. \MR{3184122}

\bibitem[MQ18]{MaQuy-Nagoya18}
Linquan Ma and Pham~Hung Quy, \emph{Frobenius actions on local cohomology
  modules and deformation}, Nagoya Math. J. \textbf{232} (2018), 55--75.
  \MR{3866500}

\bibitem[MS97]{MehtaSrinivas-Asian97}
V.~B. Mehta and V.~Srinivas, \emph{A characterization of rational
  singularities}, Asian J. Math. \textbf{1} (1997), no.~2, 249--271.
  \MR{1491985}

\bibitem[MS12]{SchenzelMahmood}
Waqas Mahmood and Peter Schenzel, \emph{On invariants and endomorphism rings of
  certain local cohomology modules}, J. Algebra \textbf{372} (2012), 56--67.
  \MR{2990000}

\bibitem[MSTW01]{MSTW}
Mircea Musta\c{t}\u{a}, Gregory~G. Smith, Harrison Tsai, and Uli Walther,
  \emph{{$\calD$}-modules on smooth toric varieties}, J. Algebra \textbf{240}
  (2001), no.~2, 744--770. \MR{1841355}

\bibitem[MSV14]{MillerSinghVarbaro14}
Lance~Edward Miller, Anurag~K. Singh, and Matteo Varbaro, \emph{The {$F$}-pure
  threshold of a determinantal ideal}, Bull. Braz. Math. Soc. (N.S.)
  \textbf{45} (2014), no.~4, 767--775. \MR{3296192}

\bibitem[MSW21]{MaSinghWalther-Dutta}
Linquan Ma, Anurag~K. Singh, and Uli Walther, \emph{Koszul and local
  cohomology, and a question of {D}utta}, Math. Z. \textbf{298} (2021),
  no.~1-2, 697--711. \MR{4257105}

\bibitem[Mus87]{Musson-tori}
Ian~M. Musson, \emph{Rings of differential operators on invariant rings of
  tori}, Trans. Amer. Math. Soc. \textbf{303} (1987), no.~2, 805--827.
  \MR{902799}

\bibitem[MZ14]{MaZhangEulerian}
Linquan Ma and Wenliang Zhang, \emph{Eulerian graded {$\mathscr D$}-modules},
  Math. Res. Lett. \textbf{21} (2014), no.~1, 149--167. \MR{3247047}

\bibitem[Nan04]{Nang-Japan04}
Philibert Nang, \emph{{$\mathscr D$}-modules associated to the determinantal
  singularities}, Proc. Japan Acad. Ser. A Math. Sci. \textbf{80} (2004),
  no.~5, 74--78. \MR{2062805}

\bibitem[Nan08]{Nang-Adv08}
\bysame, \emph{On a class of holonomic {$\mathscr D$}-modules on {$M_n(\mathbb
  C)$} related to the action of {${\rm GL}_n(\mathbb C)\times {\rm
  GL}_n(\mathbb C)$}}, Adv. Math. \textbf{218} (2008), no.~3, 635--648.
  \MR{2414315}

\bibitem[Nan12]{Nang-JA12}
\bysame, \emph{On the classification of regular holonomic {$\mathscr
  D$}-modules on skew-symmetric matrices}, J. Algebra \textbf{356} (2012),
  115--132. \MR{2891125}

\bibitem[Nar63]{Narita-PCPS63}
Masao Narita, \emph{A note on the coefficients of {H}ilbert characteristic
  functions in semi-regular local rings}, Proc. Cambridge Philos. Soc.
  \textbf{59} (1963), 269--275. \MR{146212}

\bibitem[NnB13]{Betancourt-IJM13}
Luis N\'{u}\~{n}ez Betancourt, \emph{Local cohomology modules of polynomial or
  power series rings over rings of small dimension}, Illinois J. Math.
  \textbf{57} (2013), no.~1, 279--294. \MR{3224571}

\bibitem[NnBSW19]{BetancourtSpiroffWitt}
Luis N\'{u}\~{n}ez Betancourt, Sandra Spiroff, and Emily~E. Witt,
  \emph{Connectedness and {L}yubeznik numbers}, Int. Math. Res. Not. IMRN
  (2019), no.~13, 4233--4259. \MR{3978438}

\bibitem[NnBWZ16]{BetancourtWittZhang-survey}
Luis N\'{u}\~{n}ez Betancourt, Emily~E. Witt, and Wenliang Zhang, \emph{A
  survey on the {L}yubeznik numbers}, Mexican mathematicians abroad: recent
  contributions, Contemp. Math., vol. 657, Amer. Math. Soc., Providence, RI,
  2016, pp.~137--163. \MR{3466449}

\bibitem[Oak97]{Oaku-Duke}
Toshinori Oaku, \emph{An algorithm of computing {$b$}-functions}, Duke Math. J.
  \textbf{87} (1997), no.~1, 115--132. \MR{1440065}

\bibitem[Ogu73]{Ogus-LCDAV}
Arthur Ogus, \emph{Local cohomological dimension of algebraic varieties}, Ann.
  of Math. (2) \textbf{98} (1973), 327--365. \MR{506248}

\bibitem[OT99]{OT1}
Toshinori Oaku and Nobuki Takayama, \emph{An algorithm for de {R}ham cohomology
  groups of the complement of an affine variety via {$D$}-module computation},
  J. Pure Appl. Algebra \textbf{139} (1999), no.~1-3, 201--233, Effective
  methods in algebraic geometry (Saint-Malo, 1998). \MR{1700544}

\bibitem[OT01]{OT2}
\bysame, \emph{Algorithms for {$D$}-modules---restriction, tensor product,
  localization, and local cohomology groups}, J. Pure Appl. Algebra
  \textbf{156} (2001), no.~2-3, 267--308. \MR{1808827}

\bibitem[Pan21]{Pandey-2005.03250}
Vaibhav Pandey, \emph{Cohomological dimension of ideals defining {V}eronese
  subrings}, Proc. Amer. Math. Soc. \textbf{149} (2021), no.~4, 1387--1393.
  \MR{4242298}

\bibitem[Pel88]{Pellikaan}
Ruud Pellikaan, \emph{Projective resolutions of the quotient of two ideals},
  Nederl. Akad. Wetensch. Indag. Math. \textbf{50} (1988), no.~1, 65--84.
  \MR{934475 (89c:13015)}

\bibitem[Per20]{Perlmann-JA20}
Michael Perlman, \emph{Equivariant {$\mathcal{D}$}-modules on
  {$2\times2\times2$} hypermatrices}, J. Algebra \textbf{544} (2020), 391--416.
  \MR{4027737}

\bibitem[PQ19]{PolstraQuy-JA19}
Thomas Polstra and Pham~Hung Quy, \emph{Nilpotence of {F}robenius actions on
  local cohomology and {F}robenius closure of ideals}, J. Algebra \textbf{529}
  (2019), 196--225. \MR{3938859}

\bibitem[PS73]{PeskineSzpiro-IHES73}
C.~Peskine and L.~Szpiro, \emph{Dimension projective finie et cohomologie
  locale. {A}pplications \`a la d\'{e}monstration de conjectures de {M}.
  {A}uslander, {H}. {B}ass et {A}. {G}rothendieck}, Inst. Hautes \'{E}tudes
  Sci. Publ. Math. (1973), no.~42, 47--119. \MR{374130}

\bibitem[PS19]{PuthenpurakalSingh2019}
Tony~J. Puthenpurakal and Jyoti Singh, \emph{On derived functors of graded
  local cohomology modules}, Math. Proc. Cambridge Philos. Soc. \textbf{167}
  (2019), no.~3, 549--565. \MR{4015650}

\bibitem[PS20]{PolstraSimpsonFPurityDeformsQGorenstein}
Thomas Polstra and Austyn Simpson, \emph{{$F$}-purity deforms in
  $\mathbb{Q}$-gorenstein rings}, arXiv:2009.13444, 2020.

\bibitem[Put14]{PuthenpurakalInjectiveResolutionLC}
Tony~J. Puthenpurakal, \emph{On injective resolutions of local cohomology
  modules}, Illinois J. Math. \textbf{58} (2014), no.~3, 709--718. \MR{3395959}

\bibitem[Put15]{PuthenpurakalNagoya2015}
\bysame, \emph{De {R}ham cohomology of local cohomology modules: the graded
  case}, Nagoya Math. J. \textbf{217} (2015), 1--21. \MR{3343837}

\bibitem[Put19]{Puthenpurakal-Vietnam19}
\bysame, \emph{Bockstein cohomology of associated graded rings}, Acta Math.
  Vietnam. \textbf{44} (2019), no.~1, 285--306. \MR{3935302}

\bibitem[QS17]{QuyShimomoto-Adv17}
Pham~Hung Quy and Kazuma Shimomoto, \emph{{$F$}-injectivity and {F}robenius
  closure of ideals in {N}oetherian rings of characteristic {$p>0$}}, Adv.
  Math. \textbf{313} (2017), 127--166. \MR{3649223}

\bibitem[Rai16]{Raicu-Compositio16}
Claudiu Raicu, \emph{Characters of equivariant {$\mathcal{D}$}-modules on
  spaces of matrices}, Compos. Math. \textbf{152} (2016), no.~9, 1935--1965.
  \MR{3568944}
  
  \bibitem[RW14]{RaicuWeyman-ANT14}
Claudiu Raicu and Jerzy Weyman, \emph{Local cohomology with support in generic
  determinantal ideals}, Algebra Number Theory \textbf{8} (2014), no.~5,
  1231--1257. \MR{3263142}

\bibitem[RW16]{RaicuWeyman-JLMS16}
\bysame, \emph{Local cohomology with support in ideals of symmetric minors and
  {P}faffians}, J. Lond. Math. Soc. (2) \textbf{94} (2016), no.~3, 709--725.
  \MR{3614925}

\bibitem[RWW14]{RaicuWeymanWitt-Adv14}
Claudiu Raicu, Jerzy Weyman, and Emily~E. Witt, \emph{Local cohomology with
  support in ideals of maximal minors and sub-maximal {P}faffians}, Adv. Math.
  \textbf{250} (2014), 596--610. \MR{3122178}

\bibitem[{Rei}14]{Reich2}
Thomas {Reichelt}, \emph{{Laurent Polynomials, GKZ-hypergeometric Systems and
  Mixed Hodge Modules}}, Compositio Mathematica \textbf{(150)} (2014),
  911--941.

\bibitem[RS20]{ReiSe-Hodge}
Thomas Reichelt and Christian Sevenheck, \emph{Hypergeometric {H}odge modules},
  Algebr. Geom. \textbf{7} (2020), no.~3, 263--345.
  
  \bibitem[RSW21]{RSW-lambda}
Thomas Reichelt, Morihiko Saito, and Uli Walther, \emph{Dependence of
  {L}yubeznik numbers of cones of projective schemes on projective embeddings},
  Selecta Math. (N.S.) \textbf{27} (2021), no.~1, Paper No. 6, 22. \MR{4202748}

\bibitem[RSSW21]{RSSW}
Thomas Reichelt, Mathias Schulze, Christian Sevenheck, and Uli Walther,
  \emph{Algebraic aspects of hypergeometric differential equations}, Beitr.
  Algebra Geom. \textbf{62} (2021), no.~1, 137--203. \MR{4249859}

\bibitem[RW]{RW-weight}
Thomas Reichelt and Uli Walther, \emph{Weight filtrations on {GKZ}-systems},
  Preprint arXiv:1809.04247.

\bibitem[RWZ21]{RWZ-lambda}
Thomas Reichelt, Uli Walther, and Wenliang Zhang, \emph{On {L}yubeznik type
  invariants}, arXiv:2106.04457, Topology Appl., to appaear.

\bibitem[Rob76]{RobertsApplicationsDualizingComplexes1976}
Paul Roberts, \emph{Two applications of dualizing complexes over local rings},
  Ann. Sci. \'{E}cole Norm. Sup. (4) \textbf{9} (1976), no.~1, 103--106.
  \MR{399075}

\bibitem[RS79]{RoloffStueckrad}
Hartmut Roloff and J\"{u}rgen St\"{u}ckrad, \emph{Bemerkungen \"{u}ber
  {Z}usammenhangseigenschaften und mengentheoretische {D}arstellung projektiver
  algebraischer {M}annigfaltigkeiten}, Wiss. Beitr. Martin-Luther-Univ.
  Halle-Wittenberg M \textbf{12} (1979), 125--131, Beitr\"{a}ge zur Algebra und
  Geometrie, 8. \MR{571359}

\bibitem[Rot09]{Rotman}
Joseph~J. Rotman, \emph{An introduction to homological algebra}, second ed.,
  Universitext, Springer, New York, 2009. \MR{2455920}

\bibitem[Sai80]{KSaito-logforms}
Kyoji Saito, \emph{Theory of logarithmic differential forms and logarithmic
  vector fields}, J. Fac. Sci. Univ. Tokyo Sect. IA Math. \textbf{27} (1980),
  no.~2, 265--291. \MR{586450}

\bibitem[Sai90]{SaitoMHM}
Morihiko Saito, \emph{Mixed {H}odge modules}, Publ. Res. Inst. Math. Sci.
  \textbf{26} (1990), no.~2, 221--333.

\bibitem[Sai09]{Saito-ASPM09}
\bysame, \emph{On {$b$}-function, spectrum and multiplier ideals}, Algebraic
  analysis and around, Adv. Stud. Pure Math., vol.~54, Math. Soc. Japan, Tokyo,
  2009, pp.~355--379. \MR{2499561}

\bibitem[Sai16a]{Saito-1609.04801}
\bysame, \emph{{B}ernstein--{S}ato polynomials for projective hypersurfaces
  with weighted homogeneous isolated singularities}, Preprint arXiv:1609.04801,
  2016.

\bibitem[Sai16b]{Saito-BS-arr}
\bysame, \emph{Bernstein-{S}ato polynomials of hyperplane arrangements},
  Selecta Math. (N.S.) \textbf{22} (2016), no.~4, 2017--2057. \MR{3573952}

\bibitem[Sai17]{Saito-1703.05741}
\bysame, \emph{Roots of {B}ernstein--{S}ato polynomials of certain homogeneous
  polynomials with two-dimensional singular loci}, Pure Appl. Math. Q. \textbf{16} (2020), no. 4, 1219--1280. \MR{4180246}
  
  
\bibitem[Sch82a]{SchenzelCohomologicalAnnihiltors1982}
Peter Schenzel, \emph{Cohomological annihilators}, Math. Proc. Cambridge
  Philos. Soc. \textbf{91} (1982), no.~3, 345--350. \MR{654081}

\bibitem[Sch82b]{Schenzel-Dualisierende}
\bysame, \emph{Dualisierende {K}omplexe in der lokalen {A}lgebra und
  {B}uchsbaum-{R}inge}, Lecture Notes in Mathematics, vol. 907,
  Springer-Verlag, Berlin-New York, 1982, With an English summary. \MR{654151}

\bibitem[Sch98]{SchenzelUseofLocalCohomology}
\bysame, \emph{On the use of local cohomology in algebra and geometry}, Six
  lectures on commutative algebra ({B}ellaterra, 1996), Progr. Math., vol. 166,
  Birkh\"{a}user, Basel, 1998, pp.~241--292. \MR{1648667}

\bibitem[Sch07]{Schwede-Compositio07}
Karl Schwede, \emph{A simple characterization of {D}u {B}ois singularities},
  Compos. Math. \textbf{143} (2007), no.~4, 813--828. \MR{2339829}

\bibitem[Sch09a]{Schenzel-OnConnectedness}
Peter Schenzel, \emph{On connectedness and indecomposibility of local
  cohomology modules}, Manuscripta Math. \textbf{128} (2009), no.~3, 315--327.
  \MR{2481047}

\bibitem[Sch09b]{SchenzelPAMS2009}
\bysame, \emph{On endomorphism rings and dimensions of local cohomology
  modules}, Proc. Amer. Math. Soc. \textbf{137} (2009), no.~4, 1315--1322.
  \MR{2465654}

\bibitem[Sch09c]{Schwede-AJM09}
Karl Schwede, \emph{{$F$}-injective singularities are {D}u {B}ois}, Amer. J.
  Math. \textbf{131} (2009), no.~2, 445--473. \MR{2503989}

\bibitem[Sch10]{SchenzelArchMath2010}
Peter Schenzel, \emph{Matlis duals of local cohomology modules and their
  endomorphism rings}, Arch. Math. (Basel) \textbf{95} (2010), no.~2, 115--123.
  \MR{2674247}

\bibitem[Sch11a]{SchenzelJA2011}
\bysame, \emph{On the structure of the endomorphism ring of a certain local
  cohomology module}, J. Algebra \textbf{344} (2011), 229--245. \MR{2831938}

\bibitem[Sch11b]{Schwede-TAMS11}
Karl Schwede, \emph{Test ideals in non-{$\mathbb{Q}$}-{G}orenstein rings},
  Trans. Amer. Math. Soc. \textbf{363} (2011), no.~11, 5925--5941. \MR{2817415}

\bibitem[Ser55]{Serre-FAC}
Jean-Pierre Serre, \emph{Faisceaux alg\'{e}briques coh\'{e}rents}, Ann. of
  Math. (2) \textbf{61} (1955), 197--278. \MR{68874}

\bibitem[Ser56]{Serre-GAGA}
\bysame, \emph{G\'{e}om\'{e}trie alg\'{e}brique et g\'{e}om\'{e}trie
  analytique}, Ann. Inst. Fourier (Grenoble) \textbf{6} (1955/56), 1--42.
  \MR{82175}

\bibitem[Sha69]{Sharp-Cousin}
Rodney~Y. Sharp, \emph{The {C}ousin complex for a module over a commutative
  {N}oetherian ring}, Math. Z. \textbf{112} (1969), 340--356. \MR{263800}

\bibitem[Sha77]{Sharp-Cousin+lc}
\bysame, \emph{Local cohomology and the {C}ousin complex for a commutative
  {N}oetherian ring}, Math. Z. \textbf{153} (1977), no.~1, 19--22. \MR{442062}

\bibitem[Sha07a]{Sharp-TAMS07}
\bysame, \emph{Graded annihilators of modules over the {F}robenius skew
  polynomial ring, and tight closure}, Trans. Amer. Math. Soc. \textbf{359}
  (2007), no.~9, 4237--4258. \MR{2309183}

\bibitem[Sha07b]{Sharp-PAMS07}
\bysame, \emph{On the {H}artshorne-{S}peiser-{L}yubeznik theorem about
  {A}rtinian modules with a {F}robenius action}, Proc. Amer. Math. Soc.
  \textbf{135} (2007), no.~3, 665--670. \MR{2262861}

\bibitem[Sin99a]{Singh-JPAA99}
Anurag~K. Singh, \emph{Deformation of {$F$}-purity and {$F$}-regularity}, J.
  Pure Appl. Algebra \textbf{140} (1999), no.~2, 137--148. \MR{1693967}

\bibitem[Sin99b]{Anurag-AJM99}
\bysame, \emph{{$F$}-regularity does not deform}, Amer. J. Math. \textbf{121}
  (1999), no.~4, 919--929. \MR{1704481}

\bibitem[Sin00]{Singh-p-torsion-elements}
\bysame, \emph{{$p$}-torsion elements in local cohomology modules}, Math. Res.
  Lett. \textbf{7} (2000), no.~2-3, 165--176. \MR{1764314}

\bibitem[SKK73]{KashiwaraKawaiSato1973}
Mikio Sato, Takahiro Kawai, and Masaki Kashiwara, \emph{Microfunctions and
  pseudo-differential equations}, Hyperfunctions and pseudo-differential
  equations ({P}roc. {C}onf., {K}atata, 1971; dedicated to the memory of
  {A}ndr\'{e} {M}artineau), 1973, pp.~265--529. Lecture Notes in Math., Vol.
  287. \MR{0420735}

\bibitem[Smi97a]{Smith-AJM97}
Karen~E. Smith, \emph{{$F$}-rational rings have rational singularities}, Amer.
  J. Math. \textbf{119} (1997), no.~1, 159--180. \MR{1428062}

\bibitem[Smi97b]{Smith-FujitaConjecture-JAG1997}
\bysame, \emph{Fujita's freeness conjecture in terms of local cohomology}, J.
  Algebraic Geom. \textbf{6} (1997), no.~3, 417--429. \MR{1487221}

\bibitem[Spe78]{Speiser-TAMS78}
Robert Speiser, \emph{Projective varieties of low codimension in characteristic
  {$p>0$}}, Trans. Amer. Math. Soc. \textbf{240} (1978), 329--343. \MR{491703}

\bibitem[SS04]{SinghSwanson-ass-primes}
Anurag~K. Singh and Irena Swanson, \emph{Associated primes of local cohomology
  modules and of {F}robenius powers}, Int. Math. Res. Not. (2004), no.~33,
  1703--1733. \MR{2058025}

\bibitem[SS18]{SchenzelSimon-NotNoetherian}
Peter Schenzel and Anne-Marie Simon, \emph{Completion, \v{C}ech and local
  homology and cohomology}, Springer Monographs in Mathematics, Springer, Cham,
  2018, Interactions between them. \MR{3838396}

\bibitem[SST00]{SST00}
Mutsumi Saito, Bernd Sturmfels, and Nobuki Takayama, \emph{Gr\"{o}bner
  deformations of hypergeometric differential equations}, Algorithms and
  Computation in Mathematics, vol.~6, Springer-Verlag, Berlin, 2000.
  \MR{1734566}

\bibitem[ST98]{SturmfelsTakayama98}
Bernd Sturmfels and Nobuki Takayama, \emph{Gr\"{o}bner bases and hypergeometric
  functions}, Gr\"{o}bner bases and applications ({L}inz, 1998), London Math.
  Soc. Lecture Note Ser., vol. 251, Cambridge Univ. Press, Cambridge, 1998,
  pp.~246--258. \MR{1708882}

\bibitem[ST17]{SrinivasTakagi-Adv17}
Vasudevan Srinivas and Shunsuke Takagi, \emph{Nilpotence of {F}robenius action
  and the {H}odge filtration on local cohomology}, Adv. Math. \textbf{305}
  (2017), 456--478. \MR{3570141}

\bibitem[Ste87]{Steenbrink-survey}
J.~Steenbrink, \emph{Mixed {H}odge structures and singularities: a survey},
  G\'eom\'etrie alg\'ebrique et applications, {III} ({L}a {R}\'abida, 1984),
  Travaux en Cours, vol.~24, Hermann, Paris, 1987, pp.~99--123. \MR{907936
  (89c:14004)}

\bibitem[Ste19a]{Steiner-JA}
Avi Steiner, \emph{{$A$}-hypergeometric modules and {G}auss-{M}anin systems},
  J. Algebra \textbf{524} (2019), 124--159. \MR{3904304}

\bibitem[Ste19b]{Steiner-JPAA}
\bysame, \emph{Dualizing, projecting, and restricting {GKZ} systems}, J. Pure
  Appl. Algebra \textbf{223} (2019), no.~12, 5215--5231. \MR{3975063}

\bibitem[Sti98]{Stienstra}
Jan Stienstra, \emph{Resonant hypergeometric systems and mirror symmetry},
  Integrable systems and algebraic geometry ({K}obe/{K}yoto, 1997), World Sci.
  Publ., River Edge, NJ, 1998, pp.~412--452.

\bibitem[SV79]{SchmittVogel}
Thomas Schmitt and Wolfgang Vogel, \emph{Note on set-theoretic intersections of
  subvarieties of projective space}, Math. Ann. \textbf{245} (1979), no.~3,
  247--253. \MR{553343}

\bibitem[SV86]{StueckradVogel}
J\"{u}rgen St\"{u}ckrad and Wolfgang Vogel, \emph{Buchsbaum rings and
  applications}, Springer-Verlag, Berlin, 1986, An interaction between algebra,
  geometry and topology. \MR{881220}

\bibitem[SW05]{SinghWalther-Segre}
Anurag~K. Singh and Uli Walther, \emph{On the arithmetic rank of certain
  {S}egre products}, Commutative algebra and algebraic geometry, Contemp.
  Math., vol. 390, Amer. Math. Soc., Providence, RI, 2005, pp.~147--155.
  \MR{2187332}

\bibitem[SW07]{SinghWalther-pure}
\bysame, \emph{Local cohomology and pure morphisms}, Illinois J. Math.
  \textbf{51} (2007), no.~1, 287--298. \MR{2346198}

\bibitem[SW08a]{SchulzeWalther-slopes}
Mathias Schulze and Uli Walther, \emph{Irregularity of hypergeometric systems
  via slopes along coordinate subspaces}, Duke Math. J. \textbf{142} (2008),
  no.~3, 465--509.

\bibitem[SW08b]{SinghWalther-TAMS08}
Anurag~K. Singh and Uli Walther, \emph{A connectedness result in positive
  characteristic}, Trans. Amer. Math. Soc. \textbf{360} (2008), no.~6,
  3107--3119. \MR{2379789}

\bibitem[SW09]{SchulzeWalther-ekdi}
Mathias Schulze and Uli Walther, \emph{Hypergeometric $\mathcal{D}$-modules and
  twisted {G}au\ss-{M}anin systems}, J. Algebra \textbf{322} (2009), no.~9,
  3392--3409.

\bibitem[SW11]{SinghWalther-Bockstein}
Anurag~K. Singh and Uli Walther, \emph{Bockstein homomorphisms in local
  cohomology}, J. Reine Angew. Math. \textbf{655} (2011), 147--164.
  \MR{2806109}

\bibitem[SW20]{SinghWalther-CiA2020}
\bysame, \emph{On a conjecture of {L}ynch}, Comm. Algebra \textbf{48} (2020),
  no.~6, 2681--2682. \MR{4107600}

\bibitem[Swi15]{Switala}
Nicholas Switala, \emph{Lyubeznik numbers for nonsingular projective
  varieties}, Bull. Lond. Math. Soc. \textbf{47} (2015), no.~1, 1--6.
  \MR{3312957}

\bibitem[Swi17a]{SwitalaCompos2017}
\bysame, \emph{On the de {R}ham homology and cohomology of a complete local
  ring in equicharacteristic zero}, Compos. Math. \textbf{153} (2017), no.~10,
  2075--2146. \MR{3705285}

\bibitem[Swi17b]{Switala-dR-complete}
\bysame, \emph{On the de {R}ham homology and cohomology of a complete local
  ring in equicharacteristic zero}, Compos. Math. \textbf{153} (2017), no.~10,
  2075--2146. \MR{3705285}

\bibitem[SZ18]{SwitalaZhangAdvMath2018}
Nicholas Switala and Wenliang Zhang, \emph{Duality and de {R}ham cohomology for
  graded {$\mathscr{D}$}-modules}, Adv. Math. \textbf{340} (2018), 1141--1165.
  \MR{3886190}

\bibitem[SZ19]{SwitalaZhangInjectiveDim}
\bysame, \emph{A dichotomy for the injective dimension of {$F$}-finite
  {$F$}-modules and holonomic {$D$}-modules}, Comm. Algebra \textbf{47} (2019),
  no.~6, 2525--2539. \MR{3957114}

\bibitem[Tor09]{Torrelli-intHom}
Tristan Torrelli, \emph{Intersection homology {$\mathcal D$}-module and
  {B}ernstein polynomials associated with a complete intersection}, Publ. Res.
  Inst. Math. Sci. \textbf{45} (2009), no.~2, 645--660. \MR{2510514}

\bibitem[Tra06]{Traves-orbi}
William~N. Traves, \emph{Differential operators on orbifolds}, J. Symbolic
  Comput. \textbf{41} (2006), no.~12, 1295--1308. \MR{2271326}

\bibitem[Tri97]{TrippDiffOperStanleyReisnerRings}
J.~R. Tripp, \emph{Differential operators on {S}tanley-{R}eisner rings}, Trans.
  Amer. Math. Soc. \textbf{349} (1997), no.~6, 2507--2523. \MR{1376559}

\bibitem[TT08]{TakagiTakahashi-MRL08}
Shunsuke Takagi and Ryo Takahashi, \emph{{$D$}-modules over rings with finite
  {$F$}-representation type}, Math. Res. Lett. \textbf{15} (2008), no.~3,
  563--581. \MR{2407232}

\bibitem[TW18]{TakagiWatanabe-Sugaku18}
Shunsuke Takagi and Kei-Ichi Watanabe, \emph{{$F$}-singularities: applications
  of characteristic {$p$} methods to singularity theory [translation of
  {MR}3135334]}, Sugaku Expositions \textbf{31} (2018), no.~1, 1--42.
  \MR{3784697}

\bibitem[Var12]{Varbaro-TAMS12}
Matteo Varbaro, \emph{On the arithmetical rank of certain {S}egre embeddings},
  Trans. Amer. Math. Soc. \textbf{364} (2012), no.~10, 5091--5109. \MR{2931323}

\bibitem[Var13]{Varbaro-Compositio13}
\bysame, \emph{Cohomological and projective dimensions}, Compos. Math.
  \textbf{149} (2013), no.~7, 1203--1210. \MR{3078644}

\bibitem[Var19]{Varbaro-CiA2019}
\bysame, \emph{Connectivity of hyperplane sections of domains}, Comm. Algebra
  \textbf{47} (2019), no.~6, 2540--2547. \MR{3957115}

\bibitem[vdB]{vanDenBergh-notes}
Michel van~den Bergh, \emph{Some generalities on quasi-coherent
  {$\mathcal{O}_X$} and {$\mathcal{D}_X$}-modules, {P}reprint},
  \url{https://hardy.uhasselt.be/personal/vdbergh/Publications/Geq.ps}.

\bibitem[vdB99]{vanDenBergh-Adv99}
\bysame, \emph{Local cohomology of modules of covariants}, Adv. Math.
  \textbf{144} (1999), no.~2, 161--220. \MR{1695237}

\bibitem[vdE85]{vdEcoker1985}
A.~van~den Essen, \emph{The cokernel of the operator $\frac{\partial}{\partial
  x_n}$ acting on a $\mathscr{D}_n$-module, {II}}, Compos. Math. \textbf{56}
  (1985), no.~2, 259--269.

\bibitem[Vog71]{Vogel-Monatsberichte71}
W.~Vogel, \emph{Eine {B}emerkung \"{u}ber die {A}nzahl von {H}yperfl\"{a}chen
  zur {D}arstellung algebraischer {V}ariet\"{a}ten}, Monatsb. Deutsch. Akad.
  Wiss. Berlin \textbf{13} (1971), 629--633. \MR{325630}

\bibitem[Vra00]{Vraciu-JA2000}
Adela Vraciu, \emph{Local cohomology of {F}robenius images over graded affine
  algebras}, J. Algebra \textbf{228} (2000), no.~1, 347--356. \MR{1760968}

\bibitem[vSW15]{vanStratenWarmt}
Duco van Straten and Thorsten Warmt, \emph{Gorenstein-duality for
  one-dimensional almost complete intersections---with an application to
  non-isolated real singularities}, Math. Proc. Cambridge Philos. Soc.
  \textbf{158} (2015), no.~2, 249--268. \MR{3310244}

\bibitem[Wal99]{W-lcD}
Uli Walther, \emph{Algorithmic computation of local cohomology modules and the
  local cohomological dimension of algebraic varieties}, J. Pure Appl. Algebra
  \textbf{139} (1999), no.~1-3, 303--321, Effective methods in algebraic
  geometry (Saint-Malo, 1998). \MR{1700548}

\bibitem[Wal00]{W-algdR00}
\bysame, \emph{Algorithmic computation of de {R}ham cohomology of complements
  of complex affine varieties}, J. Symbolic Comput. \textbf{29} (2000),
  no.~4-5, 795--839, Symbolic computation in algebra, analysis, and geometry
  (Berkeley, CA, 1998). \MR{1769667}

\bibitem[Wal01a]{W-alg-det01}
\bysame, \emph{Algorithmic determination of the rational cohomology of complex
  varieties via differential forms}, Symbolic computation: solving equations in
  algebra, geometry, and engineering ({S}outh {H}adley, {MA}, 2000), Contemp.
  Math., vol. 286, Amer. Math. Soc., Providence, RI, 2001, pp.~185--206.
  \MR{1874280}

\bibitem[Wal01b]{W-lambda}
\bysame, \emph{On the {L}yubeznik numbers of a local ring}, Proc. Amer. Math.
  Soc. \textbf{129} (2001), no.~6, 1631--1634. \MR{1814090}

\bibitem[Wal02]{W-M2book}
\bysame, \emph{{$D$}-modules and cohomology of varieties}, Computations in
  algebraic geometry with {M}acaulay 2, Algorithms Comput. Math., vol.~8,
  Springer, Berlin, 2002, pp.~281--323. \MR{1949555}

\bibitem[Wal05]{Walther-Bernstein}
\bysame, \emph{Bernstein-{S}ato polynomial versus cohomology of the {M}ilnor
  fiber for generic hyperplane arrangements}, Compos. Math. \textbf{141}
  (2005), no.~1, 121--145. \MR{2099772 (2005k:32030)}

\bibitem[Wal15]{Walther-survey}
\bysame, \emph{Survey on the {$D$}-module {$f^s$}}, Commutative algebra and
  noncommutative algebraic geometry. {V}ol. {I}, Math. Sci. Res. Inst. Publ.,
  vol.~67, Cambridge Univ. Press, New York, 2015, With an appendix by Anton
  Leykin, pp.~391--430. \MR{3525478}

\bibitem[Wal17]{Walther-zieg}
\bysame, \emph{The {J}acobian module, the {M}ilnor fiber, and the {$D$}-module
  generated by {$f^s$}}, Invent. Math. \textbf{207} (2017), no.~3, 1239--1287.
  \MR{3608290}

\bibitem[Wan20]{Wang-lambda}
Botong Wang, \emph{Lyubeznik numbers of irreducible projective varieties depend
  on the embedding}, Proc. Amer. Math. Soc. \textbf{148} (2020), no.~5,
  2091--2096. \MR{4078092}

\bibitem[Yan99]{Yan-JA99}
Zhao Yan, \emph{Minimal resultant systems}, J. Algebra \textbf{216} (1999),
  no.~1, 105--123. \MR{1694582}

\bibitem[Yan00]{Yan-JPAA00}
\bysame, \emph{An \'{e}tale analog of the {G}oresky-{M}ac{P}herson formula for
  subspace arrangements}, J. Pure Appl. Algebra \textbf{146} (2000), no.~3,
  305--318. \MR{1742346}

\bibitem[Yan03]{Yanagawa03}
Kohji Yanagawa, \emph{Stanley-{R}eisner rings, sheaves, and
  {P}oincar\'{e}-{V}erdier duality}, Math. Res. Lett. \textbf{10} (2003),
  no.~5-6, 635--650. \MR{2024721}

\bibitem[Zha07]{Zhang-highest}
Wenliang Zhang, \emph{On the highest {L}yubeznik number of a local ring},
  Compos. Math. \textbf{143} (2007), no.~1, 82--88. \MR{2295196}

\bibitem[Zha11]{Zhang-lambda}
\bysame, \emph{Lyubeznik numbers of projective schemes}, Adv. Math.
  \textbf{228} (2011), no.~1, 575--616. \MR{2822240}

\bibitem[Zha12]{YiZhangGradedFModules}
Yi~Zhang, \emph{Graded {$F$}-modules and local cohomology}, Bull. Lond. Math.
  Soc. \textbf{44} (2012), no.~4, 758--762. \MR{2967243}

\bibitem[Zha15]{ZhangRegularity-2015}
Wenliang Zhang, \emph{A note on the growth of regularity with respect to
  {F}robenius}, arXiv:1512.00049, 2015.

\bibitem[Zha17]{ZhangInjectiveDim}
\bysame, \emph{On injective dimension of {$F$}-finite {$F$}-modules and
  holonomic {$D$}-modules}, Bull. Lond. Math. Soc. \textbf{49} (2017), no.~4,
  593--603. \MR{3725482}

\bibitem[Zha21a]{Zhang2020}
\bysame, \emph{On asymptotic socle degrees of local cohomology modules}, J.
  Pure Appl. Algebra \textbf{225} (2021), no.~12, 106789. \MR{4260035}

\bibitem[Zha21b]{ZhangVanishingLCMixedChar}
\bysame, \emph{The second vanishing theorem for local cohomology modules},
  arXiv:2102.12545, submitted.

\bibitem[Zho98]{ZhouJA1998}
Caijun Zhou, \emph{Higher derivations and local cohomology modules}, J. Algebra
  \textbf{201} (1998), no.~2, 363--372. \MR{1612378}

\bibitem[Zho06]{ZhouUniformAnn-2006}
\bysame, \emph{Uniform annihilators of local cohomology}, J. Algebra
  \textbf{305} (2006), no.~1, 585--602. \MR{2264146}

\bibitem[Zho07]{ZhouUniformAnn-2007}
\bysame, \emph{Uniform annihilators of local cohomology of excellent rings}, J.
  Algebra \textbf{315} (2007), no.~1, 286--300. \MR{2344347}

\end{thebibliography}
\end{document}